\documentclass[12pt]{article}
\usepackage{amssymb,latexsym}
\title{Diophantine tori and spectral asymptotics for non-selfadjoint operators}
\author{Michael Hitrik\\Department of Mathematics \\University of California \\ Los Angeles
\\CA 90095-1555, USA\\hitrik@math.ucla.edu \and
Johannes Sj\"ostrand\\Centre de Math\'ematiques\\Laurent
Schwartz\\Ecole Polytechnique\\FR--91128 Palaiseau\\France \\and
UMR 7640 CNRS
\\johannes@math.polytechnique.fr \and San V\~u Ng\d{o}c \\ Institut Fourier \\ UMR CNRS-UJF 5582 \\ BP 74, 38402
\\Saint-Martin d'H\`eres \\ France\\san.vu-ngoc@ujf-grenoble.fr}
\date{}

\def\wrtext#1{\relax\ifmmode{\leavevmode\hbox{#1}}\else{#1}\fi}
\def\abs#1{\left|#1\right|}
\def\begeq{\begin{equation}}
\def\endeq{\end{equation}}
\def\Remark{\vskip 2mm \noindent {\em Remark}}

\def\ekv#1#2{\begeq\label{#1}#2\endeq}

\def\ably{arbitrarily}
\def\an{analytic}
\def\asy{asymptotic}
\def\bdd{bounded}

\def\ev{eigenvalue}

\def\fu{function}
\def\fy{family}

\def\hol{holomorphic}

\def\mfld{manifold}

\def\neigh{neighborhood}
\def\nondeg{non-degenerate}

\def\sa{selfadjoint}

\def\sufly{sufficiently}

\def\ufly{uniformly}

\def\wrt{with respect to}

\def\Re{{\rm Re\,}}
\def\Im{{\rm Im\,}}

\newcommand{\eps}{\epsilon}
\def\part#1{\frac{\partial}{\partial #1}}

\def\norm#1{||\,#1\,||}

\newcommand{\real}{\mbox{\bf R}}
\newcommand{\comp}{\mbox{\bf C}}
\newcommand{\z}{\mbox{\bf Z}}
\newcommand{\nat}{\mbox{\bf N}}

\renewcommand{\Re}{\mbox{\rm Re\,}}
\renewcommand{\Im}{\mbox{\rm Im\,}}

\renewcommand{\exp}{\mbox{\rm exp\,}}
\newcommand{\supp}{\mbox{\rm supp}}

\newtheorem{dref}{Definition}[section]
\newtheorem{lemma}[dref]{Lemma}
\newtheorem{theo}[dref]{Theorem}
\newtheorem{prop}[dref]{Proposition}

\newenvironment{proof}{\vspace{.3cm}\noindent{{\em Proof:}}}{\hfill$\Box$}
\begin{document}
\maketitle

\begin{abstract}
We study spectral asymptotics for small non-selfadjoint
perturbations of selfadjoint $h$-pseudodifferential operators in
dimension 2, assuming that the classical flow of the unperturbed
part possesses several invariant Lagrangian tori enjoying a
Diophantine property. We get complete asymptotic expansions for
all eigenvalues in certain rectangles in the complex
plane in two different cases: in the first case, we assume that
the strength $\eps$ of the perturbation is ${\cal O}(h^{\delta})$
for some $\delta>0$ and is bounded from below by a fixed positive
power of $h$. In the second case, $\eps$ is assumed to be
sufficiently small but independent of $h$, and we describe the
eigenvalues completely in a fixed $h$-independent domain in the
complex spectral plane.
\end{abstract}



\section{Introduction and statement of main results}
\label{section1} \setcounter{equation}{0} \setcounter{dref}{0}
Recently there has been a large number of new developments for
non-selfadjoint problems. These include semiclassical spectral
asymptotics for non-selfadjoint ope\-ra\-tors in low
dimensions~\cite{Hi2004},~\cite{Nedelec},~\cite{MeSj},~\cite{Sj04},~\cite{Shk},
resolvent estimates and pseudospectral
behavior~\cite{DSZ},~\cite{Davies},~\cite{BU}, spectral
instability questions~\cite{Ha},~\cite{Re}, and evolution problems
and decay to equilibrium for the Fokker-Planck
operator~\cite{HeSjSt}. The purpose of this work is to continue a
line of development initiated in~\cite{MeSj}, which opened up the
possibility of carrying out a spectral analysis for
non-selfadjoint operators in dimension two, that is as precise as
the corresponding analysis for selfadjoint problems in dimension
one. In~\cite{MeSj}, it was established that for a wide and stable
class of non-selfadjoint operators in dimension two, it is
possible to describe all eigenvalues individually in a fixed
domain in the complex plane, by means of a Bohr-Sommerfeld
quantization condition. The underlying reason for this result is a
version of the KAM theorem without small divisors, in a complex
domain.

The work~\cite{MeSj} has been continued in a series of papers~\cite{Sj2004},~\cite{HiSj1},
~\cite{HiSj2},~\cite{HiSj3}, all of them done in the context of small non-selfadjoint perturbations
of selfadjoint operators, with the important additional assumption
that the classical flow of the leading symbol of the unperturbed part
should be periodic in some energy shell. While the case of a periodic
classical flow is very special indeed, in the aforementioned works, we have
already given some applications of the general results to spectral
asymptotics for damped wave equations on analytic Zoll
surfaces~\cite{Sj00},~\cite{Hi02}, while barrier top resonances for
semiclassical Schr\"odinger operators have been treated in~\cite{KaKer}.

Now a classical Hamiltonian with a periodic flow can be
naturally viewed as a degenerate case of a completely integrable
symbol, and an even more general and much more interesting
dynamical situation occurs when considering a symbol that is
merely close to a completely integrable one. Continuing our
previous works, in this case it seems to be of interest to study
the spectrum of non-selfadjoint operators that are small
perturbations of a selfadjoint operator, whose classical
flow is close to a completely integrable one.  The present work is
the first one where we begin to study this problem, and when doing
so, as our starting point, we shall take a general assumption that
the real energy surface of the unperturbed leading symbol contains
several flow invariant Lagrangian tori satisfying a Diophantine
condition. According to a classical theorem of
Kolmogorov~\cite{BGGS}, the existence of such invariant tori is
guaranteed when the unperturbed symbol in question is close to a
completely integrable non-degenerate one.

We shall begin by describing the general assumptions on our operators,
which will be the same as in~\cite{HiSj1},~\cite{HiSj2},
and~\cite{HiSj3}. Let $M$ denote ${\bf R}^2$ or a compact
real-analytic \mfld{} of dimension 2. We shall let $\widetilde{M}$
stand for a complexification of $M$, so that $M=\comp^2$ in the
Euclidean case, and in the manifold case, $\widetilde{M}$ is a
Grauert tube of $M$.

\par When $M={\bf R}^2$, let
\begeq
\label{0.1}
P_\epsilon =P(x,hD_x,\epsilon ;h)
\endeq
be the Weyl quantization on ${\bf R}^2$ of a symbol $P(x,\xi ,\epsilon
;h)$ depending smoothly on $\epsilon \in{\rm neigh\,}(0,{\bf R})$ with
values in the space of \hol{} \fu{}s of $(x,\xi )$ in a tubular
\neigh{} of ${\bf R}^4$ in ${\bf C}^4$, with
\begeq
\label{0.2}
\vert P(x,\xi ,\epsilon ;h)\vert \le {\cal O}(1)m(\Re (x,\xi ))
\endeq
there. Here $m$ is assumed to be an order \fu{} on ${\bf R}^4$, in the
sense that $m>0$ and
\begeq
\label{0.3}
m(X)\le C_0\langle X-Y\rangle ^{N_0}m(Y),\ X,Y\in{\bf R}^4,\quad
C_0,\,\, N_0>0.
\endeq
We also assume that
\begeq
\label{0.4}
m\ge 1.
\endeq
We further assume that
\begeq\label{0.5}
P(x,\xi ,\epsilon ;h)\sim \sum_{j=0}^\infty  p_{j,\epsilon }(x,\xi)h^j,\ h\to 0, \endeq
in the space of such \fu{}s. We make the ellipticity assumption
\begeq\label{0.6}
\vert p_{0,\epsilon }(x,\xi )\vert \ge {1\over C}m(\Re (x,\xi )),\
\vert (x,\xi )\vert \ge C, \endeq
for some $C>0$.

\par When $M$ is a compact \mfld{}, we let $P_{\eps}$ be a
differential operator on $M$, such that for every choice of local
coordinates, centered at some point of $M$, it takes the form
\begeq\label{0.7}
P_\epsilon =\sum_{\vert \alpha \vert \le m}a_{\alpha ,\epsilon
}(x;h)(hD_x)^\alpha ,\endeq
where $a_{\alpha ,\epsilon }(x;h)$ is a smooth
\fu{} of $\epsilon $ with values in the space of \bdd{} \hol{} \fu{}s
in a complex \neigh{} of $x=0$. We further assume that
\begeq\label{0.7.5}
a_{\alpha ,\epsilon }(x;h)\sim \sum_{j=0}^\infty  a_{\alpha ,\epsilon
,j}(x)h^j,\ h\to 0, \endeq
in the space of such \fu{}s. The semi-classical principal symbol
$p_{0,\eps}$, defined on $T^*M$, takes the form
\begeq\label{0.8}
p_{0,\epsilon }(x,\xi )=\sum a_{\alpha ,\epsilon ,0}(x)\xi ^\alpha ,
\endeq
if $(x,\xi)$ are canonical coordinates on $T^*M$, and we make the ellipticity assumption
\begeq\label{0.9}
\vert p_{0,\epsilon}(x,\xi )\vert \ge {1\over C}\langle \xi \rangle ^m,\ (x,\xi
)\in T^*M,\,\vert \xi \vert \ge C,\endeq
for some large $C>0$. (Here we assume that $M$
has been equipped with some Riemannian metric, so that $\vert \xi \vert
$ and $\langle \xi \rangle =(1+\vert \xi \vert ^2)^{1/2}$ are
well-defined.)

\par Sometimes, we write $p_\epsilon $ for $p_{0,\epsilon }$ and
simply $p$ for $p_{0,0}$. Assume \begeq\label{0.10} P_{\epsilon
=0} \hbox{ is formally \sa{}.}
\endeq
In the case when $M$ is compact, we let the under\-lying Hil\-bert
spa\-ce be $L^2(M,\mu (dx))$ for some po\-si\-ti\-ve real-ana\-ly\-tic den\-si\-ty
$\mu(dx)$ on $M$.

\par Under these assumptions, $P_\epsilon $ will have discrete spectrum
in some fixed \neigh{} of $0\in{\bf C}$, when $h>0,\epsilon \ge 0$
are \sufly{} small, and the spectrum in this region will be contained
in a band $\vert \Im z\vert \le {\cal O}(\epsilon )$.

\par Assume for simplicity that (with $p=p_{\epsilon =0}$)
\begeq\label{0.11}
p^{-1}(0)\cap T^*M\hbox{ is connected,}\endeq
and let us also assume that the energy level $E=0$ is non-critical, so
that $dp\neq 0$ along $p^{-1}(0)\cap T^*M$.

Let $H_p=p'_\xi \cdot {\partial \over \partial x}-p'_x\cdot
{\partial \over \partial \xi }$ be the Hamilton field of $p$. We
introduce the following hypothesis, assumed to hold throughout
this work:
\begin{eqnarray}\label{0.12}
& & \hbox{The set}\, p^{-1}(0)\cap T^*M\, \hbox{contains finitely
many analytic}\, H_p\hbox{-invariant} \\ \nonumber & &
\hbox{Lagrangian tori }\, \Lambda_j,\,\, 1\leq j\leq L,\,
\hbox{such that each}\, \Lambda_j\,\hbox{carries real analytic} \\
\nonumber & & \hbox{coordinates}\,\, x_1,\,x_2\,
\hbox{identifying}\,\Lambda_j\,\hbox{with}\, {\bf T}^2\, \hbox{so
that along} \,\,\Lambda_j, \hbox{we have,}
\end{eqnarray}
\begeq
H_p=a_1\partial_{x_1}+a_2\partial_{x_2},
\endeq
where $a_1$, $a_2\in \real$ satisfy the Diophantine condition,
\begeq \label{0.13} \abs{a\cdot k}\geq \frac{1}{C_0\abs{k}^{N_0}},
\quad 0\neq k\in \z^2,
\endeq
for some fixed $C_0$, $N_0>0$. Here ${\bf T}^2=\real^2/2\pi \z^2$
is the standard 2-torus.

We write out the first few terms in a Taylor expansion of
$p_{\eps}$ in a \neigh{} of $p^{-1}(0)\cap T^*M$,
\begeq
\label{0.14} p_{\eps}=p+i\eps q+{\cal O}(\eps^2).
\endeq

When $0\leq K\in C^{\infty}_0(\real)$ is such that $\int
K(t)\,dt=1$ and $T>0$, we introduce a ``smoothed out'' flow
average of $q$, \begeq \label{0.15}
\langle{q}\rangle_{T,K,p}=\langle{q}\rangle_{T,K}=\int
K_T(-t)q\circ \exp(tH_p)\,dt,\quad
K_T(t)=\frac{1}{T}K\left(\frac{t}{T}\right),
\endeq
defined near $p^{-1}(0)\cap T^*M$. The standard flow average
corresponds to taking $K=1_{[-1,0]}$, and we shall then write
$\langle{q}\rangle_{T,K}=\langle{q}\rangle_T$.

Let $G_T$ be an analytic function, defined in a \neigh{} of $p^{-1}(0)\cap T^*M$, such that
$$
H_p G_T=q-\langle{q}\rangle_{T,K}.
$$
This is a convolution equation along the $H_p$--trajectories, and
as in~\cite{Sj00} and~\cite{HiSj2}, we solve it by setting \begeq
\label{0.16} G_T=\int T J_T(-t)q\circ \exp(tH_p)\,dt,\quad
J_T(t)=\frac{1}{T}J\left(\frac{t}{T}\right),
\endeq
where the function $J$ is compactly supported, smooth away from
$0$, with
$$
J'(t)=\delta(t)-K(t).
$$

Composing the principal symbol (\ref{0.14}) with the holomorphic
canonical transformation $\exp(i\eps H_{G_T})$, and conjugating
the operator $P_{\eps}$ by means of the corresponding Fourier
integral operator $U_{\eps}=e^{\frac{\eps}{h} G_T(x,hD_x)}$,
defined microlocally near $p^{-1}(0)\cap T^*M$, we may reduce our
operator to a new one, still denoted by $P_{\eps}$, which has the
principal symbol
$$
p_{\eps}\circ \exp(i\eps
H_{G_T})=p+i\eps\langle{q}\rangle_{T,K}+{\cal O}_T(\eps^2).
$$
Moreover, it is still true that $P_{\eps=0}$ is the original
selfadjoint operator. Repeating an argument, explained for example
in~\cite{Sj00}, which makes use of the sharp G\aa{}rding
inequality, we obtain a first localization of the spectrum of
$P_{\eps}$: if $z\in \comp$ in the spectrum of $P_{\eps}$ is such
that $\abs{\Re z}\leq \delta$, then as $\eps,\delta,h\rightarrow
0$,
\begeq
\label{0.16.01}
\frac{\Im z}{\eps}\in \left[ \lim_{T\rightarrow
\infty}\inf_{p^{-1}(0)}\Re \langle{q}\rangle_{T,K}-o(1),
\lim_{T\rightarrow \infty} \sup_{p^{-1}(0)}\Re
\langle{q}\rangle_{T,K}+o(1)\right].
\endeq
This estimate remains valid for $K=1_{[-1,0]}$. Let us also
notice that along the diophantine torus $\Lambda_j$, $1\leq j\leq
L$, we have uniformly, as $T\rightarrow \infty$, \begeq
\label{0.16.1} \langle{q}\rangle_T=F_j+{\cal
O}\left(\frac{1}{T}\right).
\endeq
Here $F_j$ is the mean value of $q$ over $\Lambda_j$, computed
with respect to the natural smooth flow-invariant measure on
$\Lambda_j$, with respect to which the $H_p$-flow on $\Lambda_j$
is ergodic. In the case when $K\in C^{\infty}_0(\real)$, using the rapid decay
of $\widehat{K}$, it is easy to see that (\ref{0.16.1}) improves
to
$$
\langle{q}\rangle_{T,K}=F_j+{\cal
O}\left(\frac{1}{T^{\infty}}\right).
$$
We shall assume from now on that
\begeq
\label{0.16.2}
F_1=F_2=\ldots =F_{L}\quad \wrtext{is independent of}\,\,\, j,
\endeq
and we write then $F_j=F$, $1\leq j\leq L$.

As will be explained in section 2, for each $j$, there exists a
smooth canonical transformation
$$
\kappa_{\infty,j}: {\rm neigh}(\Lambda_j,T^*M)\rightarrow {\rm
neigh}(\xi=0, T^*{\bf T}^2),
$$
mapping $\Lambda_j$ to $\xi=0$, such that \begeq \label{0.16.5}
p_{\eps}\circ \kappa_{\infty,j}^{-1}=p_{\infty,j}(\xi)+i\eps
q_j(x,\xi)+{\cal O}(\eps^2)+{\cal O}(\xi^{\infty}),
\endeq
and $p_{\infty,j}(\xi)=a\cdot \xi+{\cal O}(\xi^2)$, with $a$
satisfying (\ref{0.13}). We furthermore may assume that the energy
surface $p^{-1}_{\infty,j}(0)$ has the form $\xi_2=f_j(\xi_1)$,
for some smooth function $f_j$, with $f_j(0)=0$, $f_j'(0)\neq 0$.
Using the coordinates $\xi_1$, $\xi_2$, we define \begeq
\label{0.17} \langle{q_j}\rangle(\xi)=\frac{1}{(2\pi)^2}\int
q_j(x,\xi)\,dx.
\endeq
For each small $\delta>0$, we use the coordinate functions
$\xi_1\circ \kappa_{\infty,j}$ and $\xi_2\circ \kappa_{\infty,j}$
near $\Lambda_j$ to decompose the real energy surface as follows,
$$
p^{-1}(0)=\Omega_-(\delta)\cup \Lambda_{\delta}\cup
\Omega_+(\delta),
$$
where
$$
\Lambda_{\delta}=p^{-1}(0)\cap \bigcup_{j=1}^L \left(\xi_1\circ
\kappa_{\infty,j}\right)^{-1}((-\delta,\delta)).
$$
Here the sets $\Omega_{\pm}(\delta)$ are disjoint, compact, with
finitely many connected components, while in general they are not
invariant under the $H_p$--flow.

Recall that $F$ stands for the common value of the average of $q$
over the tori $\Lambda_j$, $1\leq j\leq L$. We introduce the
following global assumption:
\begin{eqnarray}
\label{H0} & & \hbox{There exist}\,\,N_1,\,N_2 \in \nat
 \backslash\{0\}\,\,\hbox{and a sequence}\,\,\delta=\delta_j\to 0\,\,\hbox{such that}\\
\nonumber & & \inf_{\Omega_+(\delta)}\biggl(\langle{\Re
 q}\rangle_{\delta^{-N_1},K}-\Re F\biggr) \geq \delta^{N_2},\quad \sup_{\Omega_-(\delta)}\biggl(\langle{\Re
 q}\rangle_{\delta^{-N_1},K}-\Re F\biggr) \leq -\delta^{N_2}.
\end{eqnarray}
Here $0\leq K\in C^{\infty}_0$ is as in (\ref{0.15}), and we adopt the convention that when
$N_2=1$, $\pm \delta^{N_2}$ in the right hand side of (\ref{H0}) should be replaced by
$\pm \delta/C_1$ for some $C_1>0$.

\begin{theo}
Let $\alpha_{1,j}$ and $\alpha_{2,j}$ be the fundamental
cycles in $\Lambda_j$, $1\leq j\leq L$, defined by
$$
\kappa_{\infty,j}(\alpha_{k,j})=\{x\in {\bf T}^2; x_k=0\},\quad
k=1,2.
$$
We write then $S_j=(S_{1,j},S_{2,j})$ and
$k_j=(k(\alpha_{1,j}),k(\alpha_{2,j}))$ for the values of the
actions and the Maslov indices of the cycles, respectively. Let us
make the global dynamical assumption {\rm (\ref{H0})}, and assume
that the differentials of the functions $p_{\infty,j}$ and $\Re
\langle{q_j}\rangle$, defined in {\rm (\ref{0.16.5})} and {\rm
(\ref{0.17})} are linearly independent when $\xi=0$, $1\leq j\leq
L$. Assume furthermore that $\eps={\cal O}(h^{\delta})$,
$\delta>0$, satisfies $\eps\geq h^K$, for some $K$ fixed but
arbitrarily large. Let $C>0$ be sufficiently large. Then the
eigenvalues of $P_{\eps}$ in the rectangle \begeq \label{R0}
\abs{\Re z}<\frac{h^{\delta}}{C},\quad \abs{\Im z-\eps \Re
F}<\frac{\eps h^{\delta}}{C}
\endeq
are given by
$$
P_j^{(\infty)}\left(h\left(k-\frac{k_j}{4}\right)-\frac{S_j}{2\pi},\eps;h\right)+{\cal
O}(h^{\infty}),\quad k\in \z^2,\quad 1\leq j\leq L.
$$
Here $P^{(\infty)}_j(\xi,\eps;h)$ is smooth in $\xi\in {\rm
neigh}(0,\real^2)$ and $\eps\in {\rm neigh}(0,\real)$, real-valued
for $\eps=0$, and has an asymptotic expansion in the space of such
functions,
$$
P^{(\infty)}_j(\xi,\eps;h)\sim \sum_{l=0}^{\infty} h^l
p_{j,l}^{(\infty)}(\xi,\eps),\quad 1\leq j\leq L.
$$
We have
$$
p^{(\infty)}_{j,0}(\xi,\eps)=p_{\infty,j}(\xi)+i\eps\langle{q_j}\rangle(\xi)+{\cal
O}(\eps^2).
$$
\end{theo}

Our next result treats the case when the strength of the
perturbation $\eps$ is sufficiently small but independent of $h$.
In this case, we obtain a complete spectral result in a fixed
$h$-independent domain.

\begin{theo} Let us continue to write $S_j$ and $k_j$
for the actions and Maslov indices of the fundamental cycles in
$\Lambda_j$, $1\leq j\leq L$. Assume that $h^{1/3-\delta}<\eps
\leq \eps_0\ll 1$, for some $\delta>0$. As in Theorem {\rm 1.1},
we make the assumption {\rm (\ref{H0})} and assume that the
differentials of $p_{\infty,j}(\xi)$ and $\Re
\langle{q_j}\rangle(\xi)$ are linearly independent for $\xi=0$,
$1\leq j\leq L$. Let $C>0$ be large enough. Then the eigenvalues
of $P_{\eps}$ in \begeq \label{R1} \abs{\Re z}\leq
\frac{\eps^{1/\widetilde{N}}}{C},\quad \abs{\frac{\Im z}{\eps}-\Re
F}\leq \frac{\eps^{1/\widetilde{N}}}{C}
\endeq
are given by
$$
z(j,k)\sim \sum_{n=0}^{\infty} h^n
\widetilde{p}_{j,n}^{(\infty)}\left(h\left(k-\frac{k_j}{4}\right)-\frac{S_j}{2\pi},\eps\right),\quad
k\in \z^2,\,\, 1\leq j\leq L,
$$
with
$$
\widetilde{p}_{j,n}^{(\infty)}(\xi,\eps)={\cal O}(\eps^{-2(n-1)_+
-n/\widetilde{N}}),\,\, n=0,1,2,\ldots\, , \quad 1\leq j\leq L,
$$
holomorphic for $\xi={\cal O}(\eps^{1/\widetilde{N}})$, and
$$
\widetilde{p}_{j,0}^{(\infty)}(\xi,\eps)=p_j(\xi)+i\eps
q_j(\xi,\eps)+{\cal
O}\left(\eps^{\frac{N}{\widetilde{N}}-1}\right).
$$
Here $p_j$ is real on the real domain and the differentials of
$p_j(\xi)$ and $\Re q_j(\xi,\eps)$ are linearly independent when
$\xi=\eps=0$, $1\leq j\leq L$. We have
$$
p_j(\xi)+i\eps q_j(\xi,\eps)=a\cdot \xi+i\eps F+{\cal
O}((\xi,\eps)^2),\quad a=a_j.
$$
The parameters $\widetilde{N}$ and $N/\widetilde{N}$ can be taken arbitrarily large.
\end{theo}

\Remark. As we shall see in section 2, the assumption (\ref{H0})
implies the existence of a suitable weight function, which allows
us to microlocalize the spectral problem for $P_{\eps}$ in a
rectangle of the form (\ref{R0}) or (\ref{R1}) to a small \neigh{}
of the union of the tori $\Lambda_j$, $1\leq j\leq L$---see
Proposition 2.3. Indeed, it is the existence of the weight
function that allows us to carry out the complete spectral
analysis in such rectangles. The purpose of the condition
(\ref{H0}) is to provide an explicit criterion of a purely
dynamical nature, which suffices for the construction of the
global weight. As will also be seen in section 2, in the case when
the $H_p$-flow is completely integrable, the condition (\ref{H0})
can be replaced by a slightly different, although equivalent one,
having the advantage of being easier to verify in practice---see
(\ref{1.31}) below, Proposition 2.5, and the discussion in section
7. In particular, in section 7, we show that in the completely
integrable case, the set of values $F$ in (\ref{0.16.2}) to which
Theorem 1.1 and Theorem 1.2 apply, covers the entire region
(\ref{0.16.01}), apart from a subset of a suitably small measure.
Such a result is stated in Theorem 7.6, where, using the
isoenergetic KAM theorem, we also extend it to the perturbed
situation, when the completely integrable symbol $p$ is replaced
by $p+{\cal O}(\lambda)$, with $\lambda>0$ small enough.

\vskip 2mm
The construction of quasimodes associated to invariant Diophantine
Lagrangian tori and to Cantor families of such tori has a long
tradition in the selfadjoint case, and we refer
to~\cite{Du},~\cite{CdV},~\cite{La}, and~\cite{Popov2} for the
results in this direction. Especially relevant in the present
two-dimensional context is the work of Shnirelman~\cite{Sh} in the
selfadjoint case, which contains the idea of using invariant
Lagrangian tori to split up the real energy surface into different
invariant regions. See also~\cite{CdV2}. In~\cite{Sh}, the main focus is on constructing
quasimodes associated with the gaps between the invariant tori in
a perturbative situation. In our non-selfadjoint case, the idea of
using the invariant tori as separatrices becomes more efficient
than in the standard selfadjoint setting, and in a future work we
plan to use it to study the global distribution of eigenvalues
inside the entire band corresponding to different invariant tori.

In the present work, we essentially exploit only the invariant
tori (\ref{0.12}), corresponding to a given value of the average
of the leading perturbation, and obtain complete spectral results
in the spirit of our previous works~\cite{Sj2004},~\cite{HiSj1},
and~\cite{HiSj2}. In particular, the work~\cite{HiSj2} introduced
and exploited a dynamical condition somewhat reminiscent of
(\ref{H0}). Technically speaking, many important differences
occur, however, due to the more general behavior of the classical
flow away from the union of the invariant tori (\ref{0.12}).

The plan of the paper is as follows.

In section 2, we use the global dynamical condition (\ref{H0}) to
construct a globally defined compactly supported weight function,
which allows us to microlocalize the spectral problem for
$P_{\eps}$ to a small \neigh{} of the union of the $\Lambda_j$'s,
$1\leq j\leq L$. We also give a discussion of a modified version
of (\ref{H0}), when the $H_p$-flow is completely integrable.

In section 3, we carry out a construction of a quantum Birkhoff
normal form for $P_{\eps}$, valid to an arbitrarily high order along the
invariant torus, and in $(\eps,h)$.

In section 4, we solve an appropriate Hamilton-Jacobi equation, to
be used in the spectral analysis of $P_{\eps}$ in the case when
$\eps$ is sufficiently small but independent of $h$. The ideas
used here are similar to~\cite{MeSj} and~\cite{Sj2004}.

In section 5, we justify the eigenvalue computation based on the
Birkhoff normal form construction from section 3 and prove Theorem
1.1, by solving a suitable global Grushin problem.

In section 6, we use the results of section 4 to modify the
Grushin analysis from the previous section, in order to establish
Theorem 1.2.

In section 7, we first study the completely integrable case and
verify the global dynamical condition (\ref{H0}), or rather
(\ref{1.31}) below, under some general assumptions. This is then
applied to the case of convex analytic surfaces of revolution, and
complex perturbations, close to rotationally symmetric ones. We
then exploit the isoenergetic KAM theorem and study the dynamical
condition in the case when the unperturbed symbol is close to a
completely integrable one. The results of this section are
summarized in Theorem 7.6, which together with Theorem 1.1 and
Theorem 1.2, can be considered as the final result of the present
work.

In section 8, we give an application of Theorem 1.2 to the barrier
top resonances for the semiclassical Schr\"odinger operator and obtain
an extension of the result of~\cite{KaKer} to an $h$-independent
region in the complex plane.

\vskip 2mm
\noindent
{\bf Acknowledgment}. Part of this project was conducted when the
first author vi\-si\-ted \'Ecole Polytechnique in June of 2004. He
is happy to thank its Centre de Math\'ematiques for a generous
hospitality and excellent working conditions. The partial support of
his research by the National Science Foundation under grant
DMS-0304970 is also gratefully acknowledged. All three authors have
benefited from the hospitality of the MSRI, Berkeley, in Spring 2003,
during a special semester on semiclassical analysis.

\section{Construction of the global weight} \label{section2} \setcounter{equation}{0}
In ~\cite{HiSj1} and~\cite{HiSj2}, the basic weight function was
coming from an averaging procedure along the $H_p$--flow, and was
given by $G_T$, defined in (\ref{0.16}), with $T$ being a common
period for the flow. This
weight was employed in a full \neigh{} of $p^{-1}(0)\cap T^*M$. As
we shall see in section 3, in the present case, the Birkhoff normal
form construction will in general be valid only to an infinite order along
an invariant torus, which will force us to work in small
$h$-dependent \neigh{}s of the tori. We shall therefore find it
convenient to reduce $q$ to the flow average $\langle{q}\rangle_{T,K}$
by means of $G_T$ when away from a small but fixed \neigh{} of the
union of $\Lambda_j$, $1\leq j \leq L$, while near each
$\Lambda_j$, we shall use a somewhat different sort of average,
which is close to $\langle{q}\rangle_{T,K}$, when $T$ is large.

\vskip 2mm
We keep all the assumptions of the introduction, and consider an
operator $P_{\eps}$ with the leading symbol \begeq \label{1.0}
p+i\eps q+{\cal O}(\eps^2),
\endeq
microlocally in a \neigh{} of $p^{-1}(0)\cap T^*M$.

\vskip 2mm
Let us say that a smooth multi-valued function defined on or near a torus is
grad-periodic if its gradient is single-valued.
In what follows, we shall work microlocally near a
fixed invariant torus, say $\Lambda_1$. From the assumption
(\ref{0.12}), we recall that there exist analytic and
grad-periodic coordinates on $\Lambda_1$, $x_1$ and $x_2$, which
induce an identification between $\Lambda_1$ and ${\bf T}^2$.
Applying Weinstein's tubular \neigh{} theorem, see~\cite{CdS}, or
simply following the argument, described in the beginning of
section 1 of~\cite{MeSj}, we see that we can extend $x_1$ and
$x_2$ to a tubular \neigh{} of $\Lambda_1$ and complete them with
analytic functions $\xi_1$ and $\xi_2$, defined near $\Lambda_1$
and vanishing on this set, so that we get a real analytic
canonical transformation, \begeq \label{1.1} \kappa_1: {\rm
neigh}(\Lambda_1,T^*M)\rightarrow {\rm neigh}(\xi=0,T^*{\bf T}^2),
\endeq
which maps $\Lambda_1$ to the zero section in $T^*{\bf T}^2$, and
such that when expressed in terms of the coordinates $x$ and
$\xi$, the unperturbed leading symbol $p$ of $P_{\eps=0}$ becomes
\begeq \label{1.2} p(x,\xi)=a\cdot \xi+{\cal O}(\xi^2),
\endeq
with $a$ satisfying (\ref{0.13}). It will now be convenient to perform an additional real canonical
transformation on $T^*{\bf T}^2$, in order to make $p$ independent
of $x$ to a high order in $\xi$. To this end, we observe that a
straightforward Birkhoff normal form construction, very similar to
the one described in detail in section 3, shows that there exists
a sequence of real-valued functions $G_1,G_2 \ldots,$ with $G_j$
being homogeneous of degree $j+1$ in $\xi$ and depending
analytically on $x$, such that if
$$
G\sim G_1+G_2+\ldots\, , $$ then at the level of formal Taylor
expansions in $\xi$, we have
$$
p\circ \exp(H_{G})=p_{\infty}(\xi)+{\cal O}(\xi^{\infty}),
$$
where $p_{\infty}(\xi)=a\cdot \xi+{\cal O}(\xi^2)$ depends on
$\xi$ only. See also (\ref{0.16.5}). Since we work in the analytic
category and do not wish to consider convergence questions for the
Birkhoff normal forms and associated canonical transformations (in
this connection, see~\cite{PM}), we truncate the
series at some fixed but arbitrarily large order $N$, and write
$$
p\circ \exp(H_{G_{(N)}})=p_N(\xi)+{\cal O}(\xi^{N+1}),
$$
with $G_{(N)}=G_1+G_2+\ldots +G_{N-1}$, and $p_N(\xi)=a\cdot
\xi+{\cal O}(\xi^2)$ independent of $x$. Notice that since
$G_{(N)}(\xi)={\cal O}(\xi^2)$, the analytic canonical
transformation \begeq \label{1.3} \kappa^{(N)}:=\exp(H_{G_{(N)}})
\endeq
 maps the zero section $\xi=0$ to
itself and it also preserves the action integrals along closed
loops. The composition of the transforms, $\kappa^{(N)}\circ
\kappa_1$, with $\kappa_1$ defined in (\ref{1.1}), maps
symplectically a \neigh{} of $\Lambda_1$ in $T^*M$ onto a \neigh{}
of $\xi=0$ in $T^*{\bf T}^2$, in such a way that $\Lambda_1$ is
mapped onto the zero section $\xi=0$. Implementing \begeq
\label{1.3.5} \kappa_{N,1}=\kappa^{(N)}\circ \kappa_1
\endeq
by means of a microlocally unitary Fourier integral operator with
a real phase, $U$, and conjugating $P_{\eps}$ by means of $U$, we
obtain a new operator, still denoted by $P_{\eps}$, which is
microlocally defined near the zero section in $T^*{\bf T}^2$, and
has the leading symbol \begeq \label{1.4}
p_0(x,\xi,\eps)=p(x,\xi)+i\eps q(x,\xi)+{\cal O}(\eps^2),
\endeq
with
\begeq
\label{1.5}
p(x,\xi)=p_{N}(\xi)+{\cal
O}(\xi^{N+1}),\quad p_{N}(\xi)=a\cdot \xi+{\cal O}(\xi^2).
\endeq
Here the operator $P_{\eps=0}$ is selfadjoint. Moreover, the
complete symbol of $P_{\eps}$ is a holomorphic function in a fixed
complex \neigh{} of $\xi=0$, and it depends smoothly on $\eps\in
{\rm neigh}(0,\real)$. On the operator level, $P_{\eps}$ acts in
the space of microlocally defined Floquet periodic functions on
${\bf T}^2$, $L^2_{\theta}({\bf T}^2)\subset L^2_{{\rm
loc}}(\real^2)$, elements $u$ of which satisfy \begeq
\label{1.5.1} u(x-\nu)=e^{i\theta\cdot \nu} u(x),\quad
\theta=\frac{S}{2\pi h}+\frac{k_0}{4},\quad \nu\in 2\pi \z^2.
\endeq
Here $S=(S_1,S_2)$ are the classical actions,
$$
S_j=\int_{\alpha_j}\xi\,dx,\quad j=1,2,
$$
and $\alpha_j$ form a system of fundamental cycles in $\Lambda_1$,
defined by
$$
\kappa_1(\alpha_j)=\beta_j,\quad j=1,2,\quad \beta_j=\{x\in {\bf T}^2;
x_j=0\}.
$$
The tuple $k_0=(k_0(\alpha_1),k_0(\alpha_2))$ stands for the Maslov
indices of the cycles $\alpha_j$, $j=1,2$.

\vskip 6mm \noindent In what follows we shall work with the
operator $P_{\eps}$, microlocally defined near $\xi=0$ in $T^*{\bf
T}^2$, which has the principal symbol (\ref{1.4}), (\ref{1.5}). For future
reference, let us also remark that an application of the implicit
function theorem shows that we may assume that the level set
$p_{N}=0$ is of the form
\begeq
\label{1.5.5}
\xi_2=f(\xi_1),
\endeq
for some analytic function $f$ with $f(0)=0$, $f'(0)\neq 0$.

\vskip 6mm We shall now discuss the problem of solving the
equation
\begeq
\label{1.6}
H_p G=q-r,
\endeq
considered near $\xi=0$ in $T^*{\bf T}^2$, for a suitable
remainder $r$. The torus average of the left hand side of
(\ref{1.6}) is ${\cal O}(\xi^N)$, and we shall try to make $r$
independent of $x$ at least to that order in $\xi$. Also, the
function $G$ should be analytic near $\xi=0$. When solving
(\ref{1.6}), we first consider the problem of solving \begeq
\label{1.7} H_{p_N} G=q-r.
\endeq
We analyze this problem on the level of formal Taylor series in
$\xi$, and look for $G$ in terms of a formal expansion,
$$
G\sim \sum_{k=0}^{\infty} G_k, $$ where $G_k$ is homogeneous of
degree $k$ in $\xi$. Taking also a (finite) Taylor expansion of
$p_N$ with respect to $\xi$, we easily see that
$$
H_{p_N}G\sim \sum_{n=0}^{\infty} f_n, $$ where $f_n$ is
homogeneous of degree $n$ in $\xi$ and is given by
$$
f_n=\sum_{k+l-1=n}\{p_{N,l},G_k\},
$$
with $p_{N,l}$ homogeneous of degree $l\geq 1$ in $\xi$.
Introducing also a Taylor expansion of $q$ with respect to $\xi$,
$q\sim q_0+q_1+\ldots\,,$ and using the fact that the operator $a\cdot
\partial_x$ is globally analytic hypoelliptic, we determine successively analytic
functions $G_0,G_1,\ldots$, with $G_k$ homogeneous of degree $k$
in $\xi$, such that
$$
q_n-f_n=\langle{q_n}\rangle,\quad n=0,1,\ldots.
$$
Here, as in (\ref{0.17}), we write
$$
\langle{f}\rangle(\xi)=\frac{1}{(2\pi)^2}\int f(x,\xi)\,dx,
$$
for the torus average of a smooth function $f$ defined near
$\xi=0$ in $T^*{\bf T}^2$. Taking the finite sum $G_0+G_1+\ldots
+G_N$, we get an analytic function $G$ such that
$$
H_{p_N}G=q-r,
$$
where
$$
r(x,\xi)=\langle{q}\rangle(\xi)+{\cal O}(\xi^{N+1}).
$$
To solve (\ref{1.6}), we just use that
$$
H_p G=H_{p_N}G+{\cal O}(\xi^N).
$$
We conclude therefore that for any fixed $N\in\nat$ as in
(\ref{1.5}), we can find an analytic function $G$ defined in a
fixed $N$-independent \neigh{} of $\xi=0$  and solving
(\ref{1.6}), with the remainder $r$ satisfying \begeq
\label{1.7.5} r(x,\xi)=\langle{q}\rangle(\xi)+{\cal O}(\xi^N).
\endeq
Composing $G$ with the inverse of $\kappa_{N,1}$, defined in
(\ref{1.3.5}), we get an analytic function \begeq \label{1.7.6}
G_1=G\circ \kappa_{N,1}^{-1},
\endeq
defined in a small but fixed \neigh{} of $\Lambda_1\subset T^*M$.

Coming back to the torus model, let us recall that we are working
in real symplectic coordinates $(x,\xi)$ for which (\ref{1.5})
holds. We shall consider the behavior of the $H_p$-flow near
$\xi=0$ for large, but finite times. With
$(x(t),\xi(t))=\exp(tH_p)(x(0),\xi(0))$ we have by the Hamilton
equations,
$$
\dot{\xi}(t)={\cal O}(\abs{\xi(t)}^{N+1}),
$$
and a standard "continuous induction" argument shows that for
$\abs{t}\leq {\cal O}_N(1)\abs{\xi(0)}^{-N}$, we have
$$
\xi(t)={\cal O}(\abs{\xi(0)}),
$$
provided that $\abs{\xi(0)}$ is small enough. We get, for these
times, \begeq \label{1.8} \xi(t)=\xi(0)+t{\cal
O}(\abs{\xi(0)}^{N+1}).
\endeq
Similarly, we have
\begeq
\label{1.9}
\dot{x}(t)=\frac{\partial p_N}{\partial \xi}+{\cal
O}(\abs{\xi}^N),
\endeq
and with $(x_0(t),\xi(0))=\exp(tH_{p_N})(x(0),\xi(0))$, we get for
$\abs{t}\leq {\cal O}_N(1) \abs{\xi(0)}^{-N}$, and writing $(x,\xi)$
rather than $(x(0),\xi(0))$,
\begin{equation}
\label{1.10}
\xi(t)=\xi+t{\cal O}(\abs{\xi}^{N+1}),\quad
x(t)=x_0(t)+t{\cal O}(\abs{\xi}^N).
\end{equation}
For future reference, we shall also investigate the behavior of
$\langle{q}\rangle_{T,K,p}$ near $\xi=0$, as $T\rightarrow
\infty$. In doing so, we notice that (\ref{1.10}) gives
\begeq
\label{1.10.5}
\langle{q}\rangle_{T,K,p}(x,\xi)=\langle{q}\rangle_{T,K,p_{N}}(x,\xi)+{\cal
O}(T \abs{\xi}^N), \quad T\leq {\cal O}_N(1) \abs{\xi}^{-N}.
\endeq
where we may also remark that the ${\cal O}$-term in (\ref{1.10.5}) depends only on
the Lipschitz norm of $q$. When analyzing
$\langle{q}\rangle_{T,K,p_N}$, we write
$$
\langle{q}\rangle_{T,K,p_N}(x,\xi)=\frac{1}{T}\int
K\left(\frac{-t}{T}\right) q(x+tp_N'(\xi),\xi)\,dt,
$$
and expanding $q(\cdot,\xi)$ in a Fourier series, we get
\begin{equation}
\label{1.11}
\langle{q}\rangle_{T,K,p_N}(x,\xi)=\langle{q}\rangle(\xi)+\sum_{0\neq
k\in {\rm \bf Z^2}} e^{ikx}
\widehat{q}(k,\xi)\widehat{K}(Tp'_N(\xi)\cdot k).
\end{equation}
Here $\widehat{q}(k,\xi)$ are the Fourier coefficients of $q(x,\xi)$ and
$\widehat{K}(\xi)=\int e^{-it\xi} K(t)\,dt$ is the Fourier
transform of $K$. Now the Diophantine condition (\ref{0.13}) gives
$$
\abs{p_N'(\xi)\cdot k}\geq
\frac{1}{C_0\abs{k}^{N_0}}-C_1\abs{\xi}\abs{k},\quad k\neq 0,
\quad C_1>0,
$$
and therefore
$$
\abs{p_N'(\xi)\cdot k}\geq
\frac{1}{2C_0\abs{k}^{N_0}},\,\,\wrtext{if}\,\,
C_1\abs{\xi}\abs{k}\leq \frac{1}{2C_0 \abs{k}^{N_0}},\,\, k\neq 0.
$$
Let $0\leq \chi \in C^{\infty}_0((-1,1))$ be equal to $1$ on $(-1/2,
1/2)$, and let us decompose the sum in (\ref{1.11}) as follows:
\begin{eqnarray*}
& & \langle{q}\rangle_{T,K,p_N}(x,\xi)-\langle{q}\rangle(\xi) =
\sum_{0\neq k\in {\rm \bf Z^2}}
\chi\left(2C_0C_1\abs{\xi}\abs{k}^{N_0+1}\right)e^{ix\cdot k}\widehat{q}(k,\xi)
\widehat{K}(Tp'_N(\xi)\cdot k) \\
& + & \sum_{0\neq k\in {\rm \bf Z^2}} (1-\chi)\left(2C_0
  C_1\abs{\xi}\abs{k}^{N_0+1}\right)e^{ix\cdot k}
  \widehat{q}(k,\xi)\widehat{K}(Tp'_N(\xi)\cdot k) = {\rm I}+{\rm II},
\end{eqnarray*}
with the natural definitions of ${\rm I}$ and ${\rm II}$. We get for
each $M\in \nat$,
$$
\abs{I}\leq {\cal O}_M(1)\sum_{k\neq
0}\frac{\chi\left(2C_0C_1\abs{\xi}\abs{k}^{N_0+1}\right)}{\abs{Tp'_N(\xi)\cdot
k}^M}\abs{\widehat{q}(k,\xi)}\leq {\cal O}_M(1)T^{-M}\sum_{k\neq
0}\abs{k}^{M N_0}\abs{\widehat{q}(k,\xi)},
$$
and therefore,
$$
I={\cal O}(T^{-\infty}).
$$
When treating ${\rm II}$, we only use that $\widehat{K}$ is bounded and because
of the presence of $1-\chi$, we have $\abs{k}\geq {\cal O}(1)^{-1}
\abs{\xi}^{-1/(N_0+1)}$ in the sum, so that the contribution coming from ${\rm
II}$ is ${\cal O}(\abs{\xi}^{\infty})$. We conclude that
\begin{equation}
\label{1.12}
\langle{q}\rangle_{T,K,p_N}(x,\xi)=\langle{q}\rangle(\xi)+{\cal
O}(\xi^{\infty}+T^{-\infty}).
\end{equation}

The discussion above is summarized in the following proposition.

\begin{prop}
Let $(x,\xi)$ be real symplectic coordinates so that {\rm (\ref{1.5})}
is true. Here $N$ is fixed but arbitrarily large. We then have for
$\xi$ small enough,
\begin{equation}
\label{1.13}
\langle{q}\rangle_{T,K,p}(x,\xi)=\langle{q}\rangle(\xi)+{\cal
O}\left(\xi^{\infty}+T^{-\infty}\right)+{\cal O}(T
\abs{\xi}^N),\quad T\leq {\cal O}_N(1)\abs{\xi}^{-N}.
\end{equation}
\end{prop}

\vskip 2mm
For future reference, we shall pause here
to recall some general estimates for convolutions of the form
$K_T*g$, where $g\in L^{\infty}(\real)$. The following discussion
is motivated by the fact that the flow average of $q$,
$$
\langle{q}\rangle_{T,K}(\rho)=\int K_T(-t)q(\exp(tH_p)(\rho))\,dt
$$
is a convolution in the time variable along the $H_p$-trajectory
passing through $\rho$. The starting point is that \begeq
\label{1.14new} K=K_{\eps}*K+r_{\eps},\quad
\norm{r_{\eps}}_{L^1}={\cal O}(\eps),
\endeq
for $0<\eps \leq 1$. It follows that for $0<S\leq T$, \begeq
\label{1.15new} K_T=K_S*K_T+r_{\frac{S}{T},T},
\endeq
where
$r_{\frac{S}{T},T}(t)=\frac{1}{T}r_{\frac{S}{T}}(\frac{t}{T})$, so
that
$$
\norm{r_{\frac{S}{T},T}}_{L^1}=\norm{r_{\frac{S}{T}}}_{L^1}={\cal
O}\left(\frac{S}{T}\right).
$$
If $g\in L^{\infty}(\real)$ is real-valued, we have \begeq
\label{1.16new} K_T*g=K_T*K_S*g+r_{\frac{S}{T},T}*g,
\endeq
with \begeq \label{1.17new}
\norm{r_{\frac{S}{T},T}*g}_{L^{\infty}}\leq {\cal
O}\left(\frac{S}{T}\right)\norm{g}_{L^{\infty}}.
\endeq
In particular, when $\norm{g}_{L^{\infty}}={\cal O}(1)$, we get
\begeq \label{1.18new} \inf_{{\rm \bf R}}(K_T*g)\geq \inf_{{\rm
\bf R}}(K_S*g)-{\cal O}\left(\frac{S}{T}\right),
\endeq
and \begeq \label{1.19new} \sup_{{\rm \bf R}}(K_T*g)\leq
\sup_{{\rm \bf R}}(K_S*g)+{\cal O}\left(\frac{S}{T}\right).
\endeq
Here we can be a little more precise about where to take the sup
and the inf. If $I\subset \real$ is an interval, then \begeq
\label{1.20new} \inf_I (K_T*g)\geq \inf_{I-T{\rm supp} K}
(K_S*g)-{\cal O}\left(\frac{S}{T}\right),
\endeq
\begeq \label{1.21new} \sup_I (K_T*g)\leq \sup_{I-T{\rm supp}
K}(K_S*g)+{\cal O}\left(\frac{S}{T}\right).
\endeq
We also notice that if $\widetilde{K}$ has the same properties as
$K$, then (\ref{1.14new}) can be generalized to
$$
K-K*\widetilde{K}_{\eps}={\cal O}(\eps)\quad \wrtext{in}\,\, L^1,
$$
leading to the possibility of replacing $K_S$ by $\widetilde{K}_S$
in (\ref{1.15new}), (\ref{1.16new}), and
(\ref{1.18new})--(\ref{1.21new}).

\vskip 4mm Let us return now to the operator $P_{\eps}$ with the
leading symbol (\ref{1.4}), (\ref{1.5}). We shall assume from now
on that \begeq \label{1.13.001} dp_{N}(0)=a\,\,\wrtext{and}\,\,
d\Re \langle{q}\rangle(0)\,\,\,\wrtext{are linearly independent}.
\endeq
Notice that this condition is independent of $N$. Then, possibly
after changing the sign of the $\xi_1$-coordinate and using
(\ref{1.5.5}), we get \begeq \label{1.13.01} \frac{d}{d \xi_1}
\bigg|_{\xi_1=0}\langle{\Re q}\rangle(\xi_1,f(\xi_1))>0.
\endeq

\begin{lemma}
Assume that {\rm (\ref{1.13.001})} holds true, and let $N$ in {\rm
(\ref{1.5})} be sufficiently large. There exists a constant
$C_0>0$ such that if {\rm (\ref{H0})} holds for some $1\leq
N_1,N_2$, and $\widetilde{N}_1\in \nat\backslash\{0\}$, then there
exists $\delta_0>0$ such that
\begin{equation}
\label{1.13.02}
 \inf_{\Omega_{+}(\delta)}\left(\langle{\Re
q}\rangle_{\delta^{-\widetilde{N}_1},K}-\Re F\right)\geq
\frac{\delta}{C_0},\quad 0<\delta\leq \delta_0,
\end{equation}
and
\begin{equation}
\label{1.13.03} \sup_{\Omega_{-}(\delta)}\left(\langle{\Re
q}\rangle_{\delta^{-\widetilde{N}_1},K}-\Re F\right)\leq -
\frac{\delta}{C_0},\quad 0<\delta\leq \delta_0.
\end{equation}

\end{lemma}
\begin{proof}
It suffices to prove only (\ref{1.13.02}), and we may also assume
that $F=0$. Moreover, as will be clear, the following argument
will not depend on the choice of $N$ in (\ref{1.5}), provided that
it is large enough, and to simplify the discussion notationally,
we shall take $N=\infty$. We shall work near one of the tori
$\Lambda_j$. From Proposition 2.1 we
then know that for every fixed $N_1\in \nat\backslash\{0\}$, we
have
$$
\langle{q}\rangle_{\delta^{-N_1},K}(x,\xi)=\langle{q}\rangle(\xi)+{\cal
O}(\delta^{\infty}),\quad \abs{\xi}< \delta,
$$
if $0<\delta\leq \delta(N_1)>0$. Take $\delta$ as in (\ref{H0})
and let $\widetilde{\delta}\in [\delta^{M_1},\delta^{M_0}]$, where
$1<M_0<M_1$ will be fixed later. Then
\begin{equation}
\label{1.13.1}
\inf_{\Omega_{+}(\widetilde{\delta})}\langle{\Re q}\rangle_{\widetilde{\delta}^{-\widetilde{N}_1},K}\geq
\min \left(\inf_{\Omega_{+}(\widetilde{\delta})\backslash
\Omega_{+}(2\delta)}\langle{\Re q}\rangle_{\widetilde{\delta}^{-\widetilde{N}_1},K},
\inf_{\Omega_{+}(2\delta)}
\langle{\Re q}\rangle_{\widetilde{\delta}^{-\widetilde{N}_1},K}\right),
\end{equation}
and we also get
\begin{equation}
\label{1.13.2} \inf_{\Omega_{+}(\widetilde{\delta})\backslash
\Omega_{+}(2\delta)}\langle{\Re
q}\rangle_{\widetilde{\delta}^{-\widetilde{N}_1},K}\geq
\frac{\widetilde{\delta}}{C_0},
\end{equation}
for some $C_0>0$, when $\delta$ is small enough. We may assume
that $C_0\geq 2$.

If we choose $M_0$ with $M_0\widetilde{N}_1\geq N_1$, then
$\widetilde{\delta}^{-\widetilde{N}_1}\geq \delta^{-N_1}$, and
using the general estimate (\ref{1.20new}) in the form
$$
\langle{q}\rangle_{T,K}(\rho)\geq \inf_{t\in [0,T]}
\langle{q}\rangle_{S,K}(\exp(tH_p)(\rho))-{\cal
O}\left(\frac{S}{T}\right),
$$
and the fact that in the region $\abs{\xi}<\delta$, the
$H_p$-trajectories are approximately confined to the tori
$\xi={\rm Const}$, for times bounded by any fixed inverse power of
$\abs{\xi}$, we get
\begin{equation}
\label{1.13.3} \inf_{\Omega_{+}(2\delta)}\langle{\Re
q}\rangle_{\widetilde{\delta}^{-\widetilde{N}_1},K}\geq
\inf_{\Omega_{+}(\delta)}\langle{\Re
q}\rangle_{\delta^{-N_1},K}-{\cal
O}(\delta^{-N_1}\widetilde{\delta}^{\widetilde{N}_1}) \geq
\delta^{N_2}-{\cal
O}(\delta^{-N_1}\widetilde{\delta}^{\widetilde{N}_1}).
\end{equation}
Here we should replace $\delta^{N_2}$ by $\delta/C_1$ if $N_2=1$.
Now assume that $M_0\geq N_2$, $M_0\widetilde{N}_1-N_1>N_2$, so
that $\delta^{N_2}\geq \widetilde{\delta}$,
$\delta^{-N_1}\widetilde{\delta}^{\widetilde{N}_1}\ll
\delta^{N_2}$. Then (\ref{1.13.3}) implies that
$$
\inf_{\Omega_{+}(2\delta)}\langle{\Re
q}\rangle_{\delta^{-\widetilde{N}_1},K}\geq
\frac{\widetilde{\delta}}{2}.
$$
Combining this with (\ref{1.13.1}) and (\ref{1.13.2}), we get
$$
\inf_{\Omega_{+}(\widetilde{\delta})}\langle{\Re q}\rangle_{\widetilde{\delta}^{-\widetilde{N}_1},K}\geq
\frac{\widetilde{\delta}}{C_0},
$$
for $\delta>0$ sufficiently small and for $\widetilde{\delta}\in
[\delta^{M_1},\delta^{M_0}]$, provided that
$$
M_1>M_0\geq \max\left(N_2,\frac{N_1+N_2}{\widetilde{N}_1}\right),
$$
and $\delta>0$ is small enough depending on the exponents and such
that the estimate (\ref{H0}) holds true.

We conclude that we have the estimate in (\ref{H0}) with
$N_1=\widetilde{N}_1$ and $N_2=1$ for $\delta$ replaced by
$\widetilde{\delta}\in [\delta^{M_1}, \delta^{M_0}]$. A special
case of this is when we start with $\delta$ with
$N_1=\widetilde{N}_1$, $N_2=1$, and then conclude that (\ref{H0})
holds with $\delta$ replaced by $\widetilde{\delta}$ in
$[\delta^{M_1},\delta^{M_0}]$, when $M_1>M_0\geq
(\widetilde{N}_1+1)/{\widetilde{N}_1}$. Choose $M_1=M_0^2$. Then
the argument can be iterated and we get (\ref{H0}) with
$N_1=\widetilde{N}_1$, $N_2=1$, for $\delta$ replaced by any
$\widetilde{\delta}$ in $[\delta^{M_0^2},\delta^{M_0}]\cup
[\delta^{M_0^3},\delta^{M_0^2}]\cup \ldots = (0,\delta^{M_0}]$,
and the lemma follows.
\end{proof}

\bigskip
Lemma 2.2 shows that (\ref{H0}) is independent of the choice of
$N_1$ and $N_2$. The proof also shows that it is independent of
$K\in C^{\infty}_0(\real)$. We also notice that the proof shows that the
validity of (\ref{H0}) is stable under small perturbations of $p$
that preserve the invariant tori $\Lambda_j$ and the Diophantine
condition (\ref{0.13}).

\bigskip

We shall now construct a global weight function $\widetilde{G}$.
When doing so, let us recall from (\ref{1.7.6}) the analytic
functions $G_j=G\circ \kappa_{N,j}^{-1}$, defined in small
neighborhoods of $\Lambda_j$, $1\leq j\leq L$. Here $\kappa_{N,j}$
is defined exactly as $\kappa_{N,1}$ in (\ref{1.3.5}). When
constructing the global weight $\widetilde{G}$, we shall glue
together the analytic functions $G_T$, defined in (\ref{0.16}) and
$G_j$. In doing so, it will be sufficient to work in a fixed
sufficiently small \neigh{} of $\xi=0$ in $T^*{\bf T}^2$. Let
$0\leq \chi\leq 1$, $\chi\in C^{\infty}_0(\real^2)$, be supported
in $\abs{\xi}<2$, with $\chi=1$ on $\abs{\xi}<1$. With $0<\mu\ll
1$ to be determined, we put \begeq \label{1.13.4}
\widetilde{G}=\left(1-\sum_{j=1}^L \chi\left(\frac{\xi\circ
\kappa_{N,j}}{\mu}\right)\right)G_T+\sum_{j=1}^L
\chi\left(\frac{\xi\circ \kappa_{N,j}}{\mu}\right)G_j.
\endeq
In the following calculation, in order to ease the notation, we
shall take $L=1$, in which case we may omit $\kappa_{N,1}$ from the
notation and write \begeq
\label{1.13.4.5} \widetilde{G}=(1-\chi_{\mu})G_T+\chi_{\mu}G,
\endeq
where $\chi_{\mu}(\xi)=\chi(\xi/\mu)$. We shall compute $H_p
\widetilde{G}$, and in doing so we notice first that
$$
H_p \chi_{\mu}={\cal O}(\xi^N),
$$
uniformly in $\mu>0$. We get
$$
H_p
\widetilde{G}=(1-\chi_{\mu})\left(q-\langle{q}\rangle_{T,K,p}\right)+\chi_{\mu}\left(q-\langle{q}\rangle\right)+{\cal
O}(\xi^N)\chi_{\mu}+H_p\chi_{\mu}\left(G-G_T\right),
$$
which we rewrite as
\begeq
\label{1.13.5}
H_p
\widetilde{G}=q-\langle{q}\rangle-(1-\chi_{\mu})(\langle{q}\rangle_{T,K,p}-\langle{q}\rangle)+{\cal
O}(\xi^N)\chi_{\mu}+H_p\chi_{\mu}(G-G_T).
\endeq
Now $\chi_{\mu}={\cal O}(1)$ and (\ref{0.16}) gives $G-G_T={\cal
O}(T)$, so that
$$
H_p
\widetilde{G}=q-\langle{q}\rangle-(1-\chi_{\mu})(\langle{q}\rangle_{T,K,p}-\langle{q}\rangle)+{\cal
O}(T\xi^N),\quad T\geq 1.
$$
Proposition 2.1 gives then that uniformly in $\mu>0$ we have
\begeq \label{1.13.6} H_p
\widetilde{G}=q-\langle{q}\rangle-(1-\chi_{\mu})\left({\cal
O}(T\xi^N)+{\cal O}(\xi^{\infty}+T^{-\infty})\right)+{\cal
O}(T\xi^N),\quad T\leq {\cal O}_N(1){\abs{\xi}^{-N}}.
\endeq
We let first $N$ be sufficiently large depending on $N_1$, and
take then $\delta>0$ small enough but fixed. Take then $\mu \sim
\delta$, and put
$$
T=\mu^{-N_1},
$$
In the region $\abs{\xi}< \mu$, we have $\widetilde{G}=G$, while
$\widetilde{G}=G_T$ for $\abs{\xi}>2\mu$. In
the intermediate region where $0<\chi_{\mu}<1$, we have
$\abs{\xi}\sim \mu$, and it follows that in this region
$$
q-H_p\widetilde{G}=\langle{q}\rangle(\xi)+{\cal O}(\xi^{N-N_1}).
$$

The discussion above can be summarized in the following
proposition.

\begin{prop}
Assume {\rm (\ref{1.13.001})} for each $1\leq j\leq L$, and {\rm
(\ref{H0})} for some $N_1$ and $N_2\in \nat\backslash\{0\}$. When
$1\leq j\leq L$, let $G_j$ be an analytic solution near
$\Lambda_j$ of the equation {\rm (\ref{1.6})}, with $r$ satisfying
{\rm (\ref{1.7.5})} for some $N$ large enough. Let $W_j$ be a
sufficiently small \neigh{} of $\Lambda_j$, $W_j=(\xi_1\circ
\kappa_{N,j})^{-1}((-\delta_0,\delta_0))$, $\delta_0>0$. Then
there exists a smooth compactly supported function $\widetilde{G}$
on all of $T^*M$, such that $\widetilde{G}=G_j$ in $W_j$, $1\leq
j\leq L$, and if $W=\cup_{j=1}^L W_j$, $\delta(W)>0$ is
sufficiently small, and $\widetilde{\Omega}_{\pm}(\delta(W))$
stand for the components of $\Omega_{\pm}(\delta_0)$ in
$p^{-1}\left((-\delta(W),\delta(W))\right)\backslash W$, we have
$$
\inf_{\widetilde{\Omega}_+(\delta(W))}\left(\Re(q-H_p
\widetilde{G}-F)\right) \geq \frac{1}{C(W)},
$$
and
$$
\sup_{\widetilde{\Omega}_-(\delta(W))}\left(\Re(q-H_p
\widetilde{G}-F)\right)\leq \frac{-1}{C(W)}.
$$
Here $F$ is the average of $q$ over $\Lambda_j$, $1\leq j\leq L$,
and $C(W)>0$.
\end{prop}

\bigskip
\noindent We shall now discuss a modified version of the global
assumption (\ref{H0}), available if we assume that the $H_p$-flow
is completely integrable. In particular, the results of the
following discussion will be applied in section 7, where, as an
example, we consider $p$ corresponding to the geodesic flow on
an analytic surface of revolution.

Assuming the complete integrability for the flow, we know that
there exists an analytic real-valued function $f$, which Poisson
commutes with $p$, and such that each $\Lambda_j$, $1\leq j\leq
L$, is a level set for the associated mapping, \begeq
\label{1.20} (p,f): T^*M\rightarrow \real^2,
\endeq
corresponding to a regular value. Then each torus $\Lambda_j$ is
embedded in a Lagrangian foliation of $H_p$-invariant tori, given
by
$$
\Lambda_{E,F}: p=E,\quad f=F,
$$
for $(E,F)\in {\rm neigh}((0,F_j),\real^2)$, where $(0,F_j)$ is a
regular value of (\ref{1.20}). We may then introduce the
action-angle coordinates, given by a real analytic canonical
transformation,
$$
\kappa_j: {\rm neigh}(\Lambda_j,T^*M)\rightarrow {\rm
neigh}(\xi=0,T^*{\bf T}^2),\quad 1\leq j\leq L,
$$
such that $\Lambda_j$ is mapped to the zero section in $T^*{\bf
T}^2$, and such that when expressed in terms of the coordinates
$x$ and $\xi$, the unperturbed leading symbol $p$ becomes a
function of $\xi$ only, with
\begin{equation}
\label{1.21.1} p(\xi)=a\cdot \xi+{\cal O}(\xi^2),\quad a=a_j,
\end{equation}

As in (\ref{1.5.5}), we may assume that the energy surface
$p^{-1}(0)$ takes the form $\xi_2=f_j(\xi_1)$, where $f_j$ is
analytic with $f_j(0)=0$. The tori $\Lambda_{0,F}$, for $F\in {\rm
neigh}(F_j,\real)$ then take the form $\xi_1=\mu,\,\,
\xi_2=f_j(\mu)$, for $\abs{\mu}<b$, with $0<b\ll 1$. Consider a
flow-invariant neighborhood of the union of the $\Lambda_j$ of the
form
$$
\bigcup_{\abs{\mu}<b} \bigcup_{j=1}^L
\left(\xi_1\circ\kappa_j\right)^{-1}(\mu)\cap
p^{-1}(0)=\bigcup_{\abs{\mu}<b}  \Lambda_{\mu}.
$$
Here $\Lambda_{\mu}$ is a disjoint union of $L$ flow-invariant
tori. We decompose the real energy surface as follows,
$$
p^{-1}(0)=\Lambda_{-b}\cup \bigcup_{\abs{\mu}<b} \Lambda_{\mu}
\cup \Lambda_b,
$$
where $\Lambda_{\pm b}$ are two disjoint compact flow-invariant
domains with at most finitely many connected components. This
gives a decomposition of the energy surface \begeq \label{1.21.5}
p^{-1}(0)=\bigcup _{\mu \in M} \Lambda_{\mu},
\endeq
where $M=[-b,b]$ and each $\Lambda_{\mu}$ is a compact
$H_p$-invariant set, with finitely many connected components
$\Lambda_{\mu,j}$. Assume also the continuity property:
\begin{eqnarray}
\label{1.22} & & \wrtext{For any}\,\,\, \mu_0\in M, \eps>0
\,\,\,\wrtext{there exists}\,\, \delta>0, \,\,\wrtext{such that
for } \\  \nonumber & & {\rm dist}(\mu,\mu_0)<\delta \,\,
\wrtext{we have}\,\,\, \Lambda_{\mu}\subset
\Lambda_{\mu_0,\eps}:=\{\rho\in p^{-1}(0); {\rm
dist}(\rho,\Lambda_{\mu_0})< \eps.\}
\end{eqnarray}

As in (\ref{0.15}), we put \begeq \label{1.23}
\langle{q}\rangle_{T,K}(\rho)=\int q(\exp(tH_p)(\rho))K_T(-t)\,dt.
\endeq
Let $Q_T(\mu,j)$ and $Q_T(\mu)$ be the ranges of $\Re
\langle{q}\rangle_{T,K}$ restricted to $\Lambda_{\mu,j}$ and
$\Lambda_{\mu}$, respectively. Then $Q_T(\mu)$ is a finite union
of the closed intervals $Q_{T}(\mu,j)$, and according to
(\ref{1.20new}) and (\ref{1.21new}), we have \begeq \label{1.24}
Q_T(\mu)\subset Q_S(\mu)+{\cal
O}(1)\norm{q|_{\Lambda_{\mu}}}_{L^{\infty}}\left[-\frac{S}{T},\frac{S}{T}\right],
\endeq
uniformly for $0<S\leq T$, $\mu \in M$. Let $Q_{\infty}(\mu,j)$ be
the non-empty intersection of all the $Q_T(\mu,j)$ for $T\geq 1$,
and define $Q_{\infty}(\mu)$ to be the union of all the
$Q_{\infty}(\mu,j)$, so that \begeq \label{1.25}
Q_{\infty}(\mu)=\bigcup_j \left[\lim_{T\rightarrow \infty}
\inf_{\Lambda_{\mu,j}} \Re \langle{q}\rangle_{T,K},
\lim_{T\rightarrow \infty} \sup_{\Lambda_{\mu,j}}\Re
\langle{q}\rangle_{T,K}\right].
\endeq
We may also remark that according to (\ref{1.24}),
$Q_{\infty}(\mu)\subset Q_T(\mu)$, for all $T\geq 1$. It is also
easy to see that $Q_{\infty}(\mu)$ does not depend on the choice
of $K$. Put \begeq \label{1.26} {\cal Q}_T=\{ (\mu,E); \mu\in M,
E\in Q_T(\mu)\}, \quad 0< T \leq \infty.
\endeq

\begin{lemma}
We have
\begin{enumerate}
\item ${\cal Q}_T$ is closed.

\item For every \neigh{} ${\cal U}$ of ${\cal Q}_{\infty}$ there
      exists a $T_0\in (0,\infty)$ such that ${\cal Q}_T\subset {\cal
      U}$ for $T\geq T_0$.
\end{enumerate}
\end{lemma}
\begin{proof}
We prove (1) first. It is clear that ${\cal Q}_T$ is closed for
each finite $T$. Let $(\mu_j,E_j)\in {\cal Q}_{\infty}$ be a
convergent sequence so that $\mu_j\rightarrow \mu_0\in M$,
$E_j\rightarrow E_0\in \real$. Let $\eps>0$. Then there exists
$T_0\in (0,\infty)$ such that $Q_{T_0}(\mu_0)\subset
Q_{\infty}(\mu_0)+[-\eps,\eps]$. We fix such a number. If
$\delta>0$ is small enough, we have \begeq \label{1.27}
Q_{\infty}(\mu)\subset Q_{T_0}(\mu)\subset
Q_{T_0}(\mu_0)+[-\eps,\eps]\subset
Q_{\infty}(\mu_0)+[-2\eps,2\eps],
\endeq
for ${\rm dist}(\mu,\mu_0)<\delta$. In particular, $E_0\in
Q_{\infty}(\mu_0)+[-2\eps,2\eps]$, and letting $\eps\rightarrow 0$, we
get $E_0\in Q_{\infty}(\mu_0)$, proving the closedness of ${\cal
Q}_{\infty}$.

We now establish (2). If $T=kT_0$, $k\geq 1$, we get from (\ref{1.24})
and (\ref{1.27}) that for ${\rm dist}(\mu,\mu_0)<\delta$,
$$
Q_{T}(\mu)\subset Q_{\infty}(\mu_0)+[-2\eps,2\eps]+\left[-\frac{{\cal
O}(1)}{k},\frac{{\cal O}(1)}{k}\right],
$$
and hence
$$
Q_T(\mu)\subset Q_{\infty}(\mu_0)+[-3\eps,3\eps],
$$
for $T$ large enough. This means that $\{(\mu,E)\in {\cal Q}_T;
{\rm dist}(\mu,\mu_0)<\delta\}$ is within a distance
$3\eps+\delta$ from ${\cal Q}_{\infty}$. The statement (2) now
follows by a covering argument.
\end{proof}

\medskip
If $\abs{\mu}<b$, then each connected component $\Lambda_{\mu,j}$,
$1\leq j\leq L$, of $\Lambda_{\mu}$ is diffeomorphic to ${\bf
T}^2$ in such a way that $H_p|_{\Lambda_{\mu,j}}$ becomes
$a_1\partial_{x_1}+a_2\partial_{x_2}$, $a_j\in \real$, and we have
\begeq \label{1.28} \langle{q}\rangle_{T,K}(x)=\sum_{k\in {\rm \bf
Z}^2}\widehat{K}(Ta\cdot k) \widehat{q}(k,\mu) e^{ix\cdot k},
\endeq
which, as $T\rightarrow \infty$, converges uniformly in $x$ (using
now also the smoothness of $q$), to
\begeq \label{1.29}
\langle{q}\rangle_{\infty,\mu}(x)=\sum_{a\cdot k=0}
\widehat{q}(k,\mu)e^{ix\cdot k}.
\endeq
Notice that $Q_{\infty}(\mu,j)$ is the range of this function.

When $\abs{\mu}<b$ and $1\leq j\leq L$, we shall now study the
compact interval $Q_{\infty}(\mu,j)$ in more detail. In doing so,
we remark that if we have an equation of the form (\ref{1.6}) on
$\cup_{\abs{\mu}<b}\Lambda_{\mu,j}$, with $G$ smooth and bounded,
then $Q_{\infty}(\mu,j)\subset \Re r(\Lambda_{\mu})$. In fact,
with $K=1_{[-1,0]}$, we have for $\rho\in \Lambda_{\mu,j}$,
$$
\langle{r}\rangle_T(\rho)=\langle{q}\rangle_{T}(\rho)-
\frac{1}{T}\left(G(\exp(tH_p)(\rho))-G(\rho)\right)=\langle{q}\rangle_T(\rho)+{\cal
O}\left(\frac{1}{T}\right),
$$
and it suffices to let $T\rightarrow \infty$, and use that $\Re
\langle{r}\rangle_T(\Lambda_{\mu})\subset \Re r(\Lambda_{\mu})$.

Recall now that we have seen that we can solve (\ref{1.6}) with
$r$ satisfying (\ref{1.7.5}), and we conclude that
\begin{equation}
\label{1.29.5} Q_{\infty}(\mu,j)\subset \Re
\langle{q_j}\rangle(\mu,f(\mu))+{\cal O}(\mu^N)[-1,1],\quad
\abs{\mu}<b,\quad q_j=q.
\end{equation}

As in the general case, from now on we introduce the assumption
(\ref{1.13.001}) for each $1\leq j \leq L$. In view of
(\ref{1.29.5}), this implies (possibly after shrinking $b$ again)
that for $0<\widetilde{b}<b$, \begeq \label{1.30} \inf
\bigcup_{\mu\in [\widetilde{b},b)} Q_{\infty}(\mu)>\Re
F_0+\frac{\widetilde{b}}{C},\quad \sup \bigcup_{\mu\in
(-b,-\widetilde{b}]} Q_{\infty}(\mu)<\Re
F_0-\frac{\widetilde{b}}{C},
\endeq
for some constant $C>0$. The global assumption in the completely
integrable case is then the following one: For any
$\widetilde{b}\in (0,b)$ there exists $C(\widetilde{b})>0$ such
that \begeq \label{1.31} \inf_{\mu\in [\widetilde{b},b]}
Q_{\infty}(\mu)>\Re F_0+\frac{1}{C(\widetilde{b})},\quad
\sup_{\mu\in [-b,-\widetilde{b}]} Q_{\infty}(\mu)<\Re
F_0-\frac{1}{C(\widetilde{b})}.
\endeq

Using (\ref{1.29.5}) together with Proposition 2.1 and Lemma 2.4
it is easy to see that in the completely integrable case, the
condition (\ref{1.31}) is equivalent to (\ref{H0}).

\begin{prop}
Assume that the $H_p$-flow is completely integrable, so that {\rm
(\ref{1.21.1})} holds true. Assume furthermore that {\rm
(\ref{1.13.001})} and {\rm (\ref{1.31})} are valid. Then we have
the same conclusion as in Proposition {\rm 2.3}.
\end{prop}

\bigskip
In what follows, to ease the notation, we shall drop the tilde and
write $G$ instead of $\widetilde{G}$, defined in Proposition 2.3.
Associated with $G\in C^{\infty}_0(T^*M)$, there is a globally
defined IR-manifold \begeq \label{1.31.5} \Lambda_{\eps
G}=\{\rho+i\eps H_{G}(\rho);\,\, \rho\in T^*M\}\subset
T^*\widetilde{M},
\endeq
and when acting on the Hilbert space $H(\Lambda_{\eps G})$,
associated to $\Lambda_{\eps G}$ by means of the FBI-Bargmann
transform
$$
Tu(x)=Ch^{-3/2}\int e^{\frac{i\varphi(x,y)}{h}}u(y)\,dy,\quad
\varphi(x,y)=\frac{i}{2}(x-y)^2,
$$
the operator $P_{\eps}$ gets the leading symbol
\begeq
\label{1.32}
\left(p+i\eps q+{\cal O}(\eps^2)\right)\left(\rho+i\eps
H_{G}\rho\right)=p(\rho)+i\eps\left(q-H_p G\right)(\rho)+{\cal
O}(\eps^2).
\endeq
Here we are tacitly assuming that $M=\real^2$ so that $T^*M$ is a
linear space. In the mani\-fold case, using that $G$ is analytic
near $\Lambda_j$, $1\leq j\leq L$, we define $\Lambda_{\eps
G}=\exp(i\eps H_G)(T^*M)$ in a complex \neigh{} of $\Lambda_j$,
and elsewhere we let $\widetilde{G}\in
C^{\infty}_0(T^*\widetilde{M})$ stand for an almost holomorphic
extension of $G$, and $\Lambda_{\eps G}$ is then defined as
$$
\Lambda_{\eps G}=\exp(\eps H^{{\rm Im}\sigma}_{{\rm
Re}\widetilde{G}})(T^*M). $$ Here the Hamilton vector field of
$\Re \widetilde{G}$ is computed with respect to the imaginary part
$\Im \sigma$ of the complex symplectic form $\sigma$ on
$T^*\widetilde{M}$.

It follows from Proposition 2.3 that when $m\in \Lambda_{\eps G}$
is away from a small but fixed \neigh{} of the union of the
Lagrangian tori $\widetilde{\Lambda}_j$, such that each
$\widetilde{\Lambda}_j$ is the image of $\Lambda_{j}$ in
$\Lambda_{\eps G}$ under the map $T^*M\ni \rho\mapsto m=\rho+i\eps
H_{G}(\rho)$, it is true that
$$
\abs{\Im P_{\eps}(m)-\eps \Re F_0}\geq \frac{\eps}{{\cal O}(1)},
$$
provided that $\abs{\Re P_{\eps}(m)}\leq 1/C$ for a sufficiently
large $C>0$. Moreover, microlocally near $p^{-1}(0)$, the operator
$P_{\eps}$ acting on $H(\Lambda_{\eps G})$ is unitarily equivalent
to an operator acting on $L^2(M)$, which has the leading symbol
(\ref{1.32}).

Summing up the discussion of this section, we have achieved that
microlocally near each $\widetilde{\Lambda}_j\subset \Lambda_{\eps
G}$, $1\leq j\leq L$, the operator
$$
P_{\eps}: H(\Lambda_{\eps G})\rightarrow H(\Lambda_{\eps G})
$$
is unitarily equivalent to an operator $\widetilde{P}_{\eps}$,
acting on $L^2_{\theta}({\bf T}^2)$ and defined microlocally near $\xi=0$
in $T^*{\bf T}^2$. Here $\widetilde{P}_{\eps}$ is such that
\begin{equation}
\label{1.33} \widetilde{P}_{\eps}\sim \sum_{\nu=0}^{\infty} h^{\nu}
\widetilde{p}_{\nu}(x,\xi,\eps),
\end{equation}
with $\widetilde{p}_{\nu}$ holomorphic in a fixed complex \neigh{} of
$\xi=0$ and
\begin{equation}
\label{1.34}
 \widetilde{p}_0=p_N(\xi)+i\eps
\langle{q}\rangle(\xi)+{\cal O}(\eps^2)+{\cal O}(\xi^{N+1})+\eps
{\cal O}(\xi^N),
\end{equation}
with $p_N(\xi)=a\cdot \xi+{\cal O}(\xi^2)$, and where the term
${\cal O}(\xi^{N+1})$ is real on the real domain. Here $N$ is
fixed but arbitrarily large. In what follows, we shall drop the
tildes and write $P_{\eps}$ and $p_{\nu}$, $\nu \geq 0$, instead of
$\widetilde{P}_{\eps}$ and $\widetilde{p}_{\nu}$, $\nu\geq 0$,
respectively.

\section{The normal form construction}
\label{section3} \setcounter{equation}{0}

\vskip 4mm We recall that we have reduced the analysis to the
operator $P_{\eps}$ with a complete symbol (\ref{1.33}) and a
leading symbol (\ref{1.34}), and recall also the assumption
(\ref{1.13.001}). Our goal in this section is to construct a
quantum Birkhoff normal form for $P_{\eps}$---see
also~\cite{CdV},~\cite{Popov2}, and~\cite{Sj92}. We may write
\begeq \label{2.1}
p_0(x,\xi,\eps)=p_{0,1}(x,\xi,\eps)+p_{0,2}(x,\xi,\eps)+\ldots \,+
p_{0,N}(x,\xi,\eps)+{\cal O}((\eps,\xi)^{N+1}),
\endeq
where $p_{0,j}$ is homogeneous of degree $j$ in $(\xi,\eps)$, so
that in particular
$$
p_{0,1}(x,\xi,\eps)=a\cdot \xi +i\eps \langle{q}\rangle(0)
$$
is independent of $x$.

We shall now remove the $x$-dependence also in the terms
$p_{0,j}$, $2\leq j\leq N$. In doing so, we consider the formal
power series
\begin{equation}
\label{2.1.5} G=\sum_{j=1}^{\infty} G_j(x,\xi,\eps),
\end{equation}
where $G_j={\cal O}((\xi,\eps)^{j+1})$ depends analytically on $x$
and is a homogeneous polynomial of degree $j+1$ in $(\xi,\eps)$.
In the sense of formal Taylor expansions, we then have
\begin{eqnarray*}
p_0\circ\exp(H_G)=p_0+\sum_{k=1}^{\infty} \frac{1}{k!} H_G^k p_0 &
= & p_0+\sum_{k=1}^{\infty}\sum_{l=1}^{\infty} \frac{1}{k!}H_G^k
p_{0,l} \\ \nonumber
 = p_0+\sum_{k=1}^{\infty}\sum_{l=1}^{\infty}\sum_{j_1=1}^{\infty}\ldots
\sum_{j_k=1}^{\infty} \frac{1}{k!} H_{G_{j_1}}\ldots
H_{G_{j_k}}p_{0,l} & = & \sum_{n=1}^{\infty} q_n.
\end{eqnarray*}
Here $q_n$ is a homogeneous polynomial of degree $n$ in
$(\xi,\eps)$. We also notice that in the sum above
$H_{G_m}p_{0,l}$ is homogeneous of degree $m+l$ in $(\xi,\eps)$.
It follows that $q_1=p_{0,1}$ is independent of $x$,
$q_2=p_{0,2}+H_{G_1}p_{0,1}$, and
$q_{n+1}=p_{0,n+1}+H_{G_n}p_{0,1}+\widetilde{q}_{n+1}$, $n\geq 2$,
where $\widetilde{q}_{n+1}$ depends only on $G_1,\ldots G_{n-1}$.
We are therefore able to determine $G_1,G_2\ldots\,G_{N-1}$
successively by solving the cohomological equations \begeq
\label{2.2} H_{p_{0,1}}G_n=a\cdot
\partial_x
G_n=p_{0,n+1}+\widetilde{q}_{n+1}-\langle{p_{0,n+1}+\widetilde{q}_{n+1}}\rangle,\quad
1\leq n\leq N-1.
\endeq
In this way we achieve that all the $q_n$, $n\leq N$, are
independent of $x$.

\Remark. It follows from the construction that the function $G_1$ is
independent of $\xi$ and is of the form $G_1(x,\xi,\eps)=\eps^2 h(x)$
for some analytic function $h$. Furthermore, it is clear that
$$
G_j={\cal O}_j(\eps^2),\quad 1\leq j \leq N-1.
$$

\Remark. In the construction above, we could also have used a slightly different argument and looked for
$G$ as a formal power series in $\eps$ alone,
$$
G\sim G_1(x,\xi)+\eps G_2(x,\xi)+\ldots\,
$$
Then composing the symbol $p_0$ with the canonical transformation
$\exp(\eps^2 H_G)$ and repeating the previous arguments, we would
have been able to determine each $G_j$ modulo ${\cal
O}(\xi^{N+1})$, and the final result would have been the same.

\vskip 4mm Summarizing the discussion so far, we get the following
result.

\begin{prop}
Let $p_0(x,\xi,\eps)=p_N(\xi)+i\eps\langle{q}\rangle(\xi)+{\cal
O}(\eps^2)+{\cal O}(\xi^{N+1})+\eps {\cal O}(\xi^{N})$ be an
analytic function defined near $\xi=0$ in $T^*{\bf T}^2$,
depending smoothly on $\eps \in {\rm neigh}(0,\real)$. Here $N\in
\nat$ is fixed but can be taken arbitrarily large. Assume that
$$
p_N(\xi)=a\cdot \xi+{\cal O}(\xi^2),
$$
where $a$ satisfies {\rm (\ref{0.13})}. Then we can find analytic
functions $G_1,G_2,\ldots\,,G_{N-1}$, with $G_j(x,\xi,\eps)$ being
a homogeneous polynomial of degree $j+1$ in $(\xi,\eps)$, such
that $G_j={\cal O}_j (\eps^2)$, $1\leq j\leq N-1$, $G_1$ is
independent of $\xi$, and such that if
\begeq \label{2.2.5}
G^{(N)}=G_1+G_2+\ldots +G_{N-1},
\endeq
then
$$
p_{0}\circ \exp(H_{G^{(N)}})=p^{(N)}(\xi,\eps)+r_{N+1}(x,\xi,\eps).
$$
Here
$$
p^{(N)}(\xi,\eps)=a\cdot \xi+i\eps\langle{q}\rangle(0)+{\cal O}((\xi,\eps)^2)
$$
is independent of $x$ and $r_{N+1}(x,\xi,\eps)={\cal
O}((\xi,\eps)^{N+1})$. Writing $p^{(N)}(\xi,\eps)=p(\xi)+i\eps
q(\xi,\eps)$, where $p$ is real, we have
$$
d_{\xi} p(0)\,\,\,\Re d_{\xi} q(0,0)\,\,\,\wrtext{are linearly
independent}.
$$
\end{prop}

In section 5, we shall see that we can quantize the holomorphic
canonical transformation \begeq \label{2.2.6}
\widehat{\kappa}:=\exp(H_{G^{(N)}})
\endeq
by means of an analytic elliptic Fourier integral operator in the
complex domain. In this section we shall proceed somewhat
formally, and carrying out the corresponding conjugation of
$P_{\eps}$, we may assume from now on that we are given an
$h$-pseudo\-dif\-fe\-ren\-ti\-al operator, still denoted by
$P_{\eps}$, defined microlocally near $\xi=0$ in $T^*{\bf T}^2$,
whose full symbol has a complete asymptotic expansion, \begeq
\label{2.2.7} P_{\eps}(x,\xi,\eps;h)=p_0+hp_1+h^2p_2+\ldots,
\endeq
with all $p_j=p_j(x,\xi,\eps)$ holomorphic in a fixed complex
\neigh{} of $\xi=0$, depending smoothly on $\eps\in {\rm
neigh}(0,\real)$, and such that
\begin{equation}
\label{2.3} p_0(x,\xi,\eps)=p^{(N)}(\xi,\eps)+{\cal
O}((\xi,\eps)^{N+1}),\quad p^{(N)}(\xi,\eps)=a\cdot
\xi+i\eps\langle{q}\rangle(0)+{\cal O}((\xi,\eps)^2).
\end{equation}
Our goal now is to make the lower order terms $p_j$, $j\geq 1$, in
(\ref{2.2.7}) independent of $x$, to a high order in $\xi$ and
$\eps$. This will be achieved by means of the usual conjugation of
$P_{\eps}$ by an elliptic pseudodifferential operator of the form
$\exp(Q(x,hD_x,\eps;h))$, where $Q$ is of order 0 in $h$. Write
$Q(x,hD_x,\eps;h)\sim q_0+hq_1+\ldots\,$. Then on the operator
level, we have
\begin{eqnarray*}
& & e^Q P_{\eps}
e^{-Q}=P_{\eps}+\sum_{k=1}^{\infty}\sum_{l=0}^{\infty}\frac{1}{k!}
h^l \left({\rm ad}Q\right)^k p_l \\ \nonumber & = &
P_{\eps}+\sum_{k=1}^{\infty}\sum_{l=0}^{\infty}\sum_{j_1=0}^{\infty}\ldots
\sum_{j_k=0}^{\infty} \frac{1}{k!} h^{l+j_1+\ldots
+j_k+k}\left(\frac{1}{h}{\rm ad}(q_{j_1})\right)\ldots
\left(\frac{1}{h}{\rm ad}(q_{j_k})\right)p_l,
\end{eqnarray*}
and on the symbol level we get
\begin{equation}
\label{2.4} e^Q P_{\eps} e^{-Q}\sim \sum_{n=0}^{\infty} h^n s_n.
\end{equation}
Here $s_0=p_0$, $s_1=p_1+iH_{p_0}q_0$, and for $n\geq 1$, the term
$s_{n+1}$ has the form
$s_{n+1}=p_{n+1}+iH_{p_0}q_n+\widetilde{s}_{n+1}$, where
$\widetilde{s}_{n+1}$ depends only on $q_0,\ldots q_{n-1}$. We
wish to make the terms $s_n$, $n\leq N$, independent of $x$,
modulo ${\cal O}((\xi,\eps)^{N+1})$. Taking a Taylor expansion of
$q_0$ in $(\xi,\eps)$, and solving the cohomological equations as
before, we see that we can choose $q_0$ so that $s_1$ becomes
independent of $x$ modulo ${\cal O}((\xi,\eps)^{N+1})$. Repeating
this argument, we determine successively $q_1,q_2,\ldots$ so that
$$
s_n=s_n^{(N)}+{\cal O}((\xi,\eps)^{N+1}),\quad n=0,1,2,\ldots,
$$
with $s_n^{(N)}=s_n^{(N)}(\xi,\eps)$ independent of $x$.

With $Q^{(N)}=\sum_{n=0}^{N-1} h^n q_n$, we get
\begin{equation}
\label{2.5} e^{Q^{(N)}}P_{\eps}
e^{-Q^{(N)}}=P^{(N)}(hD_x,\eps;h)+R_{N+1}(x,hD_x,\eps;h),
\end{equation}
where the full symbol of $P^{(N)}$ is independent of $x$ and the
symbol of $R_{N+1}$ is ${\cal O}((h,\xi,\eps)^{N+1})$. Moreover,
on the symbol level we have
\begin{equation}
\label{2.6} P^{(N+1)}(\xi,\eps;h)-P^{(N)}(\xi,\eps;h)={\cal
O}((\xi,\eps,h)^{N+1}).
\end{equation}
Using this compati\-bi\-lity property, we introduce a limi\-ting
ope\-ra\-tor
\begin{equation}
\label{C}
P^{(\infty)}=P^{(\infty)}(hD_x,\eps;h),
\end{equation}
such that
$$
P^{(\infty)}(\xi,\eps;h)=P^{(N)}(\xi,\eps;h)+{\cal
O}((\xi,\eps,h)^{N+1}),
$$
for all $N$. Then $P^{(\infty)}(\xi,\eps;h)$ is well-defined
modulo ${\cal O}((\xi,\eps,h)^{\infty})$.

We summarize the discussion of this section in the following
proposition.

\begin{prop}
Let
$$
P_{\eps}\sim p_0+hp_1+\ldots\,,\quad \abs{\xi}\leq \frac{1}{{\cal
O}(1)},
$$
be such that
$$
p_0=p_{N}(\xi)+i\eps \langle{q}\rangle(\xi)+{\cal O}(\eps^2)+{\cal
O}(\xi^{N+1})+\eps {\cal O}(\xi^{N}),
$$
for some fixed integer $N$ that can be taken arbitrarily large.
Then there exist \begeq \label{2.7}
G^{(N)}(x,\xi,\eps)=\sum_{j=1}^{N-1} G_j(x,\xi,\eps),
\endeq
with $G_j={\cal O}_j(\eps^2)$, $1\leq j\leq N-1$, being homogeneous
of degree $j+1$ in $(\xi,\eps)$ and depending analy\-ti\-cally on
$x$, and \begeq \label{2.8} Q^{(N)}(x,\xi,\eps;h)=\sum_{j=0}^{N-1}
h^j q_j(x,\xi,\eps),
\endeq
where $q_j$ are analytic, such that on the operator level,
$e^{{\rm ad}Q^{(N)}}e^{\frac{i}{h}{\rm ad}G^{(N)}}P_{\eps}$ is of
the form \begeq \label{2.9}
P^{(N)}(hD_x,\eps;h)+R_{N+1}(x,hD_x,\eps;h),
\endeq
where the full symbol of $P^{(N)}(hD_x,\eps;h)$ is independent of
$x$ and
$$
R_{N+1}(x,\xi,\eps;h)={\cal O}((h,\xi,\eps)^{N+1}).
$$
The leading symbol of $P^{(N)}(hD_x,\eps;h)$ is
$p_0^{(N)}(\xi,\eps)=p^{(N)}(\xi,\eps)=a\cdot
\xi+i\eps\langle{q}\rangle(0)+{\cal O}((\xi,\eps)^2)=p(\xi)+i\eps
q(\xi,\eps)$, where
\begin{equation}
\label{2.10} d_{\xi} p(0),\,\, \Re d_{\xi}q(0,0)\quad \wrtext{are
linearly independent}.
\end{equation}
\end{prop}

\Remark. Recalling the operator $P^{\infty}(hD_x,\eps;h)$ from
(\ref{C}), we notice that formally, Proposition 3.2 leads to the
quasi-eigenvalues \begeq \label{2.11}
P^{(\infty)}\left(h\left(k-\frac{k_0}{4}\right)-\frac{S}{2\pi},\eps;h\right)+{\cal
O}(h^{\infty}),\quad k\in \z^2,
\endeq
provided that $\eps={\cal O}(h^{\delta})$ and $kh={\cal
O}(h^{\delta})$, for some fixed $\delta>0$. Here we have also used
that the space $L^2_{\theta}({\bf T}^2)$, equipped with the
$L^2$-norm over a fundamental domain of ${\bf T}^2$ has an
orthonormal basis given by \begeq \label{2.12}
e_k(x)=e^{\frac{i}{h}x\cdot\left(h(k-\frac{k_0}{4})-\frac{S}{2\pi}\right)},\quad
k\in \z^2.
\endeq

\Remark. The pseudodifferential Birkhoff normal form construction
carried out in Proposition 3.1 and Proposition 3.2 could also have
been done in one stroke, by expanding the full symbol of
$P_{\eps}$ in a sum of homogeneous terms in all the variables
$(h,\xi,\eps)$, and making a similar decomposition for the
operator $G^{(N)}$. We hope to be able to return to this idea in a
future paper.

\section{The Hamilton-Jacobi equation}
\label{sectionHJ}\setcounter{equation}{0}

The results of this section will be instrumental in carrying out a
complete spectral analysis of the operator $P_{\eps}$, in the case
when $\eps$ is sufficiently small but fixed. Here we shall do the
work that will allow us to construct an additional complex
canonical transformation which, together with the normal form
construction of section 3, will allow us to reduce the leading
symbol of $P_{\eps}$ to a function on $T^*{\bf T}^2$, which is
independent of the spatial variables. The construction will be
somewhat similar to the corresponding constructions in~\cite{MeSj}
and~\cite{Sj2004}, and as there, the basic idea is to work with
cohomological equations of $\overline{\partial}$--type, liberating
ourselves from Diophantine conditions.

Let us recall from Proposition 3.1 that on the principal symbol
level, we can reduce ourselves to the case of a symbol \begeq
\label{hj1} p(x,\xi,\eps)=p^{(N-1)}(\xi,\eps)+r_N(x,\xi,\eps)
\endeq
in a complex \neigh{} of $\xi=0$ in $T^*{\bf T}^2$. Here
$p^{(N-1)}(\xi,\eps)$ is a polynomial of degree at most $N-1$ in
$(\xi,\eps)$, and $r_N(x,\xi,\eps)={\cal O}((\xi,\eps)^N)$. Assume
more generally that $r=r_N={\cal
O}\left(\eps^M+\abs{\xi}^N\right)$, for $2\leq M\leq N$. Here
$M,N\in \nat$ can be arbitrarily large but fixed.

Write
\begeq
\label{hj2}
p^{(N-1)}(\xi,\eps)=p(\xi)+i\eps q(\xi,\eps),
\endeq
and recall from Proposition 3.1 that ($p$ is real and that)
\begeq
\label{hj3} d_{\xi}p(0),\,\,d_{\xi}\Re q(0,0)\quad \wrtext{are
linearly independent}.
\endeq

Let $\widetilde{\eps}>0$ be a small parameter and assume that
$\xi={\cal O}(\widetilde{\eps})$. We put \begeq \label{hj4}
z(\xi,\eps)=p(\xi)+i\eps q(\xi,\eps).
\endeq
We shall try to find a grad-periodic function
$\psi(x)=\psi(x,\xi)$ with $\psi'(x)={\cal O}(\widetilde{\eps})$,
so that $\varphi(x,\xi)=x\cdot \xi+\psi(x)$ solves the
Hamilton-Jacobi equation \begeq \label{hj5}
p(x,\varphi'(x),\eps)-z(\xi,\eps)=0,
\endeq
so that \begeq \label{hj6} p(\xi+\psi'_x(x))+i\eps
q(\xi+\psi'_x(x),\eps)+r_N(x,\xi+\psi'_x(x),\eps)=z(\xi,\eps).
\endeq
Here we ``linearize'' the first two terms and get \begeq
\label{hj7} \left(p'(\xi)+i\eps q'_{\xi}(\xi,\eps)\right)\cdot
\partial_x
\psi+r_N(x,\xi+\psi_x'(x),\eps)+A(\xi,\psi_x'(x),\eps)\psi_x'(x)\cdot
\psi_x'(x)=0,
\endeq
where \begeq \label{hj8} A(\xi,\eta,\eps)=\int_0^1 (1-t)
\left(p+i\eps q\right)''(\xi+t\eta)\,dt.
\endeq

In order to apply some standard results about non-linear functions of
Sobolev class functions, we write
\begeq
\label{hj9}
\psi=\widetilde{\eps}\rho,
\endeq
and (\ref{hj7}) becomes \begeq \label{hj10} \left(p'(\xi)+i\eps
q'_{\xi}(\xi;\eps)\right)\cdot \partial_x
\rho+G(x,\xi,\rho'_x,\eps,\widetilde{\eps})=0,
\endeq
where \begeq \label{hj11} G(x,\xi,\eta,\eps,\widetilde{\eps})=
\frac{1}{\widetilde{\eps}}r_N(x,\xi+\widetilde{\eps}\eta;\eps)+\widetilde{\eps}A(\xi,\widetilde{\eps}\eta,\eps)\eta\cdot
\eta.
\endeq

Let $y_j(x)=y_j(x,\xi;\eps)$, $j=1,2$, be linear functions of
$x\in \real^2$, with \begeq \label{hj12} l_j(\xi,\eps)\cdot
\partial_x y_k=\delta_{j,k},
\endeq
where $l_1=p_{\xi}'(\xi)$ and $l_2=q'_{\xi}(\xi,\eps)$. Put \begeq \label{hj13}
w^*=\frac{1}{2}\left(\eps y_1+\frac{1}{i}y_2\right),
\endeq
and
\begeq
\label{hj14}
w=\frac{1}{2}\left(\eps y_1-\frac{1}{i}y_2\right),
\endeq
so that
\begin{eqnarray}
\label{hj15} \left(p'(\xi)+i\eps q'_{\xi}(\xi,\eps)\right)\cdot
\partial_x w^*=\eps,\\ \nonumber \left(p'(\xi)+i\eps
q'_{\xi}(\xi,\eps)\right)\cdot \partial_x w=0.
\end{eqnarray}
We now look for grad-periodic solutions to (\ref{hj10}) of the
form \begeq \label{hj16} \rho=\rho_{{\rm
per}}(x,\xi,\eps,\widetilde{\eps})+b(\xi,\eps,\widetilde{\eps})w^*(x,\xi,\eps),
\endeq
where $\rho_{{\rm per}}$ is single-valued, and try to find $\rho$ as a
limit $\rho=\lim_{j\rightarrow \infty}\rho_j$, where
\begeq
\label{hj17}
\rho_0=0,
\endeq
\begeq \label{hj18} \left(p'_{\xi}(\xi)+i\eps
q'_{\xi}(\xi,\eps)\right)\cdot \partial_x \rho_j+
G\left(x,\xi,(\rho_{j-1})'_x,\eps,\widetilde{\eps}\right)=0,\quad
j\geq 1.
\endeq
We shall estimate the gradient of the solutions to (\ref{hj18}) in
the standard $H^s$-norms on ${\bf T}^2$ for $s>1$, and observe
that \begeq \label{hj19} \norm{\rho'_x}_{H^s}\sim
\norm{(\rho_{{\rm per}})'_x}_{H^s}+\abs{b},
\endeq
uniformly with respect to $\eps$, for $\rho$ of the form
(\ref{hj16}). We shall also use the fact that if
$$
\left(p'_{\xi}+i\eps q'_{\xi}\right)\cdot \partial_x u=v,
$$
with $u$, $v$ periodic, then
$$
\norm{u'_x}_{H^s}\leq \frac{{\cal O}(1)}{\eps}\norm{v}_{H^s}.
$$

We notice now that for a given function $\rho_{j-1}$ of the form
(\ref{hj16}), (\ref{hj18}) is solvable with a solution of the same
form, unique up to a constant, and that this equation can be
decomposed into the following two equations \begeq \label{hj20}
\eps b_j+{\cal
F}\left(G(\cdot,\xi,(\rho_{j-1})'_x,\eps,\widetilde{\eps})\right)(0)=0,
\endeq
and \begeq \label{hj21} \left(p'(\xi)+i\eps
q'_{\xi}(\xi,\eps)\right)\cdot \partial_x (\rho_{{\rm
per},j})+G(x,\xi,(\rho_{j-1})'_x,\eps,\widetilde{\eps})-{\cal
F}\left(G(\cdot,\xi,(\rho_{j-1})'_x,\eps,\widetilde{\eps})\right)(0)=0.
\endeq
Here we let ${\cal F}u(k)$ stand for the Fourier coefficient at $k\in
\z^2$ of the $(2\pi \z)^2$-periodic function $u$.

In order to treat (\ref{hj18}) with $j=1$, we notice that \begeq
\label{hj22}
\norm{G(\cdot,\xi,0,\eps,\widetilde{\eps})}_{H^s}={\cal
O}(1)\frac{1}{\widetilde{\eps}}
\left(\eps^M+\widetilde{\eps}^N\right).
\endeq
For $\rho_1$, we then get
\begeq
\label{hj23}
\norm{\rho_1'}_{H^s}\leq {\cal
O}(1)\frac{\eps^M+\widetilde{\eps}^N}{\eps\widetilde{\eps}}.
\endeq
We may remark that for some component of the gradient, we get a better
estimate:
$$
\norm{p'(\xi)\cdot \partial_x \rho_1}_{H^s}\leq {\cal
O}(1)\frac{\eps^M+\widetilde{\eps}^N}{\widetilde{\eps}}.
$$

For $j\geq 2$, we consider in general for $u'_x$, $v'_x={\cal O}(1)$ in
$H^s$, for some $s>1$,
\begin{eqnarray*}
& &
G(x,\xi,u'_x,\eps,\widetilde{\eps})-G(x,\xi,v'_x,\eps,\widetilde{\eps})
\\
& = &
\frac{1}{\widetilde{\eps}}\left(r_N(x,\xi+\widetilde{\eps}u'_x,\eps)-r_N(x,\xi+\widetilde{\eps}v'_x,\eps)\right)
+ [\widetilde{\eps} A(\xi,\widetilde{\eps}\eta,\eps)\eta\cdot
  \eta]_{\eta=u'_x}^{v'_x}={\rm I}+{\rm II},
\end{eqnarray*}
with the natural definitions of ${\rm I}$ and ${\rm II}$. Writing
$$
{\rm I}=\int_0^1 (\partial_{\eta}
r_N)(x,\xi+\widetilde{\eps}(tu'_x+(1-t)v'_x),\eps)\,dt\cdot
(u'_x-v'_x),
$$
we get
$$
\norm{{\rm I}}_{H^s}\leq {\cal
O}(1)\frac{1}{\widetilde{\eps}}(\eps^M+\widetilde{\eps}^N)\norm{u'_x-v'_x}_{H^s},
$$
since by Proposition 2.2 of Chapter 2 in~\cite{AlGe}, which
extends to the case of vector-valued functions,
$$
\norm{\partial_{\eta}r_N(x,\xi+\widetilde{\eps}(tu'_x+(1-t)v'_x),\eps)}_{H^s}\leq
{\cal O}(1) \frac{1}{\widetilde{\eps}}(\eps^M+\widetilde{\eps}^N).
$$
Similarly, still assuming that $u'_x$, $v'_x={\cal O}(1)$ in
$H^s$, we get
$$
\norm{{\rm II}}_{H^s}\leq {\cal
O}(1)\widetilde{\eps}(\norm{u'_x}_{H^s}+\norm{v'_x}_{H^s})\norm{u'_x-v'_x}_{H^s}.
$$
Consequently,
\begin{eqnarray}
\label{hj24} & &
\norm{G(\cdot,\xi,u'_x,\eps,\widetilde{\eps})-G(\cdot,\xi,v'_x,\eps,\widetilde{\eps})}_{H^s}
\\ \nonumber
& \leq & {\cal
O}(1)\left(\frac{1}{\widetilde{\eps}}(\eps^M+\widetilde{\eps}^N)+
\widetilde{\eps}(\norm{u'_x}_{H^s}+\norm{v'_x}_{H^s})\right)\norm{u'_x-v'_x}_{H^s}.
\end{eqnarray}

We now make the a priori assumption
\begeq
\label{hj25}
\norm{\rho_j'}_{H^s}={\cal O}(1),\quad j\geq 0,
\endeq
uniformly in $j$, $\eps$, and notice that it is valid for
$j=0,\,1$, by (\ref{hj23}), if $\eps^{M-1}/\widetilde{\eps}$,
$\widetilde{\eps}^{N-1}/\eps\leq {\cal O}(1)$. Take the difference
between the equations (\ref{hj18}) for $j+1$ and $j$: \begeq
\label{hj26} (p'(\xi)+i\eps q'_{\xi}(\xi;\eps))\cdot
\partial_x(\rho_{j+1}-\rho_j)+
(G(x,\xi,\rho'_j,\eps,\widetilde{\eps})-G(x,\xi,\rho'_{j-1},\eps,\widetilde{\eps}))=0.
\endeq
In view of (\ref{hj24}) we get
\begeq
\label{hj27}
\norm{\rho'_{j+1}-\rho_j'}_{H^s}\leq {\cal
O}(1)\left(\frac{\eps^M+\widetilde{\eps}^N}{\eps\widetilde{\eps}}+
\frac{\widetilde{\eps}}{\eps}(\norm{\rho'_j}_{H^s}+\norm{\rho'_{j-1}}_{H^s})\right)\norm{\rho_j'-\rho_{j-1}'}_{H^s}.
\endeq

Choose $\widetilde{\eps}$ with
\begeq
\label{hj28}
\frac{\eps^M+\widetilde{\eps}^N}{\eps{\rm
min}(\eps,\widetilde{\eps})}\ll 1,
\endeq
and notice that this requires that $M>2$. Then (\ref{hj27})
implies that \begeq \label{hj29}
\norm{\rho'_{j+1}-\rho_{j}'}_{H^s}\leq
\frac{1}{2}\norm{\rho'_j-\rho_{j-1}'}_{H^s},
\endeq
as long as
$$
\norm{\rho_j'},\,\,\norm{\rho_{j-1}'}\leq {\cal
O}(1)\frac{\eps^M+\widetilde{\eps}^N}{\widetilde{\eps}\eps},
$$
uniformly in $j$. This holds for $j=1$, and as long as we have
(\ref{hj29}), we can extend it to the next $j$-value. Then
(\ref{hj27}) implies that we have a convergent sequence:
\begeq
\label{hj30}
\norm{\rho'_{j+1}-\rho'_j}_{H^s}\leq \left({\cal
O}(1)\frac{\eps^M+\widetilde{\eps}^N}{\eps{\rm
min}(\eps,\widetilde{\eps})}\right)^{j-1}
\frac{\eps^M+\widetilde{\eps}^N}{\eps\widetilde{\eps}}.
\endeq

Summing up, we have

\begin{prop}
Assume that $\eps>0$, $\widetilde{\eps}>0$ small, are such that
\begeq \label{hj31} \frac{\eps^M+\widetilde{\eps}^N}{\eps {\rm
min}(\eps,\widetilde{\eps})}\ll 1.
\endeq
Then for $\xi={\cal O}(\widetilde{\eps})$, the Hamilton-Jacobi
equation {\rm (\ref{hj5})} has a solution $\varphi(x,\xi)=x\cdot
\xi+\widetilde{\eps}\rho(x,\xi)$, with $\rho$ of the form {\rm
(\ref{hj16})}, and
$$
\norm{\rho'_x}_{H^s}={\cal
O}((\eps^M+\widetilde{\eps}^N)/(\eps\widetilde{\eps})),\quad s>1,
$$
\end{prop}

Essentially, the same iteration scheme as above shows that the
solution to (\ref{hj5}) is unique up to a constant, if we require
that $\frac{\widetilde{\eps}}{\eps}\norm{\rho'_x}_{H^s} \ll 1$ in
addition to (\ref{hj31}).

\section{Global Grushin problem for $\eps={\cal O}(h^{\delta})$}
\label{section4} \setcounter{equation}{0}

Throughout this section, it will be assumed that $\eps={\cal
O}(h^{\delta})$ for some $\delta>0$, fixed, but arbitrarily small.
Our goal is to prove Theorem 1.1 and show that the
quasi-eigenvalues introduced in (\ref{2.11}) give all the
eigenvalues of $P_{\eps}$, modulo ${\cal O}(h^{\infty})$, in a
rectangle (\ref{R0}). In doing so, we shall only rely on the
results of sections 2 and 3.

Our starting point in this section is the operator $P_{\eps}$ with
a leading symbol (\ref{1.34}). We shall first discuss how to
implement the conjugation of $P_{\eps}$ by the analytic Fourier
integral operator $e^{-iG^{(N)}/h}$, introduced in Proposition
3.2, which quantizes the complex canonical transformation
$\widehat{\kappa}=\exp(H_{G^{(N)}})$, defined in (\ref{2.2.5}),
(\ref{2.2.6}). In doing so, let us consider the IR-manifold,
defined in a complex \neigh{} of $\xi=0$ and given by
$\exp(H_{G^{(N)}})(T^*{\bf T}^2) \subset {\widetilde{{\bf
T}}^2}\times \comp^2$. Here $\widetilde{{\bf T}}^2={\bf
T}^2+i\real^2$ stands for the standard comple\-xi\-fi\-cation of
${\bf T}^2$. It follows from Proposition 3.2 that along this
manifold we have \begeq \label{4.0} \Im x={\cal O}(\eps^2),\,\,
\Im \xi={\cal O}(\eps^2).
\endeq

Let us introduce the usual FBI-Bargmann transform
$$
Tu(x)=Ch^{-3/2}\int e^{i\varphi(x,y)/h}u(y)\,dy,\quad C>0,
$$
acting on $L^2_{\theta}({\bf T}^2)$. Here the integration is
performed over the whole of $\real^2$, with $u$ being Floquet
periodic---see also the discussion in section 3 in~\cite{MeSj}.
The phase function $\varphi(x,y)=\frac{i}{2}(x-y)^2$ is defined
for $x\in {\bf T}^2+i\real^2$, and the associated canonical
transformation $\kappa_T$ is given by \begeq \label{4.0.1}
T^*\widetilde{{\bf T}}^2 \ni (y,\eta)\mapsto
(x,\xi)=(y-i\eta,\eta)\in T^*\widetilde{{\bf T}}^2.
\endeq
The transform $\kappa_T$ maps the real phase space $T^*{\bf T}^2$ to the IR-manifold
$$
\Lambda_{\Phi_0}:\,\,\xi=\frac{2}{i}\frac{\partial
\Phi_0}{\partial x}=-\Im x,\quad \Phi_0(x)=\frac{1}{2}\left(\Im
x\right)^2,
$$
contained inside $T^*\widetilde{{\bf T}}^2$. Here the zero section
in $T^*{\bf T}^2$ corresponds to $\Im x=0$. On the transform side,
the manifold $\exp(H_{G^{(N)}})\left(T^*{\bf T}^2\right)$ is
represented by
$$
\xi=\frac{2}{i}\frac{\partial \Phi_{\eps}}{\partial x}, \quad \abs{\Im
x}\leq \frac{1}{{\cal O}(1)},
$$
where $\Phi_{\eps}$ is a smooth strictly plurisubharmonic
function, such that
$$
\nabla \left(\Phi_{\eps}-\Phi_0\right)={\cal O}(\eps^2).
$$
By choosing the undetermined constant in $\Phi_{\eps}$ suitably,
we can even arrange that \begeq \label{4.0.2}
\Phi_{\eps}-\Phi_0={\cal O}(\eps^2).
\endeq

Let now $\chi=\chi(\Im x)\in C^{\infty}_0$, $0\leq \chi \leq 1$,
be a standard cutoff function in a \neigh{} of 0, and consider the
IR-manifold given by
$$
\Lambda_{\widetilde{\Phi}_{\eps}}:\quad \xi=\frac{2}{i}\frac{\partial
\widetilde{\Phi}_{\eps}}{\partial x},
$$
where
$$
\widetilde{\Phi}_{\eps}(x)=\chi(\Im x)\Phi_{\eps}(x)+\left(1-\chi(\Im x)\right)\Phi_0(x).
$$
Here $\widetilde{\Phi}_{\eps}$ is strictly plurisubharmonic, and
if we introduce the globally defined IR-manifold \begeq
\label{4.0.3}
\Lambda_{\eps}:=\kappa_T^{-1}(\Lambda_{\widetilde{\Phi}_{\eps}})\subset
\widetilde{{\bf T}}^2\times \comp^2,
\endeq
then $\Lambda_{\eps}$ agrees with $\exp(H_{G^{(N)}})(T^*{\bf
T}^2)$ in a complex \neigh{} of the zero section and it is equal
to $T^*{\bf T}^2$ further away from this set. Moreover, since
along $\Lambda_{\widetilde{\Phi}_{\eps}}$,
$$
\Im \xi=-\frac{\partial \widetilde{\Phi}_{\eps}}{\partial {\rm Re}
x}={\cal O}(\eps^2),\quad \Re \xi=-\frac{\partial
\widetilde{\Phi}_{\eps}}{\partial {\rm Im} x}=-\Im x+{\cal
O}(\eps^2),
$$
we conclude that along $\Lambda_{\eps}$, we have \begeq
\label{3.1} \Im \xi={\cal O}\left(\eps^2\right),\quad \Im x={\cal
O}(\eps^2).
\endeq
Now the Fourier integral operator
$$
e^{-\frac{i}{h} G^{(N)}}={\cal O}(1):L^2_{\theta}({\bf
T}^2)\rightarrow H_{\theta}(\Lambda_{\eps})
$$
is such that the action of $P_{\eps}$ on the space
$H(\Lambda_{\eps})$, associated to $\Lambda_{\eps}$ by means of
the FBI transform $T$, is microlocally near
$\exp(H_{G^{(N)}})({\bf T}^2\times \{\xi=0\})$, unitarily
equivalent to the operator \begeq \label{3.1.1}
e^{\frac{i}{h}G^{(N)}}P_{\eps} e^{-\frac{i}{h}G^{(N)}}:
L^2_{\theta}({\bf T}^2)\rightarrow L^2_{\theta}({\bf T}^2).
\endeq
From Proposition 3.1 we know that the leading symbol of
(\ref{3.1.1}) is independent of $x$, modulo ${\cal
O}((\xi,\eps)^{N+1})$.

Let $0<\widetilde{\eps}\ll 1$ be an additional small parameter
such that $\widetilde{\eps} \gg {\rm max}(\eps,h)$. It is then
clear, in view of (\ref{2.10}), that along $T^*{\bf T}^2$, in a
region where $\widetilde{\eps}\leq \abs{\xi}\leq 1/{\cal O}(1)$,
we have \begeq \label{3.2} \abs{\Re P_{\eps}}\geq
\frac{\widetilde{\eps}}{{\cal O}(1)}\quad \wrtext{or}\quad
\abs{\Im P_{\eps}-\eps \Re F_0}\geq \frac{\eps
\widetilde{\eps}}{{\cal O}(1)}.
\endeq

We now notice that the fact that $P_{\eps=0}$ is selfadjoint in
$L^2_{\theta}({\bf T^2})$ implies that the symbol of
$$
\Im P_{\eps}=\frac{P_{\eps}-P_{\eps}^*}{2i},
$$
taken in the operator sense on $H(\Lambda_{\eps})$, is ${\cal
O}(\eps)+{\cal O}(\eps h)$. It follows then from the property
(\ref{3.1}) of the IR-deformation given by $\Lambda_{\eps}$ that
along $\Lambda_{\eps}$, outside any $\widetilde{\eps}$-\neigh{} of
the Lagrangian torus $\exp(H_{G_N})({\bf T}^2\times \{0\})$, the
estimates (\ref{3.2}) still hold true.

Corresponding to the manifold $\Lambda_{\eps}$ on the torus side,
we get a globally defined IR-manifold $\Lambda\subset
T^*\widetilde{M}$, which is $\eps$-close to $T^*M$ everywhere,
agrees with that set near infinity, and in a complex \neigh{} of
$\Lambda_1$, it is obtained by replacing $\exp(i\eps H_G)\circ
\kappa_1^{-1}\circ(\kappa^{(N)})^{-1}(T^*{\bf T}^2)$ there by
\begeq \label{3.2.1} \exp(i\eps H_G)\circ \kappa_1^{-1}\circ
(\kappa^{(N)})^{-1}\circ \exp(H_{G^{(N)}})(T^*{\bf
T}^2)=\exp(i\eps H_G)\circ \kappa_1^{-1}\circ
\kappa_0^{-1}(\Lambda_{\eps}).
\endeq
Here we recall that the real analytic canonical transformations
$\kappa_1$ and $\kappa^{(N)}$ have been defined in (\ref{1.1}) and
(\ref{1.3}), respectively. (We have also written here $\exp(i\eps
H_G)$ for the complex canonical transformation identifying
$\Lambda_{\eps G}$ and $T^*M$ in a \neigh{} of $\Lambda_1$.)

Using Propositions 2.3 and 2.5, and taking into account also the
conjugation by the analytic pseudo-differential operator
$\exp(Q^N(x,hD_x,\eps;h))$, with $Q^N$ defined in (\ref{2.8}), we
arrive at the following result.

\begin{prop}
Keep all the general assumptions from the introduction, and in
particular, {\rm (\ref{0.12})}. We write $F$ to denote the mean
value of $q$ over the Diophantine tori $\Lambda_j$, $1\leq j\leq
L$, and let us make the global dynamical assumption {\rm
(\ref{H0})}. In the case when the $H_p$-flow is completely
integrable, we make the assumption {\rm (\ref{1.31})}. Recall also
the real analytic canonical transformations $\kappa_j$, $1\leq
j\leq L$, and $\kappa^{(N)}$, where $\kappa_j$ is defined in {\rm
(\ref{1.1})} for $j=1$, and $\kappa^{(N)}$ is defined in {\rm
(\ref{1.3})}. Then the composed transform
$$
\kappa^{(N)}\circ \kappa_j: {\rm neigh}(\Lambda_j,T^*M)\rightarrow
{\rm neigh}(\xi=0,T^*{\bf T}^2)
$$
maps $\Lambda_j$ to $\xi=0$, and has the property that when
expressed in terms of the coordinates $x$ and $\xi$ on the torus
side, the leading symbol $p_{\eps}$ of $P_{\eps}$ becomes
$$
p_{N,j}(\xi)+i\eps q_j(x,\xi)+{\cal O}(\eps^2)+{\cal
O}(\xi^{N+1}),
$$
with $p_{N,j}(\xi)=a\cdot \xi+{\cal O}(\xi^2)$. Here $N\geq 1$ is
a fixed positive integer that can be taken arbitrarily large.
Using these $(x,\xi)$-coordinates, we then define for $\xi$ small,
$$
\langle{q_j}\rangle(\xi)=\frac{1}{(2\pi)^2}\int q_j(x,\xi)\,dx,
$$
and assume that $dp_{N,j}(0)=a_j=a$ and $\Re d\langle{q_j}\rangle(0)$
are linearly independent, $1\leq j\leq L$. Let
$0<\widetilde{\eps}\ll 1$ be such that $\widetilde{\eps}\gg {\rm
max}(\eps,h)$. Then there exists a globally defined IR-manifold
$\Lambda\subset T^*\widetilde{M}$ and $L$ smooth Lagrangian tori
$\widehat{\Lambda}_1,\ldots,\,\widehat{\Lambda}_L\subset \Lambda$,
such that when $\rho\in \Lambda$ is away from any
$\widetilde{\eps}$-\neigh{} of $\cup_{j=1}^L \widehat{\Lambda}_j$
in $\Lambda$, we have
$$
\abs{\Re P_{\eps}(\rho)}\geq \frac{\widetilde{\eps}}{{\cal
O}(1)}\quad \wrtext{or}\quad \abs{\Im P_{\eps}-\eps \Re F}\geq
\frac{\eps \widetilde{\eps}}{{\cal O}(1)}.
$$
The manifold $\Lambda$ is $\eps$-close to $T^*M$ and agrees with
it outside a compact set. We have
$$
P_{\eps}={\cal O}(1): H(\Lambda,m)\rightarrow H(\Lambda),
$$
For each $j$ with $1\leq j\leq L$, there exists an elliptic
Fourier integral operator
$$
U_j={\cal O}(1): H(\Lambda)\rightarrow L^2_{\theta}({\bf T}^2),
$$
such that microlocally near $\widehat{\Lambda}_j$, we have
$$
UP_{\eps}=\left(P^{(N)}_j(hD_x,\eps;h)+R_{N+1,j}(x,hD_x,\eps;h)\right)U.
$$
Here $P^{(N)}_j(hD_x,\eps;h)+R_{N+1,j}(x,hD_x,\eps;h)$ is defined
microlocally near $\xi=0$ in $T^*{\bf T}^2$, the full symbol of
$P^{(N)}_j$ is independent of $x$, and
$$
R_{N+1,j}(x,\xi,\eps;h)={\cal O}((h,\xi,\eps)^{N+1}).
$$
The leading symbol of $P^{(N)}_j(hD_x,\eps;h)$ has the form
$$
a\cdot \xi+i\eps F+{\cal O}((\eps,\xi)^2)=p_j(\xi)+i\eps
q_j(\xi,\eps),
$$
where, as before, $p_j$ is real on the real domain, and
$d_{\xi}p_j(0)$ and $\Re d_{\xi} q_j(0,0)$ are linearly independent.
We have
$$
\Im P^{(N)}_j(\xi,\eps;h)=\eps \Re q_j(\xi,\eps)+{\cal O}(\eps h).
$$
\end{prop}

In what follows we shall choose $\widetilde{\eps}\gg \eps$ so that
\begeq \label{3.3} h^{1/2-\delta}<\widetilde{\eps}={\cal
O}(h^{\delta}).
\endeq
As we shall see, later the lower bound on $\widetilde{\eps}$  will
have to be strengthened.

Our goal is to describe the spectrum of $P_{\eps}$ in a rectangle
of the form \begeq \label{3.4} \abs{\Re
z}<\frac{\widetilde{\eps}}{C},\quad \abs{\Im z-\eps \Re
F}<\frac{\eps \widetilde{\eps}}{C},
\endeq
when $C>0$ is a sufficiently large but fixed constant. To this
end, let us introduce the quasi-eigenvalues coming from
Proposition 5.1---see also (\ref{2.11}),
\begeq \label{3.4.0.1}
z(j,k):=P^{(N)}_j\left(h\left(k-\frac{k_j}{4}\right)-\frac{S_j}{2\pi},\eps;h\right)+{\cal
O}(h^{\delta(N+1)}),\quad 1\leq j\leq L,\quad k\in \z^2,
\endeq
with $h(k-k_j/4)-S_j/2\pi={\cal O}(\widetilde{\eps})$. In what
follows we shall work under the assumption that $\eps$ is bounded
from below by an arbitrary but fixed positive power of $h$,
$\eps\geq h^K$, $K\gg 1$ is fixed. Let us also remark that it
follows from Proposition 5.1 that the distance between two
neighboring quasi-eigenvalues is $\geq \eps h/{\cal O}(1)$, provided
that $N$ is large enough.

\medskip
When $z\in \comp$ is in the rectangle (\ref{3.4}), let us consider
the equation \begeq \label{3.5} (P_{\eps}-z)u=v,\quad u\in
H(\Lambda).
\endeq
At first, we shall derive an a priori bound for the part of $u$
concentrated away from the set $\cup_{j=1}^L \widehat{\Lambda}_j$.

In what follows we shall write that a function $a\in
C^{\infty}(\Lambda)$ (also depending on $h$) is in the symbol
class $S_{\widetilde{\eps}}(1)$ if uniformly on $\Lambda$,
$$
\nabla^{m} a={\cal O}_{\alpha}(\widetilde{\eps}^{-m}),\quad m \geq
0.
$$
We introduce then a smooth partition of unity on the manifold
$\Lambda$, \begeq \label{3.6} 1=\sum_{j=1}^L
\chi_j+\psi_{1,+}+\psi_{1,-}+\psi_{2,+}+\psi_{2,-}.
\endeq
Here $0\leq \chi_j\in C^{\infty}_0(\Lambda)\cap
S_{\widetilde{\eps}}(1)$ is a cut-off function to an
$\widetilde{\eps}$-\neigh{} of $\widehat{\Lambda}_j$, $1\leq j\leq
L$, such that on the operator level we have, \begeq \label{3.7}
[P_{\eps},\chi_j]={\cal O}(h^{(N+1)\delta}): H(\Lambda)\rightarrow
H(\Lambda).
\endeq
To construct $\chi_j$ satisfying (\ref{3.7}), we use Proposition
5.1 to pass to the torus model and take there a symbol of the form
$\chi_0(\xi/\widetilde{\eps})$, where $\chi_0\in
C^{\infty}_0(T^*{\bf T}^2)$ is supported in a small \neigh{} of
$\xi=0$ in $T^*{\bf T}^2$, and $\chi_0=1$ near $\xi=0$.
Conjugating $\chi_0(hD_x/\widetilde{\eps})$ by means of the
microlocal inverse of the operator $U_j$ of Proposition 5.1, we
obtain the cut-off function $\chi_j$ with the required properties.

The functions $0\leq \psi_{1,\pm}\in S_{\widetilde{\eps}}(1)$ in
(\ref{3.6}) are chosen so that $\pm \Re P_{\eps}\geq
\widetilde{\eps}/{\cal O}(1)$ in the support of $\psi_{1,\pm}$,
respectively. Finally, we have $\pm \left(\Im P_{\eps}-\eps \Re
F\right) \geq \eps \widetilde{\eps}/{\cal O}(1)$ in the support
of $\psi_{2,\pm}\in C^{\infty}_0(\Lambda)\cap
S_{\widetilde{\eps}}(1)$, and as in the case of $\chi_j$, we arrange
so that in the operator norm,
\begeq
\label{3.7.5}
A [P_{\eps},\psi_{2,\pm}]={\cal O}(h^{(N+1)\delta}),
\endeq
where $A$ is a microlocal cut-off to a region where $\abs{\Re P_{\eps}}<\widetilde{\eps}/{\cal O}(1)$.

When $u$ satisfies (\ref{3.5}), we shall prove that \begeq
\label{3.8} \norm{\left(1-\sum_{j=1}^L \chi_j\right)u}\leq
\frac{{\cal O}(1)}{\eps \widetilde{\eps}}\norm{v}+{\cal
O}(h^{(N+1)\delta-K-1})\norm{u},
\endeq
where we shall choose $N$ such that $(N+1)\delta-K-1\gg 1$, and
provided that the lower bound on $\widetilde{\eps}$ in
(\ref{3.3}) is suitably strengthened. Here the norms are taken
in $H(\Lambda)$. When establishing (\ref{3.8}), we shall first
prove that
\begeq
\label{3.9} \norm{\psi_{1,+}u}\leq
\frac{{\cal O}(1)}{\widetilde{\eps}}\norm{v}+{\cal
O}(h^{\infty})\norm{u},
\endeq
which is essentially an elliptic estimate. To prove (\ref{3.9}),
we may follow the method of~\cite{HiSj1}. (See also~\cite{Sj92}.)
When $M\in \nat$, we let
$$
\psi_{1,+}=:\psi_0 \prec \psi_1 \prec \ldots \prec \psi_M
$$
be a nested sequence of functions in $S_{\widetilde{\eps}}(1)$,
supported in a region where $\Re P_{\eps}\geq
\widetilde{\eps}/{{\cal O}(1)}$. Near the support of $\psi_j$ we
have $\Re(P_{\eps}-z)\geq \widetilde{\eps}/{\cal O}(1)$. Using
$h/\widetilde{\eps}^2$ as a new semiclassical parameter that appears
naturally on the operator level when reducing
$S_{\widetilde{\eps}}(1)$ to $S_1(1)$ by dilation, and
applying the sharp G\aa{}rding inequality we get \begin{eqnarray}
\label{3.10}
 \Re ((P_{\eps}-z)\psi_j u|\psi_j u) & \geq
& \left(\frac{\widetilde{\eps}}{{\cal O}(1)}-{\cal
O}(1)\frac{h}{\widetilde{\eps}^2}\right)\norm{\psi_j u}^2-{\cal
O}(h^{\infty})\norm{u}^2 \\ \nonumber & \geq &
\frac{\widetilde{\eps}}{{\cal O}(1)}\norm{\psi_j u}^2-{\cal
O}(h^{\infty})\norm{u}^2, \quad j=0,\ldots\, ,M.
\end{eqnarray}
provided that $h/\widetilde{\eps}^3\ll 1$, and we shall even
strengthen this further by assuming that
\begeq
\label{3.11}
\frac{h}{\widetilde{\eps}^3}\leq h^{\delta}.
\endeq
Now using that $\psi_j(1-\psi_{j+1})={\cal O}(h^{\infty})$ in the
operator norm, we see that the absolute value of the left hand
side in (\ref{3.10}) does not exceed
$$
{\cal O}(1) \norm{v}\norm{\psi_j u}+{\cal
O}\left(\frac{h}{\widetilde{\eps}^2}\right)\norm{\psi_{j+1}
u}^2+{\cal O}(h^{\infty})\norm{u}^2.
$$
Here we have also used that $[P_{\eps},\psi_j]={\cal
O}(h/\widetilde{\eps}^2)$. We get
$$
\frac{\widetilde{\eps}}{{\cal O}(1)}\norm{\psi_j u}^2 \leq {\cal
O}(1)\norm{v}\norm{\psi_j u}+{\cal
O}\left(\frac{h}{\widetilde{\eps}^2}\right)\norm{\psi_{j+1}u}^2+{\cal
O}(h^{\infty})\norm{u}^2,
$$
and therefore,
$$
\norm{\psi_j u}^2\leq \frac{{\cal
O}(1)}{\widetilde{\eps}^2}\norm{v}^2+{\cal
O}\left(\frac{h}{\widetilde{\eps}^3}\right)\norm{\psi_{j+1}u}^2+{\cal
O}(h^{\infty})\norm{u}^2.
$$
Taking into account (\ref{3.11}), we obtain
$$
\norm{\psi_j u}^2 \leq \frac{{\cal
O}(1)}{\widetilde{\eps}^2}\norm{v}^2+{\cal
O}(h^{\delta})\norm{\psi_{j+1}u}^2+{\cal O}(h^{\infty})\norm{u}^2,
$$
and combining these estimates for $j=0,1,\ldots\, ,M$, we get
$$
\norm{\psi_{1,+}u}^2\leq \frac{{\cal
O}(1)}{\widetilde{\eps}^2}\norm{v}^2+{\cal
O}_M(h^{M\delta})\norm{\psi_M u}^2+{\cal O}(h^{\infty})\norm{u}^2.
$$
The estimate (\ref{3.9}) follows, and by repeating this
argument, we get the same estimate also for $\psi_{1,-}u$. We next
notice that $\Im P_{\eps}={\cal O}(\eps)$ on $\Lambda$, and near
$\supp\,\psi_{2,+}$ we have $\Im P_{\eps}-\eps \Re F\geq \eps\widetilde{\eps}/{\cal
O}(1)$. An application of the sharp G\aa{}rding inequality as in
(\ref{3.10}) gives then for $z$ in (\ref{3.4}),
\begeq
\label{3.11.1}
\Im ((P_{\eps}-z)\psi_{2,+}u|\psi_{2,+}u)\geq
\frac{\eps\widetilde{\eps}}{{\cal O}(1)}\norm{\psi_{2,+}u}^2-{\cal
O}(h^{\infty})\norm{u}^2.
\endeq
The left hand side of (\ref{3.11.1}) does not exceed
$$
{\cal
O}(1)\norm{v}\,\norm{\psi_{2,+}u}+\norm{[P_{\eps},\psi_{2,+}]u}\,\norm{\psi_{2,+}u},
$$
and combining this with (\ref{3.11.1}) we get
\begeq
\label{3.11.2}
\norm{\psi_{2,+}u}\leq \frac{{\cal
O}(1)}{\widetilde{\eps}\eps}\norm{v}+\frac{{\cal
O}(1)}{\widetilde{\eps}\eps}\norm{[P_{\eps},\psi_{2,+}]u}+{\cal
O}(h^{\infty})\norm{u}.
\endeq
Now using (\ref{3.7.5}) together with (\ref{3.9}) and (\ref{3.11}) we get
$$
\norm{[P_{\eps},\psi_{2,+}]u}\leq {\cal O}(h^{\delta})\norm{v}+{\cal
O}(h^{(N+1)\delta})\norm{u},
$$
and we conclude from (\ref{3.11.2}) that
$$
\norm{\psi_{2,\pm}u}\leq \frac{{\cal
O}(1)}{\eps\widetilde{\eps}}\norm{v}+{\cal
O}(h^{(N+1)\delta-K-1})\norm{u}.
$$
The estimate (\ref{3.8}) follows.

\medskip
In the spirit of~\cite{HiSj1}, as a warm-up exercise, we shall now
prove that $z$ in the rectangle (\ref{3.4}) avoids the union of
$\eps h/{\cal O}(1)$-neighborhoods of the quasi-eigenvalues
(\ref{3.4.0.1}), then $P_{\eps}-z: H(\Lambda,m)\rightarrow
H(\Lambda)$ is invertible. When doing so, we write, for $1\leq
j\leq L$,
$$
(P_{\eps}-z)\chi_j u=\chi_j v+[P_{\eps},\chi_j]u,
$$
where the $H(\Lambda)$-norm of the commutator term in the right
hand side is ${\cal O}(h^{(N+1)\delta})\norm{u}$. Applying the
operator $U_j$ of Proposition 5.1 we get
\begeq \label{3.12}
(P^{(N)}_j+R_{N+1,j}-z)U_j\chi_j u=U_j\chi_j v+T_{N,j}u,
\endeq
where
$$
T_{N,j}={\cal O}(h^{(N+1)\delta}): H(\Lambda)\rightarrow
L^2_{\theta}({\bf T}^2).
$$
Using the fact that $R_{N+1,j}(x,\xi,\eps;h)={\cal
O}((h,\eps,\xi)^{N+1})$ together with the localization properties
of $U_j\chi_j$, we can rewrite (\ref{3.12}) as
$$
(P^{(N)}_j-z)U_j\chi_j u=U_j \chi_j v+T_{1,N,j}u,
$$
where $T_{1,N,j}$ has the same bound as $T_{N,j}$. Taking an expansion
in Fourier series,
$$
v(x)=\sum_{k\in {\rm \bf Z^2}} \widehat{v}(k-\theta) e_k(x),\quad
v\in L^2_{\theta},
$$
where $e_l$ are defined in (\ref{2.12}), we next see that the
operator $P_j^{(N)}-z$ acting on $L^2_{\theta}({\bf T}^2)$ is
invertible, microlocally near $\xi=0$, with a microlocal inverse
of the norm ${\cal O}(1/\eps h)$, provided that $z$ in (\ref{3.4})
avoids the union of the $\eps h/{\cal O}(1)$-neighborhoods of the
quasi-eigenvalues (\ref{3.4.0.1}). It follows that for $1\leq j\leq
N$,
$$
\norm{\chi_j u}\leq \frac{{\cal O}(1)}{\eps h}\norm{v}+{\cal
O}(h^{(N+1)\delta-K-1})\norm{u},\,\, \eps \geq h^K,\,\,
(N+1)\delta-K-1\gg 1.
$$
Combining this estimate with (\ref{3.8}), we see that the
operator $P_{\eps}-z: H(\Lambda,m)\rightarrow H(\Lambda)$ is
injective, hence bijective, since general arguments based on the
ellipticity at infinity (\ref{0.6}), (\ref{0.9}) show that it
is a Fredholm operator of index zero.

\vskip 4mm We shall now discuss the setup of the global Grushin
problem for $P_{\eps}-z$, in the Hilbert space $H(\Lambda)$, which
will be well-posed for $z$ varying in the rectangle (\ref{3.4}).
When doing so, we will continue to assume that $\eps\geq h^K$, when $K\gg 1$ is fixed.

With $z(j,k)$, being defined in (\ref{3.4.0.1}), let us introduce
for $1\leq j\leq L$,
$$
M_j=\# \left\{z(j,k); \abs{\Re
z(k)}<\frac{\widetilde{\eps}}{C},\quad \abs{\Im z-\eps \Re
F_0}<\frac{\widetilde{\eps}{\eps}}{C} \right\},
$$
where $C$ is sufficiently large. Then $M_j={\cal
O}(\widetilde{\eps}^2h^{-2})$, and we let $k(j,1),\ldots
k(j,M_j)\in \z^2$ be the corresponding lattice points, so that
$$
h\left(k(j,l)-\frac{k_0}{4}\right)-\frac{S}{2\pi}={\cal
O}(\widetilde{\eps}),\quad 1\leq l\leq M_j,\,\,1\leq j\leq L.
$$
We introduce the operator
$$
R_+: H(\Lambda)\rightarrow \comp^{M_1}\times \cdots \times
\comp^{M_L},
$$
given by
$$
R_+u(j)(l)=(U_j\chi_j u|e_{k(j,l)}),\quad 1\leq j\leq L, \quad
1\leq l\leq M_j.
$$
Here $e_{k(j,l)}$ is as in (\ref{2.12}) and the scalar product in
the definition of $R_+$ is taken in $L^2_{\theta}({\bf T}^2)$.
Define next
$$
R_-: \comp^{M_1}\times \cdots \times \comp^{M_L} \rightarrow
H(\Lambda)
$$
by
$$
R_-u_-=\sum_{j=1}^L \sum_{l=1}^{M_j} u_-(j)(l)U_j^{-1}e_{k(j,l)}.
$$
Here $U_j^{-1}$ is a microlocal inverse of $U_j$. Let us recall
that $e_k(x)$ is microlocally concentrated in the region of
$T^*{\bf T}^2$ where $\xi\sim h(k-\frac{k_0}{4})-\frac{S}{2\pi}$.
It follows therefore that for $1\leq j\leq L$, \begeq \label{3.13}
\chi_j R_-=\sum_{l=1}^{M_j} u_-(j)(l)U_j^{-1} e_{k(j,l)}+{\cal
O}(h^{\infty}): \comp^{M_1}\times \cdots \times
\comp^{M_L}\rightarrow H(\Lambda).
\endeq

We shall next check that when $z\in \comp$ varies in the rectangle
(\ref{3.4}), with an increased value of $C$, the Grushin problem
\begin{eqnarray}
\label{3.14}{\cases{(P_{\eps}-z)u+R_-u_-=v, \cr R_+u=v_+ \cr}}
\end{eqnarray}
has a unique solution $(u,u_-)\in H(\Lambda,m)\times
\left(\comp^{M_1}\times \cdots \times \comp^{M_L}\right)$ for
every $(v,v_+)\in H(\Lambda)\times \left(\comp^{M_1}\times \cdots
\times \comp^{M_L}\right)$. Moreover, we have an a priori estimate
\begeq \label{3.15} \norm{u}+\norm{u_-}\leq \frac{{\cal
O}(1)}{\eps \widetilde{\eps}}\left(\norm{v}+{\norm{v_+}}\right),
\endeq
where the norms of $u$ and $v$ are taken in $H(\Lambda)$, and the
norms of $u_-$ and $v_+$ are taken in $\comp^{M_1}\times \cdots
\times \comp^{M_L}$. Indeed, we first notice that in view of
(\ref{3.8}) and (\ref{3.13}),
\begeq \label{3.16}
\norm{\left(1-\sum_{j=1}^L\chi_j\right)u}\leq \frac{{\cal
O}(1)}{\eps \widetilde{\eps}}\norm{v}+{\cal
O}(h^{(N+1)\delta-K-1})\left(\norm{u}+\norm{u_-}\right).
\endeq
On the other hand, applying $\chi_j$ and then $U_j$, $1\leq j\leq
L$, to the first equation of (\ref{3.14}), we obtain
\begin{eqnarray}
{\cases{(P^{(N)}_j-z)U_j\chi_j
u+\sum_{l=1}^{M_j}u_-(j)(l)e_{k(j,l)}=U_j\chi_j v+w_j, \cr
(U_j\chi_j u|e_{k(j,l)})=v_+(j)(l),\quad 1\leq l\leq M_j. \cr}}
\end{eqnarray}
Here the $L^2_{\theta}({\bf T}^2)$-norm of $w_j$ is ${\cal
O}(h^{(N+1)\delta})\left(\norm{u}+\norm{u_-}\right)$. It follows
that for $1\leq j\leq L$,
$$
\norm{\chi_j u}+\norm{u_-(j)}\leq \frac{{\cal O}(1)}{\eps
\widetilde{\eps}}\left(\norm{v}+\norm{v_+(j)}\right)+\frac{{\cal
O}(h^{(N+1)\delta})}{\eps
\widetilde{\eps}}\left(\norm{u}+\norm{u_-}\right).
$$
Here, as in (\ref{3.16}), we shall choose $N$ so large that
$\eps \widetilde{\eps}\geq \eps h \geq h^{K+1}\gg h^{(N+1)\delta}$, which together with (\ref{3.16}) implies the
injectivity, and hence the well-posedness of the Grushin problem
(\ref{3.14}). The bound (\ref{3.15}) follows.

The solution to (\ref{3.14}) is given by
$$
u=Ev+E_+v_+,\quad u_-=E_-v+E_{-+}v_+,
$$
and we recall (see~\cite{SjZw} and further references given there)
that the eigenvalues of $P_{\eps}$ in the rectangle (\ref{3.4})
are precisely the values $z$ for which the matrix $E_{-+}\in {\cal
L}(\comp^{M_1}\times \cdots \times \comp^{M_L},\comp^{M_1}\times
\cdots \times \comp^{M_L})$ is non-invertible. Arguing precisely
as in~\cite{MeSj} and~\cite{HiSj1}, we find that modulo an error
term ${\cal O}(h^{(N+1)\delta})$, $E_{-+}$ is a block-diagonal
matrix with the blocks $E_{-+}(z)(j)\in {\cal
L}(\comp^{M_j},\comp^{M_j})$, $1\leq j\leq L$, given by
$$
E_{-+}(z)(j)(m,n)=(z-z(j,k(j,m)))\delta_{mn},\,\, 1\leq m\leq n \leq M_j.
$$
\medskip
The discussion above is summarized in the following theorem.

\begin{theo}
Let $F$ stand for the mean value of $q$ along the Diophantine tori
$\Lambda_j$, $1\leq j\leq L$. When $\alpha_{1,j}$ and
$\alpha_{2,j}$ are the fundamental cycles  in $\Lambda_j$, $1\leq
j\leq L$, we write $S_j=(S_{1,j},S_{2,j})$ and
$k_j=(k(\alpha_{1,j}),k(\alpha_{2,j}))$ for the actions and the
Maslov indices of the cycles, respectively. Assume furthermore
that $\eps={\cal O}(h^{\delta})$, $\delta>0$ satisfies $\eps\geq
h^K$, for some $K$ fixed but arbitrarily large. Let
$\widetilde{\eps}>0$ be an additional small parameter such that
$\widetilde{\eps}\gg \eps$ and
$$
h^{1/3-\delta}<\widetilde{\eps}={\cal O}(h^{\delta}).
$$
We next make the global dynamical assumption {\rm (\ref{H0})}, and
assume that the differentials of the functions $p_{N,j}(\xi)$ and
$\Re \langle{q_j}\rangle(\xi)$, defined in Proposition {\rm 5.1},
are linearly independent, when $\xi=0$, $1\leq j\leq L$. In the
case when the $H_p$-flow is completely integrable, instead of {\rm
(\ref{H0})} we assume {\rm (\ref{1.31})}. Let $C>0$ be
sufficiently large. Then the eigenvalues of $P_{\eps}$ in the
rectangle \begeq \label{R} \abs{\Re
z}<\frac{\widetilde{\eps}}{C},\quad \abs{\Im z-\eps \Re
F_0}<\frac{\widetilde{\eps}\eps}{C}
\endeq
are given by
$$
P^{(N)}_j\left(h\left(k-\frac{k_j}{4}\right)-\frac{S_j}{2\pi},\eps;h\right)+{\cal
O}(h^{(N+1)\delta}),\quad k\in \z^2,\quad 1\leq j\leq L.
$$
Here $N\in \nat$ is such that $(N+1)\delta-K-1\gg 1$, and
$$
P_j^{(N)}(\xi,\eps;h)=p_{N,j}(\xi)+i\eps
\langle{q_j}\rangle(\xi)+{\cal O}(\eps^2)+{\cal O}(h).
$$
\end{theo}

Let us finally recall from (\ref{C}) the limiting operator
$P_j^{(\infty)}(hD_x,\eps;h)$, $1\leq j\leq L$, well-defined
modulo ${\cal O}((h,\eps,\xi)^{\infty})$ and such that for each
$N\in \nat$,
$$
P^{(\infty)}_j(\xi,\eps;h)=P^{(N)}_j(\xi,\eps;h)+{\cal
O}((h,\eps,\xi)^{N+1}).
$$
Then it follows from Theorem 5.2 that the eigenvalues of
$P_{\eps}$ in the domain (\ref{R}) are of the form
$$
P^{(\infty)}_j\left(h\left(k-\frac{k_j}{4}\right)-\frac{S_j}{2\pi},\eps;h\right)+{\cal
O}(h^{\infty}),\quad k\in \z^2,\quad 1\leq j\leq L.
$$
This completes the proof of Theorem 1.1.

\Remark. In the end of this section, we would like to mention that
the method of Grushin reduction, exploited in this section, has a
very long tradition and is closely related and essentially
equivalent to the so-called Feshbach projection
method---see~\cite{DJa} for a recent use of the latter in the
context of spectral analysis of the Pauli-Fierz Hamiltonian in
quantum field theory. We also refer to~\cite{SjZw} for a
systematic presentation of the Grushin method and of some of its
applications. In particular, in~\cite{SjZw}, it is explained how
the Feshbach method fits into the framework of Grushin problems.

\section{Asymptotic expansion of eigenvalues for large $\eps$}
\label{sectionGr}\setcounter{equation}{0}

In this section we let $\eps$ be sufficiently small but
independent of $h$ and the purpose is to compute all the
eigenvalues of the operator $P_{\eps}$ in a domain independent of
$h$, thereby proving Theorem 1.2. In doing so, we shall use the
results of section 4, and we shall also see that a suitable
quantum Birkhoff normal form construction will allow us to extend
the result to cover a certain $h$-dependent range of values of
$\eps$.

Let us recall from Proposition 3.1 that the leading symbol
(\ref{1.34}) of the operator (\ref{1.33}) can be reduced to the
normal form \begeq \label{gr1}
p(x,\xi,\eps)=p_{N-1}(\xi,\eps)+r_N(x,\xi,\eps),
\endeq
with $r_N(x,\xi,\eps)={\cal O}((\eps,\xi)^N)$, and
$p_{N-1}(\xi,\eps)=p(\xi)+i\eps q(\xi,\eps)$, with $p(\xi)$ real,
and $d_{\xi}p(0)$, $d_{\xi}\Re q(0,0)$ linearly independent. Then
from Proposition 4.1 we recall that the Hamilton-Jacobi equation
\begeq \label{gr2} p(x,\varphi'_x,\eps)-p_{N-1}(\xi,\eps)=0
\endeq
has the solution \begeq \label{gr3} \varphi(x,\xi,\eps)=x\cdot
\xi+\widetilde{\eps}\rho(x,\xi,\eps,\widetilde{\eps}),
\endeq
for $\xi={\cal O}(\widetilde{\eps})$, $\eps\leq \widetilde{\eps}$,
$\widetilde{\eps}^N/\eps^2\ll 1$. Here $\norm{\rho'_x}_{H^s}={\cal
O}(\widetilde{\eps}^{N-1}/\eps)$, $s>1$. The construction in
section 4 clearly works the same way for $x$ complex with
$\abs{\Im x}<1/{\cal O}(1)$ and for $\xi$ complex with $\xi={\cal
O}(\widetilde{\eps})$, if we work in $H^s$-spaces on each torus $\Im
x={\rm Const}$. If we normalize the choice of the
grad-periodic function $\rho$ by putting $\rho(0,\xi)=0$, then
$\rho(x,\xi)$ becomes holomorphic in $x$ and $\xi$, for $\abs{\Im
x}<1/{\cal O}(1)$, $\xi={\cal O}(\widetilde{\eps})$. Write \begeq
\label{gr4} \varphi(x,\xi,\eps)=x\cdot \xi+\psi(x,\xi,\eps),
\endeq
so that \begeq \label{gr5} \abs{\psi(x,\xi,\eps)}={\cal
O}\left(\frac{\widetilde{\eps}^N}{\eps}\right),
\endeq
for $\abs{\Im x}\leq 1/{\cal O}(1)$, $\xi={\cal
O}(\widetilde{\eps})$. From the Cauchy inequalities we get \begeq
\label{gr6}
\partial_x^{\alpha}\partial_{\xi}^{\beta} \psi={\cal
O}\left(\frac{\widetilde{\eps}^{N-\abs{\beta}}}{\eps}\right)
\endeq
in the same region.

Consider the actions (independent of $x$):
\begin{eqnarray}
\label{gr7} \eta_j & = & \frac{1}{2\pi}(\varphi(x+2\pi
e_j,\xi,\eps)-\varphi(x,\xi,\eps)) \\ \nonumber & = & \xi_j +
\frac{1}{2\pi}(\psi(x+2\pi
e_j,\xi,\eps)-\psi(x,\xi,\eps))=\xi_j+{\cal
O}\left(\frac{\widetilde{\eps}^N}{\eps}\right),
\end{eqnarray}
where $e_1=(1,0)$, $e_2=(0,1)$. By the implicit function theorem and
the Cauchy inequalities, we see that the map $\xi\mapsto \eta(\xi)$,
defined for $\xi={\cal O}(\widetilde{\eps})$, is bijective with an
inverse $\xi(\eta)$, and
\begeq
\label{gr8}
\eta(\xi)=\xi+{\cal
O}\left(\frac{\widetilde{\eps}^N}{\eps}\right),\quad
\xi(\eta)=\eta+{\cal O}\left(\frac{\widetilde{\eps}^N}{\eps}\right).
\endeq
Again derivatives of these maps can be estimated by means of the
Cauchy inequalities. Now write $\varphi(x,\eta,\eps)$ instead of
$\varphi(x,\xi,\eps)$. We have \begeq \label{gr9}
\varphi(x,\eta,\eps)=x\cdot \eta+{\cal
O}\left(\frac{\widetilde{\eps}^N}{\eps}\right),\quad \abs{\Im
x}<1/{\cal O}(1),\,\,\abs{\Re x}={\cal
O}(1),\,\,\,\abs{\eta}<{\cal O}(\widetilde{\eps}),
\endeq
and $\varphi(x,\eta,\eps)-x\cdot \eta$ is single-valued.

Consider the canonical transformation \begeq \label{gr10}
\kappa=\kappa_{\eps}:
(\partial_{\eta}\varphi(x,\eta,\eps),\eta)\mapsto (x,\partial_x
\varphi(x,\eta,\eps)).
\endeq
Differentiating the identity (\ref{gr7}), we get
$$
\partial_{\eta_k} \varphi(2\pi
e_j+x,\eta,\eps)=\partial_{\eta_k}\varphi(x,\eta,\eps)+2\pi
\delta_{j,k},
$$
so it is clear that $\kappa_{\eps}$ is a map from a complex \neigh{} of the
form $\abs{\Im y}<1/{\cal O}(1)$, $\eta={\cal O}(\widetilde{\eps})$ of
the zero section in $T^*{\bf T}^2$ onto another \neigh{} which
contains a \neigh{} of the same form.

By construction, $p\circ\kappa_{\eps}(y,\eta)=p_{N-1}(\xi,\eps)$, and
since $\xi=\eta+{\cal O}(\widetilde{\eps}^N/\eps)$, we get, expanding
the right hand side as a function of $(\eps,\eta)$,
\begeq
\label{gr11}
p\circ \kappa_{\eps}(y,\eta)=p_{N-1}(\eta,\eps)+{\cal
O}\left(\frac{\widetilde{\eps}^N}{\eps}\right)=p(\eta)+i\eps
q(\eta,\eps)+{\cal O}\left(\frac{\widetilde{\eps}^N}{\eps}\right),
\endeq
and this is a function of $(\eta,\eps)$ only.

We next want to implement $\kappa$ by a Fourier integral operator.
For that, it will be convenient to work in a fixed \neigh{} of the
zero section. Let $P=P(x,hD_x,\eps;h)$ be an
$h$-pseudodifferential operator with the leading symbol
$p(x,\xi,\eps)$ as in (\ref{gr1}) and assume, as we may, for
simplicity, that $p(0)=0$, $q(0,\eps)=0$. Put
$\widetilde{h}=h/\widetilde{\eps}$, so that \begeq \label{gr12}
\frac{1}{\widetilde{\eps}}P(x,hD_x,\eps;h)=\frac{1}{\widetilde{\eps}}P(x,\widetilde{\eps}\widetilde{h}D_x,\eps;h).
\endeq
As an $\widetilde{h}$-pseudodifferential operator,
$\widetilde{\eps}^{-1}P$ has a well-defined symbol in a fixed
\neigh{} of $\xi=0$, and the leading symbol will be \begeq
\label{gr13}
\frac{1}{\widetilde{\eps}}p(x,\widetilde{\eps}\xi,\eps)=\frac{1}{\widetilde{\eps}}p(\widetilde{\eps}\xi)+i\eps
\frac{1}{\widetilde{\eps}}q(\widetilde{\eps}\xi,\eps)+{\cal
O}(\widetilde{\eps}^{N-1}).
\endeq
Notice that $\widetilde{\eps}^{-1}p(\widetilde{\eps}\xi)$,
$\widetilde{\eps}^{-1} q(\widetilde{\eps}\xi,\eps)$ are uniformly
bounded in a fixed complex domain when
$\widetilde{\eps}\rightarrow 0$, and that \begeq \label{gr14}
d_{\xi} \frac{1}{\widetilde{\eps}}p(\widetilde{\eps}\xi),\,\,
d_{\xi}\frac{1}{\widetilde{\eps}}q(\widetilde{\eps}\xi,0)\quad
\wrtext{are linearly independent for}\,\, \xi=0,
\endeq
uniformly when $\widetilde{\eps}\rightarrow 0$.

Let
\begeq
\label{gr15}
U u(x)=h^{-2}\int\!\!\! \int e^{\frac{i}{h}(\varphi(x,\eta,\eps)-y\cdot
\eta)}a(x,\eta;h)u(y)\,dy\,d\eta,
\endeq
be an elliptic $h$-Fourier integral operator, associated to
$\kappa$. The change of variables,
$\widetilde{h}=h/\widetilde{\eps}$,
$\eta=\widetilde{\eps}\widetilde{\eta}$, gives \begeq \label{gr16}
U u(x)=\widetilde{h}^{-2} \int\!\!\!\int
e^{\frac{i}{\widetilde{h}}
(\frac{1}{\widetilde{\eps}}\varphi(x,\widetilde{\eps}\widetilde{\eta},\eps)-y\cdot
\widetilde{\eta})}a(x,\widetilde{\eps}\widetilde{\eta};h)u(y)\,dy\,d\widetilde{\eta}.
\endeq
From (\ref{gr9}) we get
\begeq
\label{gr17}
\frac{1}{\widetilde{\eps}}\varphi(x,\widetilde{\eps}\widetilde{\eta},\eps)=x\cdot
\widetilde{\eta}+{\cal
O}\left(\frac{\widetilde{\eps}^{N-1}}{\eps}\right),\,\,\abs {\Im
x}<\frac{1}{{\cal O}(1)},\quad \abs{\widetilde{\eta}}={\cal O}(1).
\endeq
The phase is therefore uniformly non-degenerate and the
corresponding canonical transformation is \begeq \label{gr18}
\widetilde{\kappa}:
(\varphi'_{\eta}(x,\widetilde{\eps}\widetilde{\eta},\eps),\widetilde{\eta})\rightarrow
(x,\frac{1}{\widetilde{\eps}}
\varphi'_x(x,\widetilde{\eps}\widetilde{\eta},\eps)).
\endeq

The conjugated operator $\widetilde{P}=U^{-1}
\frac{1}{\widetilde{\eps}}P U$ is a uniformly well-behaved
$\widetilde{h}$-pseu\-do\-dif\-fe\-ren\-tial operator with leading symbol
obtained from (\ref{gr11}) by division by $\widetilde{\eps}$ and
substitution $\eta=\widetilde{\eps}\widetilde{\eta}$:
\begeq
\label{gr19}
\frac{1}{\widetilde{\eps}}p(\widetilde{\eps}\widetilde{\eta})+i\eps
\frac{1}{\widetilde{\eps}}q(\widetilde{\eps}\widetilde{\eta},\eps)+{\cal
O}\left(\frac{\widetilde{\eps}^{N-1}}{\eps}\right)=:\widetilde{p}(\widetilde{\eta},\widetilde{\eps},\eps),
\endeq
and this is independent of $y$.

We now want to further simplify $\widetilde{P}$ by conjugation by an
elliptic $\widetilde{h}$-pseudo\-diffe\-ren\-tial operator, $e^A$, where $A$
is an $\widetilde{h}$-pseudodifferential operator of order $0$. The
construction will be uniform in $\widetilde{\eps}$ with a power
degeneration in $\eps$, that we shall control. First we recall that
$e^A \widetilde{P}e^{-A}=e^{{\rm ad}_A}\widetilde{P}=\sum
\frac{1}{k!}{\rm ad}_A^k \widetilde{P}$. Let the full symbol of $A$ be
$\sum_{k=0}^{\infty} \widetilde{h}^k a_k$. Then
\begin{eqnarray*}
e^A \widetilde{P}e^{-A}& = &
\sum_{k=0}^{\infty}\sum_{l=0}^{\infty}\sum_{j_1=0}^{\infty}\ldots
\sum_{j_k=0}^{\infty}\frac{1}{k!}
\widetilde{h}^{j_1+\ldots\,+j_k+l+k}\left(\frac{1}{\widetilde{h}}{\rm
ad}_{a_{j_1}}\right)\ldots \left(\frac{1}{\widetilde{h}}{\rm
ad}_{a_{j_k}}\right)(\widetilde{p}_l) \\
& = & \sum_{n=0}^{\infty} \widetilde{h}^n s_n,
\end{eqnarray*}
with $s_0=\widetilde{p}_0$,
$s_1=\frac{1}{i}H_{a_0}\widetilde{p}_0+\widetilde{p}_1=
iH_{\widetilde{p}_0}a_0+\widetilde{p}_1,\ldots\,$,$s_{n+1}=iH_{\widetilde{p}_0}a_n+\widetilde{s}_{n+1}$,
where $\widetilde{s}_{n+1}$ only depends on $a_0,\ldots,\,a_{n-1}$ and
is a sum of coefficients for $\widetilde{h}^{n+1}$ for terms
$$
\frac{1}{k!}\widetilde{h}^{j_1+\ldots\,+j_k+l+k}\left(\frac{1}{\widetilde{h}}{\rm
ad}_{a_{j_1}}\right)\ldots \left(\frac{1}{\widetilde{h}}{\rm
ad}_{a_{j_k}}\right)(\widetilde{p}_l),
$$
with $j_1+\ldots\,+j_k+l+k\leq n+1$, $j_1,\ldots\, j_k<n$. By
solving $\eps$-degenerated Cauchy-Riemann equations, we see that
the $a_j$ can be successively chosen so that $s_j$ are independent
of $y$.

Assume by induction that $\nabla a_j={\cal O}(\eps^{-1-2j})$, for
$j<n-1$ (in a complex domain, so that we have the same estimates
for the derivatives of $\nabla a_j$). Then the general term in
$\widetilde{s}_{n+1}$ is \begeq \label{gr19.1} {\cal O}(1)
\eps^{-1-2j_1}\ldots \eps^{-1-2j_k}={\cal
O}(1)\left(\frac{1}{\eps}\right)^{2(j_1+\ldots\,+j_k)+k}.
\endeq
Here
$$
2(j_1+\ldots\,+j_k)+k=2(j_1+\ldots\,+j_k+k)-k\leq
2(n+1-l)-k=2n+2-2l-k,
$$
and so the quantity (\ref{gr19.1}) is ${\cal O}(1)\eps^{-2n}$
except possibly when $2l+k<2$, i.e. when $k=l=0$ or when $k=1$,
$l=0$. In the first case we get the coefficient for
$\widetilde{h}^{n+1}$ in $\widetilde{p}_0$ which is $0$. In the
second case, we get the coefficient for $\widetilde{h}^{n+1}$ in
$\widetilde{h}^{j_1+1}(\frac{1}{\widetilde{h}}{\rm
ad}_{a_{j_1}})(\widetilde{p}_0)$ with $j_1<n$, which is ${\cal
O}(1)\eps^{-1-2j_1}$. Here $1+2j_1\leq 2n$. Thus
$\widetilde{s}_{n+1}={\cal O}(\eps^{-2n})$, in a complex domain.
We can choose $a_n$ grad-periodic, with $iH_{\widetilde{p}_0}a_n=
-\widetilde{s}_{n+1}+\langle{s_{n+1}(\cdot,\widetilde{\eta})}\rangle$
and with $\nabla a_n={\cal O}(\eps^{-1-2n})$. This completes the
induction step and we conclude that we can find $a_k$ with $\nabla
a_k={\cal O}(\eps^{-1-2k})$ in a fixed complex \neigh{} of
$\widetilde{\eta}=0$ such that if
$$
A^{(N_1)}=\sum_{k=0}^{N_1-1} \widetilde{h}^k a_k,
$$
then \begeq \label{gr19.5}
\widetilde{P}^{(N_1)}:=e^{A^{(N_1)}}U^{-1}\frac{1}{\widetilde{\eps}}P
U e^{-A^{(N_1)}}=\sum_{n=0}^{\infty} \widetilde{h}^n
\widetilde{p}_n^{(N_1)},
\endeq
where \begeq \label{gr19.6} \widetilde{p}_n^{(N_1)}={\cal
O}(\eps^{-2(n-1)_+}),
\endeq
$\widetilde{p}_0^{(N_1)}=
\widetilde{p}(\widetilde{\eta},\widetilde{\eps},\eps)$ in
(\ref{gr19}) and
$\widetilde{p}_n^{(N_1)}=\widetilde{p}_n^{(\infty)}$ is
independent of $y$ and $N_1$, for $n\leq N_1$. Here $N_1$ can be
taken arbitrarily large.

We shall next look at the quasi-modes. In view of the estimates
(\ref{gr19.6}) on $\widetilde{p}_n^{(N_1)}$, we first assume that
\begeq \label{gr20} \frac{\widetilde{h}}{\eps^2}\ll 1.
\endeq
This gives us the quasi-eigenvalues of $\widetilde{\eps}^{-1} P$,
\begeq \label{gr21} \sum_{n=0}^{N_1} \widetilde{h}^n
\widetilde{p}_n^{(\infty)}\left(\widetilde{h}\left(k-\frac{k_0}{4}\right)-\frac{S}{2\pi\widetilde{\eps}},
\widetilde{\eps},\eps\right) +{\cal
O}\left(\frac{\widetilde{h}^{N_1+1}}{\eps^{2(N_1+1)}}\right),
\endeq
for $k\in \z^2$ with
$\widetilde{h}\left(k-\frac{k_0}{4}\right)-\frac{S}{2\pi\widetilde{\eps}}={\cal
O}(1)$. From (\ref{gr21}) we get the ``leading" values
$$
\widetilde{p}_0^{(\infty)}\left(\widetilde{h}\left(k-\frac{k_0}{4}\right)-\frac{S}{2\pi\widetilde{\eps}},
\widetilde{\eps},\eps\right)
=\widetilde{p}\left(\widetilde{h}\left(k-\frac{k_0}{4}\right)-\frac{S}{2\pi\widetilde{\eps}},\widetilde{\eps},\eps\right),
$$
and from (\ref{gr19}) we infer that the distance between the
neighboring ``leading" values is $\geq \eps \widetilde{h}/{\cal
O}(1)$. The $\widetilde{\eta}$-gradient of $\widetilde{h}^n
\widetilde{p}_n^{(\infty)}(\widetilde{\eta},\widetilde{\eps},\eps)$
is ${\cal O}(1)\widetilde{h}^n\eps^{-2(n-1)}$ for $1\leq n\leq
N_1$, and this is $\ll \eps$, since by (\ref{gr20}),
$\widetilde{h}^n
\eps^{-2(n-1)}=(\frac{\widetilde{h}}{\eps^2})^n\eps^2\ll \eps^2$.
It is therefore clear that the distance between neighboring values
of
$$
\sum_{n=0}^{N_1} \widetilde{h}^n
\widetilde{p}_n^{(\infty)}\left(\widetilde{h}\left(k-\frac{k_0}{4}\right)-\frac{S}{2\pi\widetilde{\eps}},
\widetilde{\eps},\eps\right)
$$
is $\geq \eps\widetilde{h}/{\cal O}(1)$, and (\ref{gr21}) gives
``distinct'' quasi-eigenvalues if this minimal separation is $\gg
\widetilde{h}^{N_1+1}/\eps^{2(N_1+1)}$. We therefore must
strengthen (\ref{gr20}) to
$\widetilde{h}^{N_1+1}/\eps^{2(N_1+1)}\ll \eps\widetilde{h}$,
i.e.: \begeq \label{gr22}
\left(\frac{\widetilde{h}}{\eps^2}\right)^{N_1+1}\ll \eps
\widetilde{h},
\endeq
or \begeq \label{gr23}
\left(\frac{h}{\widetilde{\eps}\eps^2}\right)^{N_1+1}\ll
\frac{\eps h}{\widetilde{\eps}}.
\endeq
This will be satisfied if we choose $\eps\geq h^{\delta}$ for some
fixed $\delta\in (0,1/3)$, $\widetilde{\eps}\geq \eps$, and finally $N_1$ large enough.

Now in (\ref{gr12}) we make the substitution
$\widetilde{\eps}\rightarrow \mu \widetilde{\eps}$,
$\widetilde{h}\rightarrow \widetilde{h}/\mu$, $\mu\sim 1$, and
obtain the isospectral operator
$$
\frac{1}{\widetilde{\eps}}P(x,hD_x)=\frac{1}{\widetilde{\eps}}P(x,\widetilde{\eps}\widetilde{h}D_x)=
\mu \frac{1}{\mu \widetilde{\eps}}P\left(x,\mu
\widetilde{\eps}\frac{\widetilde{h}}{\mu}D_x\right),
$$
with the quasi-eigenvalues
$$
\mu \sum_{n=0}^{N_1} \left(\frac{\widetilde{h}}{\mu}\right)^n
\widetilde{p}_n^{(\infty)}\left(\frac{\widetilde{h}}{\mu}\left(k-\frac{k_0}{4}\right)-\frac{S}{2\pi
\widetilde{\eps}\mu}, \mu\widetilde{\eps},\eps\right)+{\cal
O}\left(\frac{\widetilde{h}^{N_1+1}}{\mu^{N_1}
\eps^{2(N_1+1)}}\right).
$$
We deduce as in~\cite{MeSj} that for each $n=1,\ldots\, ,N_1$,
\begeq \label{gr24} \mu^{1-n}
\widetilde{p}_n^{(\infty)}\left(\frac{\eta}{\mu},\mu\widetilde{\eps},\eps\right)=
\widetilde{p}_n^{(\infty)}(\eta,\widetilde{\eps},\eps).
\endeq
We use this to define $\widetilde{p}_n^{(\infty)}(\eta,1,\eps)$ by
\begeq \label{gr25} \widetilde{p}_n^{(\infty)}(\eta,1,\eps)=
\widetilde{\eps}^{1-n}\widetilde{p}_n^{\infty}\left(\frac{\eta}{\widetilde{\eps}},\widetilde{\eps},\eps\right).
\endeq

Now recall from Proposition 4.1 that we have the conditions
$\widetilde{\eps}\geq \eps$, $\widetilde{\eps}^{N}/\eps^2\ll 1$:
\begeq \label{gr26} \eps\leq \widetilde{\eps}\ll \eps^{2/N}.
\endeq
Then $\widetilde{p}_n^{(\infty)}(\eta,\widetilde{\eps},\eps)$ is a
well-defined analytic function for $\abs{\eta}<{\cal O}(1)$ and by
choosing $\mu$ to be of the order of magnitude $\eps^{2/N}$, we
see that $\widetilde{p}_n^{(\infty)}(\eta,1,\eps)$ is well-defined
and holomorphic for $\abs{\eta}\ll \eps^{2/N}$.

Now restrict $\eps$ by imposing \begeq \label{gr27} \eps\geq
h^{1/3-\delta},\quad \wrtext{for some}\,\,\delta>0.
\endeq
Then
$$
\frac{h}{\eps^2 \widetilde{\eps}}\leq \frac{h}{\eps^3}\leq
h^{3\delta},
$$
and the remainder in (\ref{gr21}) can be bounded by
$$
{\cal
O}(1)\left(\frac{h}{\widetilde{\eps}\eps^2}\right)^{N_1+1}\leq
{\cal O}(1) h^{3\delta(N_1+1)}.
$$
Combining (\ref{gr21}) with (\ref{gr25}) we get that the
quasi-eigenvalues of $P_{\eps}$ take the form \begeq \label{gr28}
\sum_{n=0}^{N_1} h^n
\widetilde{p}_n^{(\infty)}\left(h\left(k-\frac{k_0}{4}\right)-\frac{S}{2\pi},1,\eps\right)+{\cal
O}(h^{3\delta(N_1+1)}),
\endeq
for $k\in \z^2$ with $\abs{h(k-\frac{k_0}{4})-\frac{S}{2\pi}}\ll
\eps^{2/N}$.

\begin{theo}
Assume that $h^{1/3-\delta}<\eps \leq \eps_0\ll 1$, for some
$\delta>0$, and let $\widetilde{\eps}$ be an additional small
parameter such that $\widetilde{\eps}\geq \eps$ and \begeq
\label{gr29} \frac{\widetilde{\eps}^{N-3}}{\eps^2}\ll 1.
\endeq
Recall that $F$ stands for the mean value of $q$ over the
invariant tori $\Lambda_j$, $1\leq j\leq L$, and that $S_j\in
\real^2$ and $k_j\in \z^2$ are the actions and Maslov indices of
the fundamental cycles in $\Lambda_j$, $1\leq j\leq L$,
respectively. In the general case we make the assumption {\rm
(\ref{H0})} and when the $H_p$-flow is completely integrable, we
assume {\rm (\ref{1.31})}. Assume finally that the differentials
of the functions $p_{N,j}(\xi)$ and $\Re
\langle{q_j}\rangle(\xi)$, defined in Proposition {\rm 5.1}, are
linearly independent when $\xi=0$. Let $C>0$ be large enough. Then
the eigenvalues of $P_{\eps}$ in
$$
\abs{\Re z}\leq \frac{\widetilde{\eps}}{C},\quad \abs{\Im z-\eps
F}\leq \frac{\eps \widetilde{\eps}}{C}
$$
are given by \begeq \label{gr30} z(j,k)=\sum_{n=0}^{N_1} h^n
\widetilde{p}_{j,n}^{(\infty)}\left(h\left(k-\frac{k_j}{4}\right)-\frac{S_j}{2\pi},1,\eps\right)+{\cal
O}(h^{3\delta(N_1+1)}),\quad k\in \z^2,\,\,\ 1\leq j\leq L,
\endeq
with $\widetilde{p}_{j,n}^{(\infty)}(\xi,1,\eps)={\cal
O}(\eps^{-2(n-1)_+})\widetilde{\eps}^{-n}$, $1\leq j\leq L$,
holomorphic for $\xi={\cal O}(\widetilde{\eps})$, and
$$
\widetilde{p}_{j,0}^{(\infty)}(\xi,\eps)=p_j(\xi)+i\eps
q_j(\xi,\eps)+{\cal
O}\left(\frac{\widetilde{\eps}^{N}}{\eps}\right).
$$
Here $p_j$ is real on the real domain, and the differentials of
$p_j(\xi)$ and $\Re q_j(\xi,\eps)$ are linearly independent when
$\xi=\eps=0$. The integers $N$ and $N_1$ in {\rm (\ref{gr29})} and
{\rm (\ref{gr30})} can be taken arbitrarily large.
\end{theo}

When proving Theorem 6.1, we assume for simplicity that $F=0$ and
recall the leading symbol of $P_{\eps}$ from (\ref{1.34}), \begeq
\label{gr37} p_{N-1}(\xi)+i\eps \langle{q}\rangle(\xi)+{\cal
O}(\eps^2)+\eps{\cal O}(\xi^{N-1})+{\cal O}(\xi^{N}).
\endeq
(Here we have replaced $N$ by $N-1$.) Recall also the canonical
transformation $\widehat{\kappa}$ defined in (\ref{2.2.5}) and
(\ref{2.2.6}), and that we have a second canonical transformation
given by (\ref{gr10}),
$$
\kappa_{\eps}:
(\partial_{\eta}\varphi(x,\eta,\eps),\eta)\rightarrow (x,
\partial_x \varphi(x,\eta,\eps))
$$
defined in a complex domain, with the restriction that
$\abs{\eta}={\cal O}(\widetilde{\eps})$. Moreover,
$\varphi(x,\eta,\eps)=x\cdot \eta+{\cal
O}(\widetilde{\eps}^{N}/\eps)$, and we conclude that
$\kappa_{\eps}=1+{\cal O}(\widetilde{\eps}^{N-1}/\eps)$, so that
the differential of $\kappa_{\eps}$ is equal to ${\cal
O}(\widetilde{\eps}^{N-2}/\eps)$, while the higher order
differentials satisfy $\nabla^k \kappa_{\eps}={\cal
O}(\widetilde{\eps}^{N-1-k}/\eps)$, $k\geq 2$.

As in section 5, we can represent the real phase space $L_0=T^*{\bf
T}^2$ by the IR-manifold $\xi=\frac{2}{i}\frac{\partial
\Phi_0}{\partial x}(x)$, $\Phi_0(x)=\frac{1}{2}(\Im x)^2$, with
$\xi=0$ corresponding to $\Im x=0$, and $x$ varying in a \neigh{}
of $\{x; \Im x=0\}$ in ${\bf T}^2+i\real^2$. Then
$L_{\eps}=\kappa_{\eps}(L_0)$ is represented by
$$
\xi=\frac{2}{i}\frac{\partial \Phi_{\eps}}{\partial x},\quad
\abs{\Im x}<{\cal O}(\widetilde{\eps}),
$$
with $\Phi_{\eps}-\Phi_0$, $\nabla \Phi_{\eps}-\nabla \Phi_0={\cal
O}(\widetilde{\eps}^{N-1}/\eps)$,
$\nabla^k(\Phi_{\eps}-\Phi_0)={\cal
O}(\widetilde{\eps}^{N-k}/\eps)$, $k\geq 1$.

Consider the IR-manifold $\widetilde{L}_{\eps}$ represented by
$$
\xi=\frac{2}{i}\frac{\partial  \widetilde{\Phi}_{\eps}}{\partial
x}(x),\quad \widetilde{\Phi}_{\eps}(x)=\chi\left(\frac{\Im
x}{\widetilde{\eps}}\right)\Phi_{\eps}(x)
+\left(1-\chi\left(\frac{\Im
x}{\widetilde{\eps}}\right)\right)\Phi_0(x),
$$
where $\chi\in C^{\infty}_0(\real^2)$ is a standard cut-off to a
\neigh{} of $0$, so that $\widetilde{L}_{\eps}$ coincides with
$L_{\eps}$ near $\xi=0$ and with $L_0$ further away from this set.
Moreover,
\begeq \label{gr39}
\nabla(\widetilde{\Phi}_{\eps}-\Phi_0),\quad
\nabla(\widetilde{\Phi}_{\eps}-\Phi_{\eps})={\cal
O}\left(\frac{\widetilde{\eps}^{N-2}}{\eps}\right).
\endeq
On the other hand, from section 5 we know that $\abs{\Im
P_{\eps}\circ \widehat{\kappa}}$ along $L_0$ is of the order of
magnitude $\eps \widetilde{\eps}$ in the region where
$\abs{\xi}\sim \widetilde{\eps}$ and $\abs{\Re
P_{\eps}}/{\widetilde{\eps}}$ is small. In view of (\ref{gr39}) we
then have the same conclusion for $\Im P_{\eps}\circ
\widehat{\kappa}$ along ${\widetilde{L}_{\eps}}$, provided that
$\widetilde{\eps}^{N-2}/\eps\ll \eps \widetilde{\eps}$, i.e. for
$$
\frac{\widetilde{\eps}^{N-3}}{\eps^2}\ll 1,
$$
which is precisely (\ref{gr29}).

The final IR-manifold $\widehat{\Lambda}_{\eps}$ that we choose is
then given by
$$
\exp(i\eps H_G)\circ (\kappa_1^{-1})\circ (\kappa^{(N)})^{-1}
\circ \widehat{\kappa}(\widetilde{L}_{\eps})
$$
in a complex \neigh{} of $\Lambda_1$, and we do a similar
modification near the other tori $\Lambda_j$, $2\leq j\leq
L$---see also (\ref{3.2.1}). Here we are also using that
$\widehat{\kappa}(\widetilde{L}_{\eps})$ coincides with
$\widehat{\kappa}(L_0)$ in the region corresponding to $\xi\notin
\supp \chi\left(\frac{\cdot}{\widetilde{\eps}}\right)$, $\xi={\cal
O}(\widetilde{\eps})$. Then, along $\widehat{\Lambda}_{\eps}$,
$\abs{\Im P_{\eps}} \geq \eps \widetilde{\eps}/{\cal O}(1)$
outside any $\widetilde{\eps}$-\neigh{} of the union of the
$\Lambda_j$'s, $1\leq j\leq L$, in the energy slice $\abs{\Re
P_{\eps}}\leq \widetilde{\eps}/{\cal O}(1)$, and moreover the
holomorphic canonical transformation $\widehat{\kappa}\circ
\kappa_{\eps}$ maps a region $\abs{\xi}<{\cal
O}(\widetilde{\eps})$ in $T^*{\bf T}^2$ onto the intersection of
$\widehat{\Lambda}_{\eps}$ with a similar \neigh{} of $\Lambda_1$.

Using these facts and the normal form (\ref{gr28}), we get Theorem
6.1 by solving a globally well-posed Grushin problem as in section
5. We then take $N_1\rightarrow \infty$, and put
$\widetilde{\eps}=\eps^{1/\widetilde{N}}$, where $\widetilde{N}$
can be arbitrarily large. This completes the proof of Theorem 1.2.

\section{The dynamical condition and the isoenergetic KAM theorem}
\label{section5} \setcounter{equation}{0}

\subsection{The completely integrable case.} In the beginning of
this section, we shall discuss the case when the $H_p$-flow is
completely integrable. To be precise, as in section 2, we shall
assume that there exists an analytic real valued function $f$ on
$T^*M$ such that $H_p f=0$. The sets $\Lambda_a=p^{-1}(0)\cap
f^{-1}(a)$, $a\in J$, then form a singular foliation of
$p^{-1}(0)$, and we get a corresponding decomposition of the
energy surface, \begeq \label{7.0} p^{-1}(0)=\bigcup_{a\in J}
\Lambda_{a}.
\endeq
Here $J\subset \real$ is a compact interval and we know that for
each $a\in J$, $\Lambda_a$ is a compact $H_p$-invariant set, such
that we have the continuity property (\ref{1.22}). We now
introduce a global assumption that for each $a\in J$, except for a
finite set $S\subset J$, the set $\Lambda_a$ is a finite disjoint
union of analytic flow-invariant Lagrangian tori, depending
analytically on $a$, $\Lambda_a=\cup_{j=1}^ {N(a)} \Lambda_{a,j}$.
Here $N(a)$ is a locally constant bounded function on $J\backslash
S$. As we shall see, this assumption is satisfied in the case when
$M$ is a convex analytic surface of revolution, with the
exceptional set $S$ in that case corresponding to the endpoints of
$J$, and with $N(a)=1$. In order to simplify the discussion
notationally, in what follows we shall assume that each
$\Lambda_a$, $a\in J\backslash S$, is connected.

\Remark. In some situations, it turns out to be more appropriate
to assume that $J$ in (\ref{7.0}) is a circle, or, more generally,
a connected graph with finitely many vertices and edges, with $S$
then being the set of vertices---see
also~\cite{CdVSan},\cite{San}. Assuming that $J$ is a graph, the
following discussion goes through with minor modifications---see
also the end of this section.

\vskip 2mm
It is now classical (see~\cite{Ar}) that locally, near any regular
torus $\Lambda_a$ in $T^*M$, we may introduce analytic
action-angle coordinates $(x_1,x_2;\xi_1,\xi_2)$ with $x_j\in S^1$
and $\xi_j\in {\rm neigh}(0,\real)$, $j=1,2$, so that when
expressed in terms of these coordinates, $p$ becomes a function of
$\xi_1$ and $\xi_2$, $p=p(\xi_1,\xi_2)$. The rotation number of
$\Lambda_a$, $\omega(a)$, is then defined as the ratio of the
frequencies \begeq \label{7.0.1} [\partial_{\xi_2}p :
\partial_{\xi_1}p],
\endeq
viewed as an element of the real projective line. We shall assume
that $\omega(a)$ depends analytically on $a\in J\backslash S$ and
that $a\mapsto \omega(a)$ is not identically constant on any
connected component of $J\backslash S$. Following~\cite{Ar}, we
recall then that the torus $\Lambda_a\subset p^{-1}(0)$ satisfies
the isoenergetic condition when the map $a\mapsto \omega(a)$ is a
local diffeomorphism.

When $q$ is a bounded holomorphic function in a tubular \neigh{}
of $T^*M$, that is for simplicity assumed to be real on the real
domain, we shall consider the
limit of \begeq \label{7.1}
\langle{q}\rangle_T=\frac{1}{T}\int_0^T q\circ
\exp(tH_p)\,dt,\quad T\rightarrow \infty,
\endeq
taken along an invariant torus $\Lambda_a$, $a\in J\backslash S$. All the results below will
be valid if $\langle{q}\rangle_T$ is replaced by
$\langle{q}\rangle_{T,K}$, with $K\in C^{\infty}_0$.

When analyzing (\ref{7.1}), we switch to the action-angle variables, so that
$\Lambda_a$ becomes the standard torus ${\bf T}^2$, $H_p$ along
$\Lambda_a$ becomes $b_1\partial_{x_1}+b_2\partial_{x_2}$, $b_1$, $b_2\in
\real$, and
$$
q=q(x_1,x_2).
$$
Then as in (\ref{1.28}), we get \begeq \label{7.1.1}
\langle{q}\rangle_T(x)=\frac{1}{T}\int_0^T q(x+tb)\,dt=\sum_{k\in
{\rm \bf z}^2} \widehat{q}(k)e^{ikx}
\widehat{1_{[0,1]}}(-Tb\cdot k).
\endeq
It follows as in (\ref{1.29}) that the limit of
$\langle{q}\rangle_T$ along $\Lambda_a$ is equal to \begeq
\label{7.1.2} \widehat{q}(0)+\sum_{0\neq k,\, b\cdot k=0}
\widehat{q}(k) e^{ikx}.
\endeq
We conclude that when $\Lambda_a$ is irrational (in the sense that
its rotation number $\omega(a)=b_2/b_1$ is irrational), the limit
is equal to the average of $q$  over the torus,
$\langle{q}\rangle_{\Lambda_a}$, computed with respect to the
natural absolutely continuous measure on $\Lambda_a$, with respect
to which the $H_p$-flow is ergodic. In the rational case, we write
$\omega(a)=m/n$, where $m\in \z$, $n\in \nat$, are relatively
prime, and introduce
$$
k(\omega(a)):=\abs{m}+\abs{n}.
$$
Using the smoothness of $q$, it follows then from (\ref{7.1.2})
that the limit of (\ref{7.1}) along a rational torus is equal to
$$
\langle{q}\rangle_{\Lambda_a}+{\cal
O}\left(\frac{1}{k(\omega(a))^{\infty}}\right).
$$
Here the integration measure used in the definition of the torus
average is absolutely continuous and comes from the diffeomorphism
between $\Lambda_a$ and ${\bf T}^2$, given by the action-angle
coordinates.

Let us summarize the discussion above in the following
proposition.

\begin{prop}
Assume that the $H_p$-flow on $T^*M$ is completely integrable and,
more precisely, that {\rm (\ref{7.0})} holds true, with
$\Lambda_a$ being an analytic invariant Lagrangian torus when
$a\in J\backslash S$, for a finite set $S$. When $q$ is a real-valued
bounded analytic function on $T^*M$ and $a\in J\backslash S$, we define
the closed interval $Q_{\infty}(a)\subset \real$ to be the range
of the limit of $\langle{q}\rangle_T$, as $T\rightarrow \infty$,
restricted to the invariant torus $\Lambda_a\subset p^{-1}(0)$. We
then have
\begin{itemize}
\item When the rotation number $\omega(a)\notin {\bf Q}$, then
$Q_{\infty}(a)=\{ \langle{q}\rangle_{\Lambda_a} \}$. \item When
$\omega(a)=m/n$, with $m\in \z$, $n\in \nat$ being relatively
prime, then
$$
Q_{\infty}(a)\subset \langle{q}\rangle_{\Lambda_a}+{\cal
O}\left(\frac{1}{k(\omega(a))^{\infty}}\right)[-1,1],\quad
k(\omega(a)):=\abs{m}+\abs{n}.
$$
In particular,
$$
\sum_{a;\, \omega(a)\in {\rm \bf Q}} \abs{Q_{\infty}(a)}<\infty.
$$
Here $\abs{I}$ stands for the length of an interval $I\subset
\real$.
\end{itemize}
\end{prop}

When $a\in S$, we also introduce the compact interval
$Q_{\infty}(a)\subset \real$, defined as in (\ref{1.25}). We shall
assume from now on that the function
$\langle{q}\rangle_{\Lambda_a}$ depends analytically on $a\in
J\backslash S$ and that it extends to a continuous function on the
whole interval $J$. We shall also assume that $a\mapsto
\langle{q}\rangle_{\Lambda_a}$ is not identically constant on any
connected component of $J\backslash S$.

 \vskip 2mm Let us next consider the behavior of the interval $Q_{\infty}(a)$
corresponding to a rational torus $\Lambda_a$ in a \neigh{} of a
Diophantine torus. In doing so, we assume that $\Lambda_{a_0}$,
$a_0\in J\backslash S$, is such that $\omega_0:=\omega(a_0)$ is
$(\alpha,d)$--Diophantine for some $\alpha>0$ and $d>0$, i.e. it
satisfies the usual Diophantine condition---cf. (\ref{0.13}),
\begin{equation}
\label{7.2}
 \abs{\omega_0-\frac{p}{q}}\geq
\frac{\alpha}{q^{2+d}},\quad p\in \z,\,\,q\in \nat.
\end{equation}
Let $\omega=\omega(a)$, $a\in {\rm neigh}(a_0,\real)$, be
rational, and write $\omega=m/n$ where $m\in \z$ and $n\in \nat$
are relatively prime. Then
$$
0=\omega-\frac{m}{n}=\omega-\omega_0+\omega_0-\frac{m}{n},
$$
and applying (\ref{7.2}) with $\frac{p}{q}=\frac{m}{n}$, we get
$$
\abs{\omega-\omega_0}\geq \frac{\alpha}{n^{2+d}}.
$$
Therefore,
$$
\frac{1}{k(\omega)}\leq {\cal
O}(1)\abs{\omega-\omega_0}^{1/(2+d)},
$$
and an application of Proposition 7.1 shows that $Q_{\infty}(a)$
is a closed interval such that
\begin{equation}
\label{7.3} Q_{\infty}(a)\subset \langle{q}\rangle_{\Lambda_a}+{\cal
O}\left(\abs{\omega(a)-\omega_0}^{\infty}\right)[-1,1].
\end{equation}
This result should be compared with (\ref{1.29.5}).

\vskip 4mm We shall next consider the set of Lagrangian tori in
$p^{-1}(0)$, whose rotation numbers are uniformly Diophantine and
that satisfy a uniform isoenergetic condition. Put for $\alpha>0$
and $d>0$ fixed, \begeq \label{7.3.5}
\Omega_{\alpha,d}=\left\{a\in J;\,\,{\rm dist}(a,S)\geq
\alpha,\,\, \abs{\omega'(a)}\geq \alpha,\,\,
\abs{\omega(a)-\frac{p}{q}}\geq \frac{\alpha}{q^{2+d}},\,\,\, p\in
\z,\,\,\, q\in \nat\right\}.
\endeq
In what follows, the dependence on $d$ will not be indicated
explicitly, and we shall write
$\Omega_{\alpha}=\Omega_{\alpha,d}$. Introduce next
\begin{equation}
\label{7.4}
 \widetilde{\Omega}_{\alpha}:=\Omega_{\alpha}\cap \{ a\in J\backslash S,
 \abs{d_a\langle{q}\rangle_{\Lambda_a}}\geq \alpha\},
\end{equation}
and notice that the measure of the complement of
$\widetilde{\Omega}_{\alpha}$ in the interval $J$, $\complement
\widetilde{\Omega}_{\alpha}$, is small, together with $\alpha$,
when $d>0$ is kept fixed. We then define the set of good values
${\cal G}_{\alpha}$ contained inside the closed interval
$\langle{q}\rangle_{\Lambda_a}\left(J\right)$, in the following
way,
\begin{equation}
\label{G}
 {\cal G}_{\alpha}=\complement
\left(\langle{q}\rangle_{\Lambda_a}\left(\complement
\widetilde{\Omega}_{\alpha}\right)\right).
\end{equation}
It follows that the complement of ${\cal G}_{\alpha}$ in
$\langle{q}\rangle_{\Lambda_a}(J)$ has a small measure, when
$\alpha$ is small and $d>0$ is kept fixed. We notice also that it
follows from the construction that when $F_0\in {\cal G}_{\alpha}$,
the pre-image
$$
\langle{q}\rangle_{\Lambda_a}^{-1}(F_0)
$$
is a finite set,
$$
\langle{q}\rangle_{\Lambda_a}^{-1}(F_0)=\left\{a_1,\ldots\, ,
a_L;\,a_j\in \widetilde{\Omega}_{\alpha}\right\},
$$
and an application of (\ref{7.3}) shows that
$$ {\rm
dist}\left(Q_{\infty}(a),F_0\right)\geq \frac{\abs{a-a_j}}{{\cal
O}(1)},\quad a\in {\rm neigh}(a_j,\real),
$$
for any $j=1,\ldots L$.

When considering the intervals $Q_{\infty}(a)$ for $a\in J$ away from the
$a_j$, $j=1,\ldots L$, we notice that an application of Lemma 2.4
shows that the set
$$
\bigcup_{a\notin {\rm neigh}(\{a_1,\ldots a_L\})} Q_{\infty}(a)
$$
is closed. Our final assumption now is that
\begeq
\label{7.4.1}
\sup_{a\in J} \abs{Q_{\infty}(a)}\,\,\,\wrtext{is
sufficiently small depending on}\,\,\,\alpha\,\,\wrtext{and}\,\,d.
\endeq
When $a$ is away from the $a_j$, $j=1,\ldots L$, using (\ref{7.4.1}), we get
$$
{\rm dist}\left(Q_{\infty}(a),F_0\right)\geq \frac{1}{{\cal
O}(1)}.
$$

We summarize the discussion above in the following theorem.

\begin{theo}
Assume that the $H_p$-flow on $T^*M$ is completely integrable and
assume that the rotation number $\omega(a)$ of the invariant tori
$\Lambda_a$ depends analytically on $a$ and is not identically
constant on any connected component of $J\backslash S$. When $q$
is a real-valued bounded analytic function on $T^*M$, define
$\langle{q}\rangle_{\Lambda_a}$, $a\in J\backslash S$, to be the
torus average of $q$ with respect to the natural smooth measure on
$\Lambda_a$, and assume that $a\mapsto
\langle{q}\rangle_{\Lambda_a}$ is an analytic function on
$J\backslash S$, which extends continuously to $J$ and which is
not identically constant on any connected component of
$J\backslash S$. When $\alpha>0$, $d>0$, let us define next the
set ${\cal G}_{\alpha}\subset \langle{q}\rangle_{\Lambda_a}(J)$
according to {\rm (\ref{G})}, {\rm (\ref{7.4})}. Then the measure
of the complement of ${\cal G}_{\alpha}$ in
$\langle{q}\rangle_{\Lambda_a}(J)$ is small when $\alpha$ is
sufficiently small and $d$ is kept fixed. Assume that
$$
\sup_{a\in J} \abs{Q_{\infty}(a)}
$$
is small enough depending on $\alpha$, $d$. When $F_0\in {\cal
G}_{\alpha}$, we have that for any \neigh{} $W$ of the finite set
$\langle{q}\rangle_{\Lambda_a}^{-1}(F_0)$ there exists a constant
$C(W)>0$ such that
$$
{\rm dist}\left(Q_{\infty}(a),F_0\right) \geq \frac{1}{C(W)},\quad
a\in \complement W.
$$
Here $Q_{\infty}(a)$, $a\in J$, is defined as in {\rm
(\ref{1.25})}, so that it is equal to the range of the limit of
$\langle{q}\rangle_T$, as $T\rightarrow \infty$, along
$\Lambda_a$, when $a\in J\backslash S$.
\end{theo}

\subsection{Surfaces of revolution.} We shall now illustrate Theorem 7.2
in the case when $M$ is an analytic surface of revolution in
$\real^3$. In doing so, we normalize $M$ so that the $x_3$-axis is
its axis of revolution, and we parametrize it by the cylinder
$[0,L]\times S^1$,
$$
[0,L]\times S^1 \ni (s,\theta)\mapsto (f(s)\cos\theta, f(s)\sin
\theta, h(s)),
$$
assuming, as we may, that the parameter $s$ is the arc-length
along the meridians, so that $(f'(s))^2+(h'(s))^2=1$. A simple
calculation then shows that in the coordinates $(s,\theta)$, the
metric on $M$ takes the form
\begin{equation}
\label{7.5}
 g=ds^2+f^2(s)d\theta^2.
\end{equation}
The functions $f$ and $h$ are assumed to be real analytic on
$[0,L]$, and from~\cite{Be} we recall that the regularity of $M$
at the poles is guaranteed by requiring that for each $k\in \nat$,
$$
f^{(2k)}(0)=f^{(2k)}(L)=0,
$$
and that $f'(0)=1$, $f'(L)=-1$, which we shall assume from now
on.

Throughout the following discussion, we shall assume furthermore
that $M$ is a simple surface of revolution, in the sense that
$0\leq f(s)$ has precisely one critical point $s_0\in (0,L)$, and
that this critical point is a non-degenerate maximum,
$f''(s_0)<0$. To fix the ideas, we shall also assume that
$f(s_0)=1$. For future reference, we notice that $s_0$ corresponds
to the equatorial geodesic $\gamma_E\subset M$ given by $s=s_0$,
$\theta\in S^1$. This is an elliptic orbit.

We now come to analyze the geodesic flow on $M$, which we shall
view as the Hamilton flow on $T^*M$ of the dual form to the metric
$g$. In doing so, we write
$$
T^*\left(M\backslash \{(0,0,g(0)),(0,0,g(L))\}\right)\simeq
T^*\left((0,L)\times S^1\right),
$$
and using (\ref{7.5}) we immediately see that the dual form to $g$
is given by
\begin{equation}
\label{7.6}
p(s,\theta,\sigma,\theta^*)=\sigma^2+\frac{(\theta^*)^2}{f^2(s)}.
\end{equation}
 Here $\sigma$ and $\theta^*$ are the dual variables to $s$ and
$\theta$, respectively. Writing out the Hamilton equations for the
$H_p$--flow, we see next that the functions $p$ and $\theta^*$ are
in involution, and we recover the well-known fact that the
geodesic flow on $M$ is completely integrable.

When $E>0$ and $\abs{F}< E^{1/2}$, $F\neq 0$, it is true that the
set defined by
$$
\Lambda_{E,F}: p=E,\,\, \theta^*=F,
$$
is an analytic Lagrangian torus sitting inside the energy surface
$p^{-1}(E)$. Geometrically, the torus $\Lambda_{E,F}$ consists of
geodesics contained between and intersecting tangentially the
parallels $s_{\pm}(E,F)$ on $M$ defined by the equation
$$
f(s_{\pm}(E,F))=\frac{\abs{F}}{E^{1/2}}.
$$
For $F=0$, the parallels reduce to the two poles and we obtain a
torus consisting of a family of meridians. The case
$\abs{F}=E^{1/2}$ is degenerate and corresponds to the equator
$s=s_0$, traversed with the two different orientations.

\vskip 2mm The principal actions of the torus $\Lambda_{E,F}$ are
given by \begeq \label{7.7}
I_1(E,F)=\int_{\sigma^2+\frac{F^2}{f^2(s)}=E}
\sigma\,ds=2\int_{s_-(E,F)}^{s_+(E,F)}
\left(E-\frac{F^2}{f^2(s)}\right)^{1/2}\,ds,
\endeq
and \begeq \label{7.8} I_2(E,F)=2\pi F.
\endeq
The functions $I_1$, $I_2$ depend analytically on $E$, $F\neq 0$.

Let us now restrict the further discussion to the energy surface
$p^{-1}(1)$. We shall derive an explicit expression for the
rotation number $\omega(a)$ of the torus
$\Lambda_a:=\Lambda_{1,a}$, $0\neq a\in (-1,1)$. In doing so, we
introduce the action coordinates $\eta_1=I_1/2\pi$, $\eta_2=I_2/2\pi$,
and notice that it follows from (\ref{7.7}) and (\ref{7.8}) that along
$p^{-1}(1)$, $\eta_1$ becomes an analytic function of $\eta_2$,
$$
\eta_1=\varphi(\eta_2):=\frac{1}{\pi}\int_{s_-(a)}^{s_+(a)}
\left(1-\frac{\eta_2^2}{f^2(s)}\right)^{1/2}\,ds,\quad \eta_2=a,\quad
s_{\pm}(a):=s_{\pm}(1,a).
$$
Differentiating next the relation $p(\varphi(\eta_2),\eta_2)=1$,
we get that the rotation number $\omega(a)$ is given by
$-\varphi'(a)$, so that
\begin{equation}
\label{7.9}
\omega(a)=\frac{a}{\pi}\int_{s_-(a)}^{s_+(a)}\frac{1}{f^2(s)}\left(1-\frac{a^2}{f^2(s)}\right)^{-1/2}\,ds.
\end{equation}
We recall that the torus $\Lambda_a\subset p^{-1}(1)$ satisfies
the isoenergetic condition precisely when the expression
$\omega'(a)$, given by
\begin{equation}
\label{7.10}
\frac{1}{\pi}\int_{s_(a)}^{s_+(a)}\frac{1}{f^2(s)}\left(1-\frac{a^2}{f^2(s)}\right)^{-1/2}\,ds+
\frac{a}{\pi}\frac{\partial}{\partial
a}\left(\int_{s_-(a)}^{s_+(a)}\frac{1}{f^2(s)}\left(1-\frac{a^2}{f^2(s)}\right)^{-1/2}\,ds\right)
\end{equation}
is non-vanishing.

From now on, we shall assume that the surface of revolution $M$ is
such that $\omega(a)$ in (\ref{7.10}) is not identically constant.
As explained in~\cite{Z}, examples of such surfaces are provided
by ellipsoids of revolution, with $\omega'(a)<0$ for oblong
ellipsoids, and $\omega'(a)>0$ for oblate ones. The separating
case of the sphere is of course degenerate with $\omega'(a)\equiv
0$. More generally, we may also remark that the isoenergetic
condition rules out the case of analytic Zoll surfaces of
revolution.

Let now $q$ be an analytic function on $M$, which we shall view as
a function on $T^*M$. We now come to discuss the long time
averages of $q$ along the $H_p$-flow, $\langle{q}\rangle_T$, $T\gg
1$. In doing so, we shall first consider the case when
$q=q_0=q_0(s)$ is a function of $s$ only. It follows from
(\ref{7.6}) that along an $H_p$-orbit in $\Lambda_a$, $0\neq a\in
(-1,1)$, we have
\begin{equation}
\label{7.11}
 \dot{s}(t)=\pm 2 \left(1-\frac{a^2}{f^2(s(t))}\right)^{1/2},
\end{equation}
and then using the rotational invariance of $q_0$, we get
$$
\langle{q_0}\rangle_T(m)=\frac{1}{T}\int_0^T
q_0(s(t))\,dt=\langle{q_0}\rangle_{\Lambda_a}+{\cal
O}\left(\frac{1}{T}\right),\quad m\in \Lambda_a, \quad
T\rightarrow \infty.
$$
Here, as we immediately compute from (\ref{7.11}),
\begin{equation}
\label{7.12}
\langle{q_0}\rangle_{\Lambda_a}=\frac{J(q_0,a)}{J(1,a)},
\end{equation}
where, in general, for an analytic function $\psi$ we write
\begin{equation}
\label{7.13}
 J(\psi,a)=\int_{s_-(a)}^{s_+(a)}
\psi(s)\frac{f(s)}{(f^2(s)-a^2)^{1/2}}\,ds,\quad
f(s_{\pm}(a))=\abs{a}.
\end{equation}
We remark that (\ref{7.12}), (\ref{7.13}) provide an explicit
description of $\langle{q_0}\rangle_{\Lambda_a}$ which, as before,
is defined as the average of $q_0$ over the invariant torus
$\Lambda_a$.

The function $\langle{q_0}\rangle_{\Lambda_a}$ depends
analytically on $a\in (-1,1)$, and it extends to a continuous function
on the whole interval $[-1,1]$. In what follows we shall
assume that the function $a\mapsto
\langle{q_0}\rangle_{\Lambda_a}$ is not identically constant.

\medskip
{\it Example}. Let
$$
q_0(s)=(f'(s))^2=(f''(s_0))^2(s-s_0)^2+{\cal O}((s-s_0)^3),\quad
s\rightarrow s_0.
$$
Since $q_0$ is minimal along the equator $s=s_0$, it is
clear that $\langle{q_0}\rangle_{\Lambda_a}$ is not identically
constant and achieves its minimum for $a=\pm 1$. Here we shall compute
a leading term in the asymptotic expansion of this function, as
$\abs{a}\to 1$. In doing so, following~\cite{Be} and~\cite{Z},
it will be convenient to make a change of variables on the
surface $M$. We introduce
$$
r=\cases{\arcsin f(s),\,\, s\in (0,s_0) \cr
        \pi-\arcsin f(s),\,\, s\in (s_0,L).}
$$
It follows then that $s=c(\cos r)$, where $c(x)$, $x\in [-1,1]$, is
given by
$$
c(x)=\cases{ (f|_{[0,s_0]})^{-1}(\sqrt{1-x^2}),\,\, x\in [0,1], \cr
             (f|_{[s_0,\pi]})^{-1}(\sqrt{1-x^2}),\,\, x\in [-1,0].}
$$
Therefore, $\sin r=f(s)$, $f'(c(\cos r))ds=\cos r dr$, and making
the change of variables in (\ref{7.13}), we get \begeq
\label{7.14} J(\psi,a)=\int_{\arcsin \abs{a}}^{\pi-\arcsin
\abs{a}} \psi(c(\cos r)) \frac{\sin r h(\cos r)}{(\sin^2
r-a^2)^{1/2}}\,dr=\int_{-\sqrt{1-a^2}}^{\sqrt{1-a^2}}\frac{\psi(c(x))h(x)}{(1-x^2-a^2)^{1/2}}\,dx.
\endeq
Here $h(x)=x/f'(c(x))$, $x\neq 0$, $h(0)=1/{\abs{f''(s_0)}}^{1/2}\neq 0$.

Taylor expanding $h$, $h(x)=h_0+h_1 x+h_2x^2+\ldots$, $
x\rightarrow 0$, we get after a simple computation,
$$
J(1,a)=\pi h_0+\frac{\pi}{2}h_2(1-a^2)+{\cal O}((1-a^2)^2),\quad
\abs{a}\rightarrow 1.
$$
Now
$$
q_0(c(x))=(f'(c(x)))^2=\frac{x^2}{h(x)^2}=\frac{1}{h_0^2}x^2+{\cal
O}(x^3),
$$
and using (\ref{7.14}), we get
$$
J(q_0,a)=\frac{\pi}{2h_0}(1-a^2)+{\cal O}((1-a^2)^2).
$$
It follows that \begeq \label{7.15}
\langle{q_0}\rangle_{\Lambda_a}=\frac{J(q_0,a)}{J(1,a)}=\frac{1}{2}(1-a^2)+{\cal
O}((1-a^2)^2)
\endeq
is not identically constant.

\Remark. The preceding example is closely related to the computations
in Appendix C in~\cite{SjCJM}.

\medskip
We shall now introduce a weak angular dependence in the
perturbation $q$, and when doing so, for an analytic function
$q_1$ on $M$, we set \begeq \label{7.15.5}
q_{\eta}(s,\theta)=q_0(s)+\eta q_1(s,\theta),\quad 0<\eta \ll 1.
\endeq
As before, we are then interested in the limit
\begin{equation}
\label{7.16}
 \lim_{T\rightarrow \infty} \frac{1}{T}\int_0^T q_{\eta}\circ
\exp(tH_p)\,dt,
\end{equation}
taken along the invariant torus $\Lambda_a$. When analyzing
(\ref{7.16}), we again switch to the action-angle variables, so
that $\Lambda_a$ becomes the standard torus ${\bf T}^2$, $H_p$
along $\Lambda_a$ becomes $b_1\partial_{x_1}+b_2\partial_{x_2}$,
$b_j \in \real$, $b_1\neq 0$, and
$$
q_{\eta}(x_1,x_2)=q_1(x_1)+\eta q_2(x_1,x_2).
$$
It follows from (\ref{7.1.2}) that (\ref{7.16}) is equal to
$$
\widehat{q}_1(0)+\eta \sum_{b\cdot k=0} \widehat{q}_2(k)e^{ix
\cdot k}.
$$

Let us summarize this discussion in the following proposition,
which is just a specialization of Proposition 7.1 to the case of
surfaces of revolution, with a perturbation of the form
(\ref{7.15.5}).

\begin{prop}
On the surface of revolution $M$, let us consider an analytic
function $q_{\eta}(s,\theta)=q_0(s)+\eta q_1(s,\theta)$, $0<\eta
\ll 1$. When $a\in (-1,1)$, we define the closed interval
$Q_{\infty,\eta}(a)\subset \real$ to be the range of the limit of
$\langle{q_{\eta}}\rangle_T$, as $T\rightarrow \infty$, restricted
to the invariant torus $\Lambda_a\subset p^{-1}(1)$. We then have
\begin{itemize}
\item When the rotation number $\omega(a)\notin {\bf Q}$, then
$Q_{\infty,\eta}(a)=\{ \langle{q_{\eta}}\rangle_{\Lambda_a} \}$.
\item When $\omega(a)=m/n$, with $m\in \z$, $n\in \nat$ being
relatively prime, then
$$
Q_{\infty,\eta}(a)\subset \langle{q_{\eta}}\rangle_{\Lambda_a}+\eta {\cal
O}\left(\frac{1}{k(\omega(a))^{\infty}}\right)[-1,1],\quad
k(\omega(a)):=\abs{m}+\abs{n}.
$$
\end{itemize}
\end{prop}

\vskip 4mm Let us now define the set $\Omega_{\alpha}$ as in
(\ref{7.3.5}), with $J=[-1,1]$ and $S=\{\pm 1\}$, and define
following (\ref{7.4}),
\begin{equation}
\label{7.17}
 \widetilde{\Omega}_{0,\alpha}:=\Omega_{\alpha}\cap \{ a\in
 (-1,1)\backslash\{0\};
 \abs{d_a\langle{q_0}\rangle_{\Lambda_a}}\geq \alpha\}.
\end{equation}
We then define the $\eta$-dependent set of good values ${\cal
G}_{\alpha,\eta}$ contained inside the closed interval
$\langle{q_{\eta}}\rangle_{\Lambda_a}\left([-1,1]\right)$, as in
(\ref{G}),
\begin{equation}
\label{G1}
 {\cal G}_{\alpha,\eta}=\complement
\left(\left(\langle{q_0}\rangle_{\Lambda_a}+\eta\langle{q_1}\rangle_{\Lambda_a}\right)\left(\complement
\widetilde{\Omega}_{0,\alpha}\right)\right).
\end{equation}
The complement of ${\cal G}_{\alpha,\eta}$ has a small measure, when
$\alpha$ is small, uniformly in $\eta$ small enough, and when
$d>0$ is kept fixed.

The general assumptions of Theorem 7.2 are verified in this case,
and we get the following result.

\begin{prop}
Assume that $M$ is a simple analytic surface of revolution for which the
rotation number defined in {\rm (\ref{7.9})} is not identically
constant. When $q_0=q_0(s)$ is an analytic function on $M$, let us
define $\langle{q_0}\rangle_{\Lambda_a}$ as in {\rm (\ref{7.12})},
{\rm (\ref{7.13})}, and assume that $a\mapsto
\langle{q_0}\rangle_{\Lambda_a}$ is not identically constant. Set
when $0\leq \eta \ll 1$,
$$
q_{\eta}(s,\theta)=q_0(s)+\eta q_1(s,\theta),
$$
and define $\langle{q_{\eta}}\rangle_{\Lambda_a}$ to be the
average of $q_{\eta}$ over the invariant torus $\Lambda_a$, with
respect to the natural smooth measure. When $\alpha>0$, let us
define the set ${\cal G}_{\alpha,\eta}\subset
\langle{q_{\eta}}\rangle_{\Lambda_a}\left([-1,1]\right)$ according
to {\rm (\ref{G1})}. Then the measure of the complement of ${\cal
G}_{\alpha,\eta}$ in
$\langle{q_{\eta}}\rangle_{\Lambda_a}\left([-1,1]\right)$ is small
when $\alpha$ and $\eta$ are sufficiently small (and $d$ is kept
fixed). When $F_0\in {\cal G}_{\alpha,\eta}$ and $\eta$ is small
enough, we have that for any \neigh{} $W$ of the finite set
$\langle{q_{\eta}}\rangle_{\Lambda_a}^{-1}(F_0)$ there exists a
constant $C(W)>0$ such that
$$
{\rm dist}\left(Q_{\infty,\eta}(a),F_0\right) \geq
\frac{1}{C(W)},\quad a\in \complement W.
$$
Here $Q_{\infty,\eta}(a)$ is the range of the limit of
$\langle{q_{\eta}}\rangle_T$, as $T\rightarrow \infty$, along
$\Lambda_a$.
\end{prop}

It follows from Proposition 7.4 that for $F_0\in {\cal
G}_{\alpha,\eta}$, the assumptions of Theorems 1.1 and 1.2 (see also
Theorems 5.2 and 6.1) are satisfied. An application of, say,
Theorem 1.1 gives then the complete asymptotic expansions of all
the eigenvalues of $-h^2\Delta +i\eps q_{\eta}$ in rectangles in
the spectral complex plane, that are of the form
$$
\abs{\Re z-1}\leq \frac{h^{\delta}}{{\cal O}(1)},\quad \abs{\Im
z-\eps \Re F_0}\leq \frac{\eps h^{\delta}}{{\cal O}(1)},\quad
\delta>0,
$$
when $\eps={\cal O}(h^{\delta})$ is bounded from below by a fixed
positive power of $h$. We refrain from formulating precisely the
final result, as it is immediately obtained from the statement of
Theorem 1.1.

\subsection{Perturbations of completely integrable systems.} We
now come to discuss the perturbed situation and before doing so,
we shall pause to recall the statement of the isoenergetic KAM
theorem.

We consider
 \begeq \label{kam1} p_{\lambda}=p_0+\lambda p_1,\quad 0<\lambda
\ll 1,
\endeq
where $p_0$ and $p_1$ are holomorphic bounded functions in a
tubular \neigh{} of $T^*M$, that are real on the real domain. We
assume that the $H_{p_0}$-flow is completely integrable, and let
$\Lambda_0\subset p_0^{-1}(0)$ be an analytic $H_{p_0}$-invariant
Diophantine Lagrangian torus. Take then the analytic action-angle
transformation \begeq \label{kam2} \kappa: {\rm neigh}(\Lambda_0,
T^*M)\rightarrow {\rm neigh}(\xi=0, T^*{\bf T}^2),
\endeq
such that $\kappa(\Lambda_0)$ is the zero section in $T^*{\bf
T}^2$, and
$$
p_0\circ \kappa^{-1}=p_0(\xi)=a\cdot \xi+{\cal O}(\xi^2),\quad
a=(a_1,a_2)\in \real^2,
$$
with $a$ satisfying the Diophantine condition (\ref{0.13}), where we
assume that $N_0>1$.
Composing $p_{\lambda}$ in (\ref{kam1})
with the inverse of $\kappa$ in (\ref{kam2}), we get a function
$$
p_{\lambda}(x,\xi)=p_0(\xi)+\lambda p_1(x,\xi),
$$
analytic in a fixed complex \neigh{} of $\xi=0$. An application of
the implicit function theorem shows that the energy surface
$p_0^{-1}(0)$ takes the form $\xi_2=f(\xi_1)$, where $f$ is an
analytic function near $0\in \real$, with $f(0)=0$,
$f'(0)=-a_1/a_2$. We assume that the isoenergetic condition holds:
\begeq \label{kam3} f''(0)\neq 0.
\endeq
The condition (\ref{kam3}) means that for $\abs{\mu}\ll 1$, the
$H_{p_0}$-invariant tori $\Lambda_{\mu}: \xi_1=\mu$,
$\xi_2=f(\mu)$ can be parametrized by the corresponding rotation
numbers $f'(\mu)$.

\begin{theo}(The isoenergetic KAM theorem).
Make the isoenergetic assumption {\rm (\ref{kam3})} and let us define
for $\alpha>0$ sufficiently small and $d>0$ fixed,
$$
K_{\alpha}=\left\{\mu \in {\rm neigh}(0,\real);\,\, \abs{f'(\mu)-\frac{p}{q}}\geq
\frac{\alpha}{q^{2+d}},\,\,\, p\in \z,\,\,\, q\in \nat \right\}.
$$
Assume that $0<\lambda\ll \alpha^2$. Then there exists a $C^{\infty}$-map \begeq \label{kam3.5}
\Psi_{\lambda}: {\bf T}^2\times {\rm neigh}(0,\real)\rightarrow
{\rm neigh}(\xi=0,T^*{\bf T}^2)
\endeq
analytic in the angular variable (uniformly with respect to the
other factor), depending analytically on $\lambda>0$, and with
${\rm neigh}(0,\real)$ and ${\rm neigh}(\xi=0,T^*{\bf T}^2)$ in
{\rm (\ref{kam3.5})} uniform in $\lambda$, such that for each $\mu
\in K_{\alpha}$, the set \begeq \label{kam4}
\Lambda_{\mu,\lambda}=\{\Psi_{\lambda}(x,\mu); x\in {\bf
T}^2\}\subset T^*{\bf T}^2
\endeq
is a uniformly analytic Lagrangian torus $\subset
p_{\lambda}^{-1}(0)$, close to $\Lambda_{\mu}\subset p_0^{-1}(0)$.
The restriction of $\exp(tH_{p_{\lambda}})$ to
$\Lambda_{\mu,\lambda}$ is conjugated to the flow of a constant
vector field with the rotation number $f'(\mu)$ on ${\bf T}^2$, by
means of an analytic diffeomorphism, depending analytically on
$\lambda$ and smoothly on $\mu\in {\rm neigh}(0,\real)$. Moreover,
relative to a sufficiently small \neigh{} of $\xi=0$ in $T^*{\bf
T}^2$, the Liouville measure of the complement of the union of the
invariant tori $\Lambda_{\mu,\lambda}$, $\mu\in K_{\alpha}$, in
$p_{\lambda}^{-1}(0)$, is ${\cal O}(\alpha)$.
\end{theo}

\Remark. The classical theorem of Kolmogorov on the existence of an
invariant torus for the perturbed Hamiltonian, close to any given
Diophantine invariant torus for the initial integrable system, is
stated and proved in great detail in~\cite{BGGS}. The isoenergetic
version of the Kolmogorov theorem is stated explicitly in~\cite{Bost},
but without a proof, and a direct argument showing how to derive the
isoenergetic version of the theorem from the usual one is described
in~\cite{Douady}. The paper~\cite{Po} gives a complete proof of the
full KAM theorem and establishes, in particular, the
$C^{\infty}$--dependence of the KAM tori on the frequencies, in the
sense of Whitney. The
statement of Theorem 7.5 above can be extracted from the
references~\cite{BC} (see, in particular, Lemma 2.3 of that paper),
and especially,~\cite{BH}. The latter paper derives
Theorem 7.5 from the ordinary KAM theorem of~\cite{Po}. Finally, we
would also like to mention an interesting article~\cite{DG}, which
provides a complete proof of the isoenergetic KAM theorem. It seems, however,
that only a continuous, rather than a Whitney--smooth, dependence of
the KAM tori on the rotation numbers, is obtained in that work.

\Remark. A more refined version of the isoenergetic KAM
theorem, which also follows from the arguments of~\cite{BH}, gives
Cantor families of Diophantine invariant tori with prescribed rotation
numbers, in energy surfaces $p_{\lambda}^{-1}(E)$, when $E\in {\rm
neigh}(0,\real)$, with an analytic dependence on $E$.

\medskip
We shall apply Theorem 7.5 to the situation described in the
beginning of this section, with $p_0=p$. Let us recall that the
energy surface $p^{-1}(0)$ is foliated by the Lagrangian tori
$\Lambda_a$, $a\in J\backslash S$, and that the rotation number of
$\Lambda_a$, $\omega(a)$, is assumed to be not identically
constant on any connected component of $J\backslash S$.
As before, let $q$ be a real-valued bounded analytic
function on $T^*M$, and define the function
$\langle{q}\rangle_{\Lambda_a}$ as the average of $q$ over the
$H_p$--invariant tori $\Lambda_a$, $a\in J\backslash S$. Assume
that the assumptions of Theorem 7.2 are satisfied and let $F_0\in
{\cal G}_{\alpha}$, with the latter set defined in (\ref{7.9}),
(\ref{G}). Then
$$
\langle{q}\rangle_{\Lambda_a}^{-1}(F_0)=\{a_1,\ldots\, , a_L; a_j\in
\widetilde{\Omega}_{\alpha}\},
$$
and
\begeq
\label{kam4.1}
{\rm dist}(Q_{\infty}(a),F_0)\geq 1/{\cal O}(1),\,\,
\wrtext{when}\,\,\,a\,\,\wrtext{is away from the}\,\, a_j,\,\,
j=1,\ldots\,, L.
\endeq
It follows from the equivalence of the assumptions (\ref{H0}) and
(\ref{1.31}), noticed in section 2, together with Lemma 2.2, that for
any sufficiently small $H_p$--invariant $\delta$-\neigh{} $W_{\delta}$
of $\cup_{j=1}^L \Lambda_{a_j}$ in $T^*M$ we have, for $\rho \in
p^{-1}((-\widetilde{\delta}(\delta),\widetilde{\delta}(\delta)))\backslash
W_{\delta}$, provided that $\widetilde{\delta}(\delta)>0$ is
small enough,
\begeq
\label{kam4.2}
\abs{\langle{q}\rangle_{\delta^{-N_1},K,p}(\rho)-F_0}\geq
\frac{\delta}{{\cal O}(1)},\quad N_1\in \nat\backslash \{0\}.
\endeq

Let now $p_1$ be a bounded analytic and real-valued function on
$T^*M$, and let us consider
$$
p_{\lambda}=p+\lambda p_1,\quad \lambda>0.
$$
An application of Theorem 7.5 shows that for $\lambda/\alpha^2$
small enough, for each $1\leq j\leq L$, there exists a smooth
family of tori $\Lambda^{(j)}_{a,\lambda}$, $a\in {\rm
neigh}(a_j,\real)$, close to $\Lambda_{a_j}$, such that if $a\in
\Omega_{\alpha,j}$, where
$$
\Omega_{\alpha,j}=\left\{ a\in {\rm
neigh}(a_j,\real),\,\,\abs{\omega(a)-\frac{p}{q}}\geq
\frac{\alpha}{2 q^{2+d}},\, p\in \z,\, q\in \nat\right \},
$$
then $\Lambda^{(j)}_{a,\lambda}\subset
p_{\lambda}^{-1}(0)$ is a Lagrangian torus, on which the
$H_{p_{\lambda}}$--flow is quasi-periodic, with the rotation number
$\omega(a)$. As before, we may remark that the measure of the set ${\rm
neigh}(a_j,\real)\backslash \Omega_{\alpha,j}$ is small when
$\alpha>0$ is small.

When $1\leq j\leq L$, we define a smooth function
of $a\in {\rm neigh}(a_j,\real)$,
$\langle{q}\rangle_{\Lambda^{(j)}_{a,\lambda}}$,
obtained by averaging $q$ over the tori $\Lambda^{(j)}_{a,\lambda}$.
Then we have, in the $C^1$-sense, as $\lambda\rightarrow 0$,
$$
\langle{q}\rangle_{\Lambda^{(j)}_{a,\lambda}}\rightarrow
\langle{q}\rangle_{\Lambda^{(j)}_{a,\lambda=0}},
$$
and therefore, when $a\in {\rm neigh}(a_j,\real)$, we get
$$
\abs{d_a\langle{q}\rangle_{\Lambda^{(j)}_{a,\lambda}}}\geq
\frac{\alpha}{2},\quad 1\leq j\leq L,
$$
provided that $\lambda$ is small enough. Following (\ref{G}), let
us define \begeq \label{E} {\cal E}_{\alpha,F_0}=\complement
\bigcup_{j=1}^L
\langle{q}\rangle_{\Lambda^{(j)}_{a,\lambda}}\left({\rm
neigh}(a_j,\real)\backslash \Omega_{\alpha,j}\right).
\endeq
Then the measure of the complement of ${\cal E}_{\alpha,F_0}$ is
small with $\alpha$, and if $F\in {\cal E}_{\alpha}$, then for
each $1\leq j\leq L$, there exists $b_j\in {\rm
neigh}(a_j,\real)$, such that $b_j\in \Omega_{\alpha,j}$ and
$\langle{q}\rangle_{\Lambda^{(j)}_{b_j,\lambda}}=F$, $1\leq j\leq
L$. Hence the sets $\Lambda^{(j)}_{b_j,\lambda}\subset
p_{\lambda}^{-1}(0)$, $1\leq j\leq L$, are uniformly Diophantine
Lagrangian tori, and the basic assumptions
(\ref{0.12})--(\ref{0.13}) are satisfied for the
$H_{p_{\lambda}}$--flow.

To be able to apply Theorems 1.1 and 1.2 to an operator satisfying the
general assumptions of the introduction and having the leading symbol
$p_{\eps,\lambda}=p+\lambda p_1+i\eps q+{\cal O}(\eps^2)$, we shall now discuss the
construction of a weight function $\widetilde{G}$ for
$p_{\eps,\lambda}$, satisfying the conclusion of Proposition
2.3. Outside a small \neigh{} of $\cup_{j=1}^L \Lambda^{(j)}_{b_j,\lambda}$, we shall take
$\widetilde{G}=G_{T,p}$, where $G_{T,p}=G_T$ is defined in
(\ref{0.16}), using the unperturbed symbol $p$. When working locally near a fixed $H_{p_{\lambda}}$--invariant
Diophantine torus $\Lambda^{(j)}_{b_j,\lambda}$, which we identify
with $\xi=0$ in $T^*{\bf T}^2$, we follow (\ref{1.13.4.5}) and set
\begeq
\label{kam5}
\widetilde{G}=(1-\chi_{\mu})G_{T,p}+\chi_{\mu} G.
\endeq
Here $G$ is defined as in (\ref{1.7.6}) and $0<\mu \ll 1$ is to be
chosen. An application of (\ref{1.13.6}) gives that uniformly in
$\mu>0$ we have
$$
q-H_{p_{\lambda}}
\widetilde{G}=\langle{q}\rangle+(1-\chi_{\mu})\left({\cal
O}(T\xi^N)+{\cal O}(\xi^{\infty}+T^{-\infty})\right)+{\cal
O}(T\xi^N)+{\cal O}_T(\lambda),
$$
provided that $\abs{T}\leq {\cal O}_N(1)\abs{\xi}^{-N}$.
Here we have also used that
$$
H_{p_{\lambda}}G_{T,p}=q-\langle{q}\rangle_{T,K,p_{\lambda}}+{\cal
O}_T(\lambda).
$$
As in section 2, we take first $N$ fixed but sufficiently large
depending on $N_1$. Choose then $\mu\sim \delta$, with $\delta>0$
sufficiently small but fixed as in (\ref{kam4.2}), and put
$T=\mu^{-N_1}$. Then in the region where $0<\chi_{\mu}<1$, we have
\begin{equation}
\label{kam5.1}
q-H_{p_{\lambda}}\widetilde{G}-F=\langle{q}\rangle(\xi)-\langle{q}\rangle(0)+{\cal
O}(\xi^{N-N_1})+{\cal O}_T(\lambda).
\end{equation}
On the other hand, away from the union of the tori
$\Lambda^{(j)}_{b_j,\lambda}$, we are in the region where
$\widetilde{G}=G_{T,p}$, and there we have, as $\lambda\rightarrow 0$,
\begin{equation}
\label{kam5.2}
q-H_{p_{\lambda}}\widetilde{G}-F=\langle{q}\rangle_{T,K,p}-F_0+{\cal
O}_T(\lambda)-o(1).
\end{equation}
Using (\ref{kam5.1}) and (\ref{kam5.2}), together with
(\ref{kam4.2}), we conclude that with this choice of
$\widetilde{G}$, the conclusion of Proposition 2.3 is satisfied
for $p_{\lambda}$, provided that $\lambda>0$ is taken sufficiently
small, depending on the parameters $N_1$ and $\delta$. Hence the
results of Theorem 1.1 and Theorem 1.2 apply in this perturbative
situation as they stand and give complete spectral results for an
operator satisfying the general assumptions of the introduction,
and with the leading symbol
$$
p_{\eps,\lambda}=p+\lambda p_1+i\eps q+{\cal O}(\eps^2),
$$
in rectangles of the form (\ref{R0}) and (\ref{R1}), for each $F\in {\cal E}_{\alpha,F_0}$.

\vskip 2mm
The discussion of this section is summarized in the following
theorem.

\begin{theo}
Let us keep all the general assumptions on the operator $P_{\eps}$
from the introduction, and write the leading symbol of $P_{\eps}$ as
$$
p_{\eps}=p+i\eps q +{\cal O}(\eps^2),
$$
near $p^{-1}(0)\cap T^*M$. Let us assume for simplicity that $q$ is
real-valued. We assume that the $H_p$--flow is
completely integrable, so that we have a decomposition {\rm
(\ref{7.0})}. When $a\in J\backslash S$, we let $\omega(a)$ stand
for the rotation number of the invariant torus $\Lambda_a$, and
assume that $a\mapsto \omega(a)$ is analytic and not identically
constant on any connected component of $J\backslash S$. When $a\in
J$, let us also introduce the compact sets $Q_{\infty}(a)\subset
\real$, defined as in {\rm (\ref{1.25})}, and recall from {\rm
(\ref{0.16.01})} that
$$
\frac{1}{\eps} \Im \left({\rm Spec}(P_{\eps})\cap \{z; \abs{\Re z}\leq \delta\}\right)
\subset \left [\inf \bigcup_{a\in J} Q_{\infty}(a)-o(1), \sup
\bigcup_{a\in J} Q_{\infty}(a)+o(1)\right],
$$
as $\eps$, $h$, $\delta\to 0$.

We then define a function $\langle{q}\rangle_{\Lambda_a}$, $a\in
J\backslash S$, obtained by averaging $q$ over the invariant tori
$\Lambda_a$, and assume that $\langle{q}\rangle_{\Lambda_a}$
depends analytically on $a\in J\backslash S$, and extends
continuously to $J$. Assume next that $a\mapsto
\langle{q}\rangle_{\Lambda_a}$ is not identically constant on any
connected component of $J\backslash S$. When $\alpha>0$ and $d>0$,
let ${\cal F}_{\alpha,d}\subset \cup_{a\in J} Q_{\infty}(a)$ be
the set obtained by removing from $\cup_{a\in J} Q_{\infty}(a)$
the following set
\begin{eqnarray*}
& & \left(\bigcup_{{\rm dist}(a,S)<\alpha} Q_{\infty}(a)\right)
\bigcup \left(\bigcup_{a\in J\backslash S,\,
\abs{\omega'(a)}<\alpha} Q_{\infty}(a)\right) \bigcup
\left(\bigcup_{a\in J\backslash
S,\,\abs{d\langle{q}\rangle_{\Lambda_a}}<\alpha}
Q_{\infty}(a)\right) \\
& & \bigcup \left(\bigcup_{a\in J\backslash S,\,\,
\omega(a)\,\wrtext{is
not}\,\,(\alpha,d)-\wrtext{Diophantine}}Q_{\infty}(a)\right).
\end{eqnarray*}
Then Theorem {\rm 1.1} applies to give complete asymptotic expansions
for all the eigenvalues of $P_{\eps}$ in a rectangle of the form
\begeq
\label{rect1}
\left[-\frac{h^{\delta}}{{\cal O}(1)}, \frac{h^{\delta}}{{\cal
O}(1)}\right]+i\eps \left[ F_0-\frac{h^{\delta}}{{\cal O}(1)},
F_0+\frac{h^{\delta}}{{\cal O}(1)}\right],\quad F_0\in {\cal
F}_{\alpha,d},
\endeq
when $\delta>0$ and $h^K\leq \eps={\cal O}(h^{\delta})$, $K\gg 1$. An
analogous statement is obtained from Theorem {\rm 1.2} in the case when
$h^{1/3-\delta}\leq \eps\ll 1$, $\delta>0$. The
measure of the complement of ${\cal F}_{\alpha,d}\subset \cup_{a\in J}
Q_{\infty}(a)$ is small, when $\alpha>0$ is small and $d$ is fixed,
provided that the measure of
$$
\left(\bigcup_{a\in J\backslash S, \omega(a)\in {\bf Q}}Q_{\infty}(a)\right)\bigcup
\left(\bigcup_{a\in S} Q_{\infty}(a)\right)
$$
is sufficiently small, depending on $\alpha$, $d$.

When $p_1$ is an analytic function in a tubular \neigh{} of
$T^*M$, real on the real domain, with $p_1(x,\xi)={\cal
O}(m(\Re(x,\xi))$ in the case when $M=\real^2$, and $p_1(x,\xi)={\cal
O}(\langle{\xi}\rangle^m)$ in the manifold case, let $P_{\eps,\lambda}$, $\lambda>0$, be an
operator satisfying all the general assumptions of the
introduction, and with the leading symbol
$$
p_{\eps,\lambda}=p_{\lambda}+i\eps q+{\cal O}(\eps^2),\quad
p_{\lambda}=p+\lambda p_1.
$$
When $F_0\in {\cal F}_{\alpha,d}$, let $\Lambda_{a_j}\subset
p^{-1}(0)$, $a_j\in J\backslash S$, $1\leq j\leq L$, be the
$(\alpha,d)$-Diophantine tori, with $\abs{\omega'(a_j)}\geq
\alpha$, and such that the average of $q$ along $\Lambda_{a_j}$ is equal
to $F_0$, $1\leq j\leq L$. When $\lambda/\alpha^2 \ll 1$, we introduce
a Cantor family of Diophantine Lagrangian tori $\Lambda^{(j)}_{a,\lambda}\subset
p_{\lambda}^{-1}(0)$, close to $\Lambda_{a_j}$, $1\leq j\leq L$,
whose existence is guaranteed by the isoenergetic KAM Theorem {\rm
7.5}. We then define a function
$a\mapsto \langle{q}\rangle_{\Lambda^{(j)}_{a,\lambda}}$, obtained by averaging
$q$ over the tori $\Lambda^{(j)}_{a,\lambda}$, and recall from Theorem
{\rm 7.5} that $\langle{q}\rangle_{\Lambda^{(j)}_{a,\lambda}}$ extends to a smooth
function on a full \neigh{} of $a_j$, $1\leq j\leq L$. Let ${\cal
E}_{\alpha,F_0}$ be the set obtained by removing from
$\cup_{j=1}^L \langle{q}\rangle_{\Lambda^{(j)}_{a,\lambda}}({\rm
neigh}(a_j,\real))$ the set given by
\begin{equation}
\bigcup_{j=1}^L
\langle{q}\rangle_{\Lambda^{(j)}_{a,\lambda}}\biggl(a\in {\rm
neigh}(a_j,\real);\, \omega(a)\,\,\wrtext{is not}\,\,
(\alpha/2,d)-\wrtext{Diophantine}\biggr).
\end{equation}
Then the measure of the complement of ${\cal E}_{\alpha,F_0}$ is
small when $\alpha$ is small and $d$ is kept fixed, and when $F\in
{\cal E}_{\alpha,F_0}$ and $h^K\leq \eps={\cal O}(h^{\delta})$,
$K\gg 1$, $\delta>0$, Theorem {\rm 1.1} applies and gives complete
asymptotic expansions for all the eigenvalues of
$P_{\eps,\lambda}$ in a rectangle of the form {\rm (\ref{rect1})},
with $F_0$ replaced by $F$, provided that $\lambda$ is
sufficiently small. A similar result holds by applying Theorem
{\rm 1.2} to $P_{\eps,\lambda}$, in the case when
$h^{1/3-\delta}\leq \eps\leq \eps_0$, $0<\eps_0\ll 1$.
\end{theo}

In the last theorem and elsewhere in the paper, we have tried to
find windows in the spectral band where all the \ev{}s of $P_{\eps}$ can be
described \asy{}ally. We shall finally discuss a reformulation of
Theorem 7.6 which permits us to describe a larger subset of the
spectral band where the conclusion of the theorem will hold
\ufly{}. At the same time, we shall let $J$ be a graph rather than
an interval (although this is not essential).
\medskip
\par\noindent \it The completely integrable case. \rm We make the same
assumptions of complete integrability as in Theorem 7.2. We also
assume that $p^{-1}(0)$ decomposes into a disjoint union
$$
p^{-1}(0)=\bigcup_{\Lambda \in J}\Lambda ,
$$
where $\Lambda $ are closed connected $H_p$-invariant sets. We
assume that $J$ is a finite connected graph with $S$ denoting the
set of vertices. We assume that the union of edges $J\setminus S$
has a natural real-\an{} structure and that every $\Lambda \in
J\setminus S$ is an \an{} torus depending \an{}ally on $\Lambda $
\wrt{} that structure.

\par We identify each edge of $J$ \an{}ally with a real \bdd{} interval and this
determines a distance on $J$ in the natural way. As in (\ref{1.22}), we assume the
continuity property:
\begin{eqnarray}
\label{1} && \hbox{For every }\Lambda _0\in J\hbox{ and every
}\epsilon >0,\ \exists\, \delta >0,\ \hbox{such that if}
\\
&&\Lambda \in J,\ {\rm dist\,}(\Lambda ,\Lambda _0)<\delta ,\hbox{
then }\Lambda \subset \{ \rho \in p^{-1}(0);\, {\rm dist\,}(\rho
,\Lambda _0)<\epsilon \} .\nonumber
\end{eqnarray}
For $\Lambda \in J\setminus S$, let $\omega (\Lambda )$ be the
corresponding rotation number, depending \an{}ally on $\Lambda $.
Assume \ekv{2} { \omega (\Lambda )\hbox{ is not identically
constant on any open edge.} } Assume for simplicity that $q$ is
real-valued. (In the general case, simply replace $q$ below by
$\Re q$.) We also assume \ekv{3} { \langle q\rangle  (\Lambda
)\hbox{ is not identically constant on any open edge.} } Here
$\langle q\rangle (\Lambda )$ denotes the average of
${q_\vert}_{\Lambda }$. As before, we also assume that $\langle{q}\rangle$
extends continuously to $J$.

\par For every $\Lambda \in J$, define the closed interval $Q_\infty
(\Lambda )$ as in Proposition 7.1, so that $Q_\infty (\Lambda )=\{
\langle q\rangle (\Lambda )\}$, when $\Lambda \in J\setminus S$
and $\omega (\Lambda )\not\in{\bf Q}$. Again, by Lemma 2.4 we know
that $\{(\Lambda ,z);\, \Lambda \in K,\, z\in Q_\infty (\Lambda )\}$ is
closed whenever $K$ is a closed subset of $J$. In particular,
$\bigcup_{\Lambda \in K}Q_\infty (\Lambda )$ is compact for every
closed subset $K\subset J$. As we have seen in the case of
surfaces of revolution in subsection 7.2, it may happen that
$$
\bigcup_{\Lambda \in \omega ^{-1}({\bf Q})\cup S}Q_\infty (\Lambda
)
$$
has a small measure compared to that of $\langle q\rangle
(J\setminus S)$.

\par We shall first define the \ufly{} good values in ${\bf R}$.  Fix
$d>0$.  Given $\alpha ,\beta ,\gamma >0$, we say that $r\in{\bf
R}$ is $(\alpha ,\beta ,\gamma )$-good if the following hold:
\smallskip\par\noindent
-- $r$ is not in the union of all $Q_\infty (\Lambda )$, with
${\rm dist\,}(\Lambda ,S)\le \alpha $ or with $\omega (\Lambda )$
not $(\alpha ,d)$-Diophantine (in the sense of (\ref{7.2})).
\smallskip
\par\noindent
-- If $r=\langle q\rangle (\Lambda )$, then $\vert d_\Lambda
\langle q\rangle (\Lambda )\vert ,\, \vert d_\Lambda \omega
(\Lambda )\vert \geq \alpha .$\smallskip

\par\noindent
-- Let $\langle q\rangle ^{-1}(r)=\{ \Lambda _1,...,\Lambda _L\}$,
with $\Lambda _j\in J\setminus S$, ${\rm dist\,}(\Lambda
_j,S)>\alpha $. Then
$${\rm dist\,}\big( r,\bigcup_{\Lambda \in J,\atop {\rm dist\,} (\Lambda ,\cup_1^L\Lambda
_j)>\beta }Q_\infty (\Lambda )\, \big) >\gamma .$$\smallskip

\par If $r$ is not $(\alpha ,\beta ,\gamma )$-good we say that it is $(\alpha ,\beta
,\gamma )$-bad. Choosing successively $\alpha ,\beta ,\gamma $
\sufly{} small, we see that the measure of the set of $(\alpha
,\beta ,\gamma )$-bad values in ${\bf R}\setminus \bigcup_{\Lambda
\in \omega ^{-1}({\bf Q})\cup S}Q_\infty (\Lambda )$ can be made
\ably{} small.
\medskip
\par\noindent \it The close to completely integrable case. \rm We now
replace $p$ by $p_\lambda $ as in (\ref{kam1}) and Theorems 7.5,
7.6. Near each $\Lambda_j $ with $\langle q\rangle (\Lambda_j)$
$(\alpha ,\beta ,\gamma )$-good, we apply the isoenergetic KAM
theorem 7.5, with $\lambda $ small enough depending only on
$\alpha ,\beta ,\gamma $ and \ufly{} \wrt{} $\Lambda $, and get a
\fy{} of KAM-tori $\Lambda _{j,\lambda ,\mu }$.

\par We say that $F_0\in {\bf R}$ is $(\alpha ,\beta ,\gamma ,\lambda
)$-good if: \smallskip

\par\noindent -- It is a $\lambda $-perturbation of an $(\alpha ,\beta
,\gamma )$-good value $r$ with $\langle q\rangle ^{-1}(r)=\{
\Lambda _1,...,\Lambda _L\}$.\smallskip

\par\noindent -- There exist KAM-tori $\Lambda _{j,\lambda }$ close to
$\Lambda _j$, that are $(\alpha /2,d)$-Diophantine, such that
$F_0=\langle q\rangle _{\Lambda _{j,\lambda }}$ for $j=1,...,L$.
\smallskip

\par Then the set of $(\alpha ,\beta ,\gamma ,\lambda )$-good values in
any \bdd{} interval has a measure that tends to that of the set of
$(\alpha ,\beta ,\gamma )$-good values in the same interval, when
$\lambda $ tends to zero. Moreover the conclusions of Theorems
1.1 and 1.2 hold \ufly{} when $\lambda $ is small enough depending on
$\alpha ,\beta ,\gamma $, and $F_0$ is $(\alpha ,\beta ,\gamma
,\lambda )$-good.

\section{Barrier top resonances: the non-resonant case}
\label{section8} \setcounter{equation}{0}
The purpose of this section is to illustrate Theorem 1.2 by applying
it to the problem of distribution of semiclassical resonances for the Schr\"odinger
operator. The discussion here will be analogous to the
corresponding treatments in~\cite{MeSj},~\cite{HiSj1},
and~\cite{HiSj2}, and for this reason, the following presentation will
be somewhat less detailed.

Consider \begeq \label{res1} P=-h^2 \Delta +V(x),\quad
p(x,\xi)=\xi^2+V(x),\quad (x,\xi)\in T^*\real^2,
\endeq
where $V$ is an analytic potential satisfying satisfying the
general assumptions of section 7 of~\cite{HiSj1}, that allow us to
define the resonances of $P$ in a fixed sector in the fourth
quadrant. Let $V(0)=E_0$, $V'(0)=0$ and assume that we are in the
barrier top case, so that $V''(0)<0$. In this section, we shall
also assume for simplicity that $V$ is an even function. The Taylor expansion
of $p$ in a suitable system of linear symplectic coordinates then takes
the form
\begeq
\label{res2}
p(x,\xi)-E_0=\sum_{j=1}^2
\frac{\lambda_j}{2}(\xi_j^2-x_j^2)+p_4(x)+\ldots,
\endeq
where $\lambda_j>0$, $j=1,2$, and $p_4$ is a homogeneous polynomial of
degree 4.

We shall assume that we are in the non-resonant case so that
\begeq \label{res3.6} \lambda\cdot k\neq 0,\quad 0\neq k\in \z^2.
\endeq
In this case, the result of~\cite{KaKer} gives all resonances of $P$ in
the open disc $D(E_0,h^{\delta})=\{z\in \comp; \abs{z-E_0}<
h^{\delta}\}$, for any $\delta>0$. They are given by
$$
z_k=E_0+f(h(k-\theta_0)h;h),\quad k\in \nat^2,
$$
where $\theta_0\in \left(\frac{1}{2}\z\right)^2$ is fixed, and
$f(\theta;h)$ is a smooth function of $\theta\in {\rm
neigh}(0,\real^2)$, with $f(\theta;h)\sim f_0(\theta)+h
f_1(\theta)+\ldots$. We have
$$
f_0(\theta)=-i\lambda\cdot \theta+{\cal O}(\theta^2).
$$
The purpose of this section is to show how Theorem 1.2 allows us
to go further away from the real axis and obtain a description of
resonances that are at a distance $\sim \eps_0$, $0<\eps_0\ll 1$,
away from the real axis.

\medskip
\noindent
As in~\cite{HiSj1},~\cite{HiSj2}, in order to study
the resonances of $P$ near $E_0$, we perform the complex
 scaling, which near $(0,0)$ is given by
$x=e^{i\pi/4}\widetilde{x}$, $\xi=e^{-i\pi/4}\widetilde{\xi}$,
$\widetilde{x}$, $\widetilde{\xi}\in \real^2$.
Using (\ref{res2}) and dropping the tildes from the notation, we get a new operator with the leading symbol
\begeq \label{res3}
\frac{1}{i}\left(p_2(x,\xi)-ip_4(x)+\ldots\right)=:\frac{1}{i}q(x,\xi),\quad
(x,\xi)\rightarrow (0,0).
\endeq
Here \begeq \label{res3.5} p_2(x,\xi)=\sum_{j=1}^2
\frac{\lambda_j}{2}(x_j^2+\xi_j^2)
\endeq
is the harmonic oscillator. We shall be interested in eigenvalues
$E$ of the operator
$$
Q(x,hD_x;h)=q(x,hD_x)+{\cal O}(h),
$$
with $\abs{E}\sim \eps^2$, $0<\eps \ll 1$. An additional
rescaling $x=\eps y$, $h^{\delta}\leq \eps \leq 1$, $0<\delta<1/2$, gives
$$
\frac{1}{\eps^2}Q(x,hD_x;h)=\frac{1}{\eps^2}Q(\eps(y,\widetilde{h}D_y);h),\quad
\widetilde{h}=\frac{h}{\eps^2}\ll 1.
$$
Viewed as an $\widetilde{h}$-pseudodifferential operator, $\eps^{-2}
Q(x,hD_x;h)$ has the leading symbol
$$
p_2(x,\xi)-i\eps^2p_4(x)+{\cal O}(\eps^4),
$$
to be considered for $(x,\xi)$ in a bounded region.

The $H_{p_2}$--flow is completely integrable, and as in~\cite{HiSj1},
we introduce the action-angle coordinates $I_j\geq 0$ and $\tau_j\in
\real/2\pi \z$, $j=1,2$, for $p_2$, given by
$$
x_j=\sqrt{2I_j}\cos \tau_j,\quad \xi_j=-\sqrt{2I_j}\sin \tau_j.
$$
We also remark that the non-resonant condition (\ref{res3.6}) implies
that the $H_{p_2}$--flow is ergodic on each invariant torus
$\Lambda_{\mu}\subset p_2^{-1}(1)$, given by
$$
I_1=\mu,\quad
I_2=\frac{1}{\lambda_2}-\frac{\lambda_1}{\lambda_2}\mu,\quad \mu
\notin \left\{0,\frac{1}{\lambda_1}\right \},
$$
Therefore, the flow average
$$
\langle{p_4}\rangle_T =\frac{1}{T} \int_0^T p_4\circ
\exp(tH_{p_2})\,dt
$$
along the torus $\Lambda_{\mu}$ converges to the torus average of
$p_4$, $\langle{p_4}\rangle=\langle{p_4}\rangle_{\infty}$, as
$T\rightarrow \infty$.

When computing the flow and torus averages, it will be convenient
to work in the symplectic coordinates $(y,\eta)$ given by
$$
y=\frac{1}{\sqrt{2}}(x-i\xi),\quad
\eta=\frac{1}{i\sqrt{2}}(x+i\xi).
$$
In these coordinates $p_2=\sum_{j=1}^2 i\lambda_j y_j \eta_j$, and the flow is given by
$$
\exp(tH_{p_2})(y,\eta)=(e^{it\lambda_1}y_1,
e^{it\lambda_2}y_2,e^{-it\lambda_1}\eta_1,e^{-it\lambda_2}\eta_2).
$$
Writing
$$
x^{\alpha}=\sum_{0\leq k\leq \alpha} a_{k\alpha} y^k
\eta^{\alpha-k},\quad \abs{\alpha}=4,\quad
a_{k\alpha}=\frac{i^{\abs{\alpha-k}}}{4} \pmatrix{\alpha \cr k},
$$
we immediately see that as $T\rightarrow \infty$, the flow average
$\langle{x^{\alpha}}\rangle_T$ converges to
$$
\langle{x^{\alpha}}\rangle=\sum_{{0\leq k\leq\alpha} \atop{
2k=\alpha}} a_{k\alpha} y^k \eta^{\alpha-k}.
$$
It follows that if
$$
p_4(x)=\sum_{\abs{\alpha}=4} v_{\alpha} x^{\alpha},
$$
then
$$
\langle{p_4}\rangle(y,\eta)=\frac{-1}{4}\left(6
v_{(4,0)}y_1^2\eta_1^2+4 v_{(2,2)} y_1 y_2 \eta_1\eta_2+6
v_{(0,4)}y_2^2 \eta_2^2\right).
$$
This function is in involution with $p_2$, and can be expressed in terms of the action variables as
follows,
$$
\langle{p_4}\rangle=\frac{1}{4}\left(6v_{(4,0)}I_1^2+4v_{(2,2)}I_1
I_2+6v_{(0,4)}I_2^2\right).
$$

Now $p_2=\lambda_1 I_1+\lambda_2 I_2$, and therefore we see that for
$I_1 I_2\neq 0$, the differentials $dp_2=(\lambda_1,\lambda_2)$ and
$d\langle{p_4}\rangle=-(3v_{(4,0)}I_1+v_{(2,2)}I_2,3v_{(0,4)}I_2+v_{(2,2)}I_1)$
are linearly dependent precisely when
\begeq
\label{res4}
(\lambda_1 v_{(2,2)}-3\lambda_2 v_{(4,0)})I_1=(\lambda_2
v_{(2,2)}-3\lambda_1 v_{(0,4)})I_2.
\endeq
When $I_1=0$ or $I_2=0$, the question of linear independence of the
differentials should be examined directly in the
$(y,\eta)$--coordinates, and we then easily see that we have the linear
dependence when $I_1 I_2=0$. We also compute that $\langle{p_4}\rangle$ restricted to
$p_2^{-1}(1)$ has the critical values $A_1=\frac{3}{2}
v_{(0,4)}\lambda_2^{-2}$, corresponding to $I_1=0$,
$A_2=\frac{3}{2} v_{(4,0)}\lambda_1^{-2}$, corresponding to
$I_2=0$, and a third value $A_3$ which occurs when the line
(\ref{res4}) intersects the quadrant $I_1>0$, $I_2>0$, i.e., when
$$
(\lambda_1 v_{(2,2)}-3\lambda_2 v_{(2,2)})(\lambda_2
v_{(2,2)}-3\lambda_1 v_{(0,4)})>0. $$

We are now in a position to apply Theorem 1.2 to the operator
$$
\frac{1}{\eps^2} Q(x,hD_x)-1,
$$
with $\eps^2$ considered as a small perturbation parameter and using
$\widetilde{h}=\frac{h}{\eps^2}$ as the new semiclassical parameter. Here we should
also remark that the Birkhoff normal form construction of section 3
goes through in the present case assuming only (\ref{res3.6}) and that
the Diophantine condition is not required---see also~\cite{Sj92} for the
classical Birkhoff construction in the non-resonant case, near a
stable equilibrium point.

\begin{theo}
Consider the operator $P$ in {\rm (\ref{res1})} with the leading
symbol
$$
p(x,\xi)-E_0=\sum_{j=1}^2
\frac{\lambda_j}{2}\left(\xi_j^2-x_j^2\right)+x_1^4+{\cal
O}(x^6),\quad (x,\xi)\rightarrow (0,0).
$$
Assume that the non-resonance condition
$\langle{\lambda,k}\rangle\neq 0$, $0\neq k\in \z^2$, holds. Let us
put $A_1=0$, $A_2=\frac{3}{2}\lambda_1^{-2}$. Assume that
$0<\delta<\frac{1}{4}$. Then the resonances of $P$ in the domain
\begeq \label{8.6}
\{z\in \comp; h^{\delta}\ll \abs{z-E_0}\leq \eps_0 \}\backslash
\bigcup_{j=1}^2 \{z\in \comp; \abs{\Re
z-E_0-A_j\abs{\Im z}^2}<\eta \abs{\Im z}^2\},
\endeq
where $0<\eps_0=\eps_0(\eta)\ll 1$ and $\eta>0$ is arbitrary small but
fixed, are given by
$$
\sim E_0-i\sum_{n=0}^{\infty} h^n \eps^{2-2n}
\widetilde{p}^{(\infty)}_n\left(\frac{h}{\eps^2}\left(k-\frac{\alpha}{4}\right)-\frac{S}{2\pi},\eps\right),\quad
k\in \z^2,
$$
with $\widetilde{p}^{(\infty)}_n(\xi,\eps)$ analytic in $\xi\in {\rm neigh}(0,\comp^2)$ and
smooth in $\eps\in {\rm neigh}(0,\real)$. We have $S\in \real^2$ and
$\alpha\in \z^2$ are fixed, and we choose
$\eps>0$ with $\abs{E-E_0}\sim \eps^2$.
\end{theo}

\Remark. Combining Theorem 8.1 with the arguments of section 7
of~\cite{MeSj}, we obtain that
the result of~\cite{KaKer} extends to a set of the
form
$$
D(E_0,\eps_0)\backslash \bigcup_{j=1}^{2} \left\{z; \Re z-E_0\in
-[A_j-\eta, A_j+\eta](\Im z)^2\right\}.
$$
Here $D(E_0,\eps_0)=\{z, \abs{z-E_0}<\eps_0\}$, with
$0<\eps_0=\eps_0(\eta)\ll 1$,
and $\eta>0$ is arbitrary but fixed.

\end{document}